\documentclass[12pt]{article}%
\usepackage{amsmath}
\usepackage{amsfonts}
\usepackage{amssymb}
\usepackage{graphicx}%
\providecommand{\U}[1]{\protect\rule{.1in}{.1in}}
\begin{document}

\title{Vector Bundles over Multipullback Quantum Complex Projective
Spaces\thanks{This work was partially supported by the grant
H2020-MSCA-RISE-2015-691246-QUANTUM DYNAMICS and the Polish government grant
3542/H2020/2016/2.}}
\date{}
\author{Albert Jeu-Liang Sheu\thanks{The author would like to thank the Mathematics
Institute of Academia Sinica for the warm hospitality and support during his
visit in the summer of 2017.}\\{\small Department of Mathematics, University of Kansas, Lawrence, KS 66045,
U. S. A.}\\{\small e-mail: asheu@ku.edu}}
\maketitle

\begin{abstract}
We work on the classification of isomorphism classes of finitely generated
projective modules over the C*-algebras $C\left(  \mathbb{P}^{n}\left(
\mathcal{T}\right)  \right)  $ and $C\left(  \mathbb{S}_{H}^{2n+1}\right)  $
of the quantum complex projective spaces $\mathbb{P}^{n}\left(  \mathcal{T}%
\right)  $ and the quantum spheres $\mathbb{S}_{H}^{2n+1}$, and the quantum
line bundles $L_{k}$ over $\mathbb{P}^{n}\left(  \mathcal{T}\right)  $,
studied by Hajac and collaborators. Motivated by the groupoid approach of
Curto, Muhly, and Renault to the study of C*-algebraic structure, we analyze
$C\left(  \mathbb{P}^{n}\left(  \mathcal{T}\right)  \right)  $, $C\left(
\mathbb{S}_{H}^{2n+1}\right)  $, and $L_{k}$ in the context of groupoid
C*-algebras, and then apply Rieffel's stable rank results to show that all
finitely generated projective modules over $C\left(  \mathbb{S}_{H}%
^{2n+1}\right)  $ of rank higher than $\left\lfloor \frac{n}{2}\right\rfloor
+3$ are free modules. Furthermore, besides identifying a large portion of the
positive cone of the $K_{0}$-group of $C\left(  \mathbb{P}^{n}\left(
\mathcal{T}\right)  \right)  $, we also explicitly identify $L_{k}$ with
concrete representative elementary projections over $C\left(  \mathbb{P}%
^{n}\left(  \mathcal{T}\right)  \right)  $.

\begin{description}
\item[Keywords:] multipullback quantum projective space; multipullback quantum
sphere; quantum line bundle; finitely generated projective module;
cancellation problem; Toeplitz algebra of polydisk; groupoid C*-algebra;
stable rank; noncommutative vector bundle

\item[AMS 2010 Mathematics Subject Classification:] 46L80; 46L85

\end{description}
\end{abstract}

\pagebreak

\section{Introduction}

Since the concept of noncommutative geometry first popularized by Connes
\cite{Conn}, many interesting examples of a C*-algebra $\mathcal{A}$ viewed as
the algebra $C\left(  X_{q}\right)  $ of continuous functions on a virtual
quantum space $X_{q}$ have been constructed with a topological or geometrical
motivation, and analyzed in comparison with their classical counterpart. For
example, quantum odd-dimensional spheres and associated complex projective
spaces have been introduced and studied by Soibelman, Vaksman, Meyer, and
others \cite{VaSo,Me} as $\mathbb{S}_{q}^{2n+1}$ and $\mathbb{C}P_{q}^{n}$ via
a quantum universal enveloping algebra approach, and by Hajac and his
collaborators including Baum, Kaygun, Matthes, Nest, Pask, Sims,
Szyma\'{n}ski, Zieli\'{n}ski, and others
\cite{BaHaMaSz,HaMaSz,HaKaZi,HaNePaSiZi} as $\mathbb{S}_{H}^{2n+1}$ and
$\mathbb{P}^{n}\left(  \mathcal{T}\right)  $ via a multi-pullback and Toeplitz
algebra approach. Actually $\mathbb{S}_{H}^{2n+1}$ is the untwisted special
case of the more general version of $\theta$-twisted spheres $\mathbb{S}%
_{H,\theta}^{2n+1}$ introduced in \cite{HaNePaSiZi}.

Motivated by Swan's work \cite{Swan}, the concept of a noncommutative vector
bundle $E_{q}$ over a quantum space $X_{q}$ can be reformulated as a finitely
generated projective (left) module $\Gamma\left(  E_{q}\right)  $ over
$C\left(  X_{q}\right)  $. Based on the strong connection approach to quantum
principal bundles \cite{Haja:sc} for compact quantum groups \cite{Woro,Wo:cm},
Hajac and his collaborators introduced quantum line bundles $L_{k}$ of degree
$k$ over $\mathbb{P}^{n}\left(  \mathcal{T}\right)  $ as some rank-one
projective modules realized as spectral subspaces $C\left(  \mathbb{S}%
_{H}^{2n+1}\right)  _{k}$ of $C\left(  \mathbb{S}_{H}^{2n+1}\right)  $ under a
$U\left(  1\right)  $-action \cite{HaNePaSiZi}. Besides having the $K_{0}%
$-group of $C\left(  \mathbb{P}^{n}\left(  \mathcal{T}\right)  \right)  $
computed, they found that $L_{k}$ is not stably free unless $k=0$, extending
earlier results for the case of $n=1$ \cite{HaMaSz,HaMaSz2006}.

It has always been an interesting but challenging task to classify finitely
generated projective modules over an algebra up to isomorphism, which goes
beyond their classification up to stable isomorphism by $K_{0}$-group and
appears in the form of so-called cancellation problem. Classically it is known
that the cancellation law holds for complex vector bundles of rank no less
than $\frac{d}{2}$ over a $d$-dimensional CW-complex, which implies that all
complex vector bundles over $\mathbb{S}^{2n+1}$ of rank $n+1$ or above are trivial.

The study of such classification problem for C*-algebras was popularized by
Rieffel \cite{Ri:dsr,Ri:ct} who introduced useful versions of stable ranks for
C*-algebras to facilitate the analysis involved. Some successes have been
achieved for certain quantum algebras \cite{Ri:ct,Ri:pm,Sh:ct,Bach,Pete}. In
particular, Peterka showed that all finitely generated projective modules over
the $\theta$-deformed $3$-spheres $S_{\theta}^{3}$ are free, and constructed
all those over $S_{\theta}^{4}$ up to isomorphism \cite{Pete}. With more
effort, the result of Bach \cite{Bach} on the cancellation law for
$\mathbb{S}_{q}^{2n+1}$ and $\mathbb{C}P_{q}^{n}$ can be strengthened to a
complete classification of finitely generated projective modules over them,
which we will address elsewhere.

With the $K_{0}$-group of $C\left(  \mathbb{P}^{n}\left(  \mathcal{T}\right)
\right)  $ known \cite{HaNePaSiZi}, it is natural to try to classify finitely
generated projective modules over $C\left(  \mathbb{P}^{n}\left(
\mathcal{T}\right)  \right)  $ and identify the line bundles $L_{k}$ among
them. In \cite{Sh:pmqpl}, a complete solution was obtained for the special
case of $n=1$.

In this paper, we use the powerful groupoid approach to C*-algebras initiated
by Renault \cite{Rena} and popularized by Curto, Muhly, and Renault
\cite{CuMu,MuRe} to study multi-variable Toeplitz C*-algebras $\mathcal{T}%
^{\otimes n}$, quantum spheres $C\left(  \mathbb{S}_{H}^{2n+1}\right)  $, and
quantum complex projective spaces $C\left(  \mathbb{P}^{n}\left(
\mathcal{T}\right)  \right)  $. Utilizing results on stable ranks of
C*-algebras obtained by Rieffel \cite{Ri:dsr}, we analyze finitely generated
projective modules over $\mathcal{T}^{\otimes n+1}$ and $C\left(
\mathbb{S}_{H}^{2n+1}\right)  $, and get those of rank higher than
$\left\lfloor \frac{n}{2}\right\rfloor +3$ and also a large class of
\textquotedblleft standard\textquotedblright\ modules classified up to
isomorphism. Furthermore, besides identifying a large portion of the positive
cone of the $K_{0}$-group $K_{0}\left(  C\left(  \mathbb{P}^{n}\left(
\mathcal{T}\right)  \right)  \right)  $, we explicitly identify the quantum
line bundles $L_{k}$ with concrete representative elementary projections.

On the other hand, there are still a lot of questions to be further
investigated, e.g. whether the cancellation law holds for low-ranked finitely
generated projective modules, and whether the more general case of $\theta
$-twisted multipullback quantum sphere $\mathbb{S}_{H,\theta}^{2n+1}$ brings
in new phenomena. Finally it is of interest to note the recent work of Farsi,
Hajac, Maszczyk, and Zieli\'{n}ski \cite{FHMZ} on $K_{0}\left(  C\left(
\mathbb{P}^{2}\left(  \mathcal{T}\right)  \right)  \right)  $, identifying its
free generators arising from Milnor modules as sums of $L_{k} $, which are
also expressed in terms of elementary projections, showing a perfect
consistency with our result.

\section{Notations}

Taking the groupoid approach to C*-algebras initiated by Renault \cite{Rena}
and popularized by the work of Curto, Muhly, and Renault \cite{CuMu,MuRe}, we
give a description of the C*-algebras $C\left(  \mathbb{S}_{H}^{2n-1}\right)
$ and $C\left(  \mathbb{P}^{n-1}\left(  \mathcal{T}\right)  \right)  $ of
\cite{HaNePaSiZi} as some concrete groupoid C*-algebras. We refer to
\cite{Rena,MuRe} for the concepts and theory of groupoid C*-algebras used
freely in the following discussion.

By abuse of notation, for any C*-algebra homomorphism $\phi:\mathcal{A}%
\rightarrow\mathcal{B}$, we denote the C*-algebra homomorphism $M_{k}\left(
\phi\right)  :M_{k}\left(  \mathcal{A}\right)  \rightarrow M_{k}\left(
\mathcal{B}\right)  $ for $k\in\mathbb{N}\equiv\left\{  1,2,3,...\right\}  $
also by $\phi$. We use $\mathcal{A}^{\times}$ to denote the set of all
invertible elements of an algebra $\mathcal{A}$, and use $\mathcal{A}^{+}$ to
denote the minimal unitization of $\mathcal{A}$. For any topological group
$G$, we use $G^{0}$ to denote the identity component of $G$, i.e. the
connected component that contains the identity element of $G$.

We denote by $M_{\infty}\left(  \mathcal{A}\right)  $ the direct limit (or the
union as sets) of the increasing sequence of matrix algebras $M_{n}\left(
\mathcal{A}\right)  $ over $\mathcal{A}$ with the canonical inclusion
$M_{n}\left(  \mathcal{A}\right)  \subset M_{n+1}\left(  \mathcal{A}\right)  $
identifying $x\in M_{n}\left(  \mathcal{A}\right)  $ with $x\boxplus0\in
M_{n+1}\left(  \mathcal{A}\right)  $ for any algebra $\mathcal{A}$, where
$\boxplus$ denotes the standard diagonal concatenation (sum) of two matrices.
So the size of an element in $M_{\infty}\left(  \mathcal{A}\right)  $ can be
taken arbitrarily large. We also use $GL_{\infty}\left(  \mathcal{A}\right)  $
to denote the direct limit of the general linear groups $GL_{n}\left(
\mathcal{A}\right)  $ over a unital C*-algebra $\mathcal{A}$ with
$GL_{n}\left(  \mathcal{A}\right)  $ embedded in $GL_{n+1}\left(
\mathcal{A}\right)  $ by identifying $x\in GL_{n}\left(  \mathcal{A}\right)  $
with $x\boxplus1\in GL_{n+1}\left(  \mathcal{A}\right)  $.

By an idempotent $P$ over a unital C*-algebra $\mathcal{A}$, we mean an
element $P\in M_{\infty}\left(  \mathcal{A}\right)  $ with $P^{2}=P$, and a
self-adjoint idempotent in $M_{\infty}\left(  \mathcal{A}\right)  $ is called
a projection over $\mathcal{A}$. Two idempotents $P,Q\in M_{\infty}\left(
\mathcal{A}\right)  $ are called equivalent, denoted as $P\sim Q$, if there
exists $U\in GL_{\infty}\left(  \mathcal{A}\right)  $ such that $UPU^{-1}=Q$.
Each idempotent $P\in M_{n}\left(  \mathcal{A}\right)  $ over $\mathcal{A}$
defines a finitely generated left projective module $E:=\mathcal{A}^{n}P$ over
$\mathcal{A}$ where elements of $\mathcal{A}^{n}$ are viewed as row vectors.
The mapping $P\mapsto\mathcal{A}^{n}P$ induces a bijective correspondence
between the equivalence classes of idempotents over $\mathcal{A}$ and the
isomorphism classes of finitely generated left projective modules over
$\mathcal{A}$ \cite{Blac}. From now on, by a module over $\mathcal{A}$, we
mean a left $\mathcal{A}$-module, unless otherwise specified.

Two finitely generated projective modules $E,F$ over $\mathcal{A}$ are called
stably isomorphic if they become isomorphic after being augmented by the same
finitely generated free $\mathcal{A}$-module, i.e. $E\oplus\mathcal{A}%
^{k}\cong F\oplus\mathcal{A}^{k}$ for some $k\geq0$. Correspondingly, two
idempotents $P$ and $Q$ are called stably equivalent if $P\boxplus I_{k}$ and
$Q\boxplus I_{k}$ are equivalent for some identity matrix $I_{k}$. The $K_{0}%
$-group $K_{0}\left(  \mathcal{A}\right)  $ classifies idempotents over
$\mathcal{A}$ up to stable equivalence. The classification of idempotents over
a C*-algebra up to equivalence, appearing as the so-called cancellation
problem, was popularized by Rieffel's pioneering work \cite{Ri:dsr,Ri:ct} and
is in general an interesting but difficult question.

The set of all equivalence classes of idempotents over a C*-algebra
$\mathcal{A}$ is an abelian monoid $\mathfrak{P}\left(  \mathcal{A}\right)  $
with its binary operation provided by the diagonal sum $\boxplus$. The image
of the canonical homomorphism from $\mathfrak{P}\left(  \mathcal{A}\right)  $
into $K_{0}\left(  \mathcal{A}\right)  $ is the so-called positive cone of
$K_{0}\left(  \mathcal{A}\right)  $.

Furthermore, it is well-known \cite{Blac} that in the above descriptions of
$\mathfrak{P}\left(  \mathcal{A}\right)  $ and $K_{0}\left(  \mathcal{A}%
\right)  $, one can restrict to the self-adjoint idempotents, called
projections over $\mathcal{A}$, and their unitary equivalence classes, which
faithfully represent the elements of $\mathfrak{P}\left(  \mathcal{A}\right)
$ and $K_{0}\left(  \mathcal{A}\right)  $.

In this paper, we use freely the basic techniques and manipulations for
$K$-theory found in \cite{Blac,Tay}.

For a Hilbert space $\mathcal{H}$, we denote the C*-algebra consisting of all
compact linear operators on $\mathcal{H}$ by $\mathcal{K}\left(
\mathcal{H}\right)  $, or simply by $\mathcal{K}$ if $\mathcal{H}$ is the
essentially unique separable infinite-dimensional Hilbert space.

In the following, we use the notations $\mathbb{Z}_{\geq k}:=\left\{
n\in\mathbb{Z}|n\geq k\right\}  $ and $\mathbb{Z}_{\geq}:=\mathbb{Z}_{\geq0}$.
In particular, $\mathbb{N}=\mathbb{Z}_{\geq1}$. We use $I$ to denote the
identity operator canonically contained in $\mathcal{K}^{+}\subset
\mathcal{B}\left(  \ell^{2}\left(  \mathbb{Z}_{\geq}\right)  \right)  $, and
\[
P_{m}:=\sum_{i=1}^{m}e_{ii}\in M_{m}\left(  \mathbb{C}\right)  \subset
\mathcal{K}%
\]
to denote the standard $m\times m$ identity matrix in $M_{m}\left(
\mathbb{C}\right)  \subset\mathcal{K}$ for any integer $m\geq0$ (with
$M_{0}\left(  \mathbb{C}\right)  =0$ and $P_{0}=0$ understood). We also use
the notation
\[
P_{-m}:=I-P_{m}\in\mathcal{K}^{+}%
\]
for integers $m>0$, and take symbolically $P_{-0}\equiv I-P_{0}=I\neq P_{0}$.

\section{Quantum spaces as groupoid C*-algebras}

Let $\mathfrak{T}_{n}:=\left.  \left(  \mathbb{Z}^{n}\ltimes\overline
{\mathbb{Z}}^{n}\right)  \right\vert _{\overline{\mathbb{Z}_{\geq}}^{n}}$ with
$n\geq1$ be the transformation group groupoid $\mathbb{Z}^{n}\ltimes
\overline{\mathbb{Z}}^{n}$ restricted to the positive \textquotedblleft
cone\textquotedblright\ $\overline{\mathbb{Z}_{\geq}}^{n}$ where
$\overline{\mathbb{Z}}:=\mathbb{Z}\cup\left\{  +\infty\right\}  $ containing
$\mathbb{Z}_{\geq}\equiv\left\{  n\in\mathbb{Z}|n\geq0\right\}  $ carries the
standard topology, and $\mathbb{Z}^{n}$ acts on $\overline{\mathbb{Z}}^{n}$
componentwise in the canonical way. From the groupoid isomorphism
\[
\left.  \left(  \mathbb{Z}^{n}\ltimes\overline{\mathbb{Z}}^{n}\right)
\right\vert _{\overline{\mathbb{Z}_{\geq}}^{n}}\cong\times^{n}\left(  \left.
\left(  \mathbb{Z}\ltimes\overline{\mathbb{Z}}\right)  \right\vert
_{\overline{\mathbb{Z}_{\geq}}}\right)
\]
and the well-known C*-algebra isomorphism $C^{\ast}\left(  \left.  \left(
\mathbb{Z}\ltimes\overline{\mathbb{Z}}\right)  \right\vert _{\overline
{\mathbb{Z}_{\geq}}}\right)  \cong\mathcal{T}$ for the Toeplitz C*-algebra
$\mathcal{T}$, we get the groupoid C*-algebra
\[
C^{\ast}\left(  \mathfrak{T}_{n}\right)  \equiv C^{\ast}\left(  \left.
\left(  \mathbb{Z}^{n}\ltimes\overline{\mathbb{Z}}^{n}\right)  \right\vert
_{\overline{\mathbb{Z}_{\geq}}^{n}}\right)  \cong\mathcal{T}^{\otimes n}%
\equiv\otimes^{n}\mathcal{T}.
\]

We consider two important nontrivial invariant open subsets of the unit space
$\overline{\mathbb{Z}_{\geq}}^{n}$ of $\mathfrak{T}_{n}$, namely,
$\mathbb{Z}_{\geq}^{n}$ the smallest one and $\overline{\mathbb{Z}_{\geq}}%
^{n}\backslash\left\{  \infty^{n}\right\}  $ the largest one, where
$\infty^{n}:=\left(  \infty,...,\infty\right)  \in\overline{\mathbb{Z}_{\geq}%
}^{n} $. By the theory of groupoid C*-algebras developed in Renault's book
\cite{Rena}, they give rise to two short exact sequences of C*-algebras%
\[
0\rightarrow C^{\ast}\left(  \left.  \left(  \mathbb{Z}^{n}\ltimes
\overline{\mathbb{Z}}^{n}\right)  \right\vert _{\mathbb{Z}_{\geq}^{n}}\right)
\cong\mathcal{K}\left(  \ell^{2}\left(  \mathbb{Z}_{\geq}^{n}\right)  \right)
\rightarrow C^{\ast}\left(  \mathfrak{T}_{n}\right)  \equiv\mathcal{T}%
^{\otimes n}\rightarrow C^{\ast}\left(  \mathfrak{G}_{n}\right)  \rightarrow0
\]
with $\mathcal{K}\left(  \ell^{2}\left(  \mathbb{Z}_{\geq}^{n}\right)
\right)  \cong\otimes^{n}\mathcal{K}\equiv\mathcal{K}^{\otimes n}$ where
\[
\mathfrak{G}_{n}:=\left.  \left(  \mathbb{Z}^{n}\ltimes\overline{\mathbb{Z}%
}^{n}\right)  \right\vert _{\overline{\mathbb{Z}_{\geq}}^{n}\backslash
\mathbb{Z}_{\geq}^{n}}%
\]
is $\mathfrak{T}_{n}$ restricted to the \textquotedblleft limit
boundary\textquotedblright\ of its unit space, and
\[
0\rightarrow C^{\ast}\left(  \left.  \left(  \mathbb{Z}^{n}\ltimes
\overline{\mathbb{Z}}^{n}\right)  \right\vert _{\overline{\mathbb{Z}_{\geq}%
}^{n}\backslash\left\{  \infty^{n}\right\}  }\right)  \rightarrow C^{\ast
}\left(  \mathfrak{T}_{n}\right)  \equiv\mathcal{T}^{\otimes n}\overset{\sigma
_{n}}{\rightarrow}C^{\ast}\left(  \left.  \left(  \mathbb{Z}^{n}%
\ltimes\overline{\mathbb{Z}}^{n}\right)  \right\vert _{\left\{  \infty
^{n}\right\}  }\right)  \cong C^{\ast}\left(  \mathbb{Z}^{n}\right)  \cong
C\left(  \mathbb{T}^{n}\right)  \rightarrow0
\]
where the quotient map $\sigma_{n}$ extends the notion of the well-known
symbol map $\sigma$ on $\mathcal{T}$ in the case of $n=1$.

Note that the open invariant set $\mathbb{Z}_{\geq}^{n}$ being dense in the
unit space $\overline{\mathbb{Z}_{\geq}}^{n}$ of $\mathfrak{T}_{n}$ induces a
faithful representation $\pi_{n}$ of $C^{\ast}\left(  \mathfrak{T}_{n}\right)
$ on $\ell^{2}\left(  \mathbb{Z}_{\geq}^{n}\right)  $ that realizes the
groupoid C*-algebra $C^{\ast}\left(  \mathfrak{T}_{n}\right)  $ and its closed
ideal $C^{\ast}\left(  \left.  \left(  \mathbb{Z}^{n}\ltimes\overline
{\mathbb{Z}}^{n}\right)  \right\vert _{\mathbb{Z}_{\geq}^{n}}\right)  $
respectively as a C*-subalgebra of $\mathcal{B}\left(  \ell^{2}\left(
\mathbb{Z}_{\geq}^{n}\right)  \right)  $ and the closed ideal $\mathcal{K}%
\left(  \ell^{2}\left(  \mathbb{Z}_{\geq}^{n}\right)  \right)  $ consisting of
all compact operators on $\ell^{2}\left(  \mathbb{Z}_{\geq}^{n}\right)  $.

In this paper, we freely identify elements of $C^{\ast}\left(  \mathfrak{T}%
_{n}\right)  \equiv\mathcal{T}^{\otimes n}$ with operators on $\ell^{2}\left(
\mathbb{Z}_{\geq}^{n}\right)  $ via the faithful representation $\pi_{n}$ and
use these two conceptually different notions interchangeably.

In \cite{HaNePaSiZi}, Hajac, Nest, Pask, Sims, and Zieli\'{n}ski defined the
(untwisted) \emph{multipullback} or \emph{Heegaard }quantum odd-dimensional
sphere $S_{H}^{2n-1}$ as the quantum space of the multipullback C*-algebra
\cite{Pe} determined by homomorphisms of the form $\mathrm{id}^{\otimes
j}\otimes\sigma\otimes\mathrm{id}^{\otimes n-j-1}$ from $\mathcal{T}^{\otimes
i}\otimes C\left(  \mathbb{T}\right)  \otimes\mathcal{T}^{\otimes n-i-1}$ with
$i\neq j$ to some $\mathcal{T}^{\otimes m}\otimes C\left(  \mathbb{T}\right)
\otimes\mathcal{T}^{\otimes k}\otimes C\left(  \mathbb{T}\right)
\otimes\mathcal{T}^{\otimes n-m-k-2}$. (Actually more general $\theta$-twisted
quantum spheres $S_{H,\theta}^{2n-1}$ are studied there.) They showed that
\[
C\left(  \mathbb{S}_{H}^{2n-1}\right)  \cong\left(  \otimes^{n}\mathcal{T}%
\right)  /\left(  \otimes^{n}\mathcal{K}\right)
\]
and hence we have
\[
C\left(  \mathbb{S}_{H}^{2n-1}\right)  \cong C^{\ast}\left(  \mathfrak{G}%
_{n}\right)
\]
identified as a groupoid C*-algebra.

With the ideal $C^{\ast}\left(  \left.  \left(  \mathbb{Z}^{n}\ltimes
\overline{\mathbb{Z}}^{n}\right)  \right\vert _{\overline{\mathbb{Z}_{\geq}%
}^{n}\backslash\left\{  \infty^{n}\right\}  }\right)  $ containing the ideal
$C^{\ast}\left(  \left.  \left(  \mathbb{Z}^{n}\ltimes\overline{\mathbb{Z}%
}^{n}\right)  \right\vert _{\mathbb{Z}_{\geq}^{n}}\right)  $, the quotient map
$\sigma_{n}$ induces a well-defined quotient map $\tau_{n}$ in the short exact
sequence%
\[
0\rightarrow C^{\ast}\left(  \left.  \left(  \mathbb{Z}^{n}\ltimes
\overline{\mathbb{Z}}^{n}\right)  \right\vert _{\overline{\mathbb{Z}_{\geq}%
}^{n}\backslash\left(  \mathbb{Z}_{\geq}^{n}\cup\left\{  \infty^{n}\right\}
\right)  }\right)  \rightarrow C\left(  \mathbb{S}_{H}^{2n-1}\right)
\cong\left(  \otimes^{n}\mathcal{T}\right)  /\left(  \otimes^{n}%
\mathcal{K}\right)  \overset{\tau_{n}}{\rightarrow}C\left(  \mathbb{T}%
^{n}\right)  \rightarrow0.
\]

\section{Stable ranks of quantum spaces}

In his seminal paper \cite{Ri:dsr}, Rieffel introduced and popularized the
notions of topological stable rank $\operatorname{tsr}\left(  \mathcal{A}%
\right)  $ and connected stable rank $\operatorname{csr}\left(  \mathcal{A}%
\right)  $ of a C*-algebra $\mathcal{A}$, which are useful tools in the study
of cancellation problems for finitely generated projective modules. Later,
Herman and Vaserstein \cite{HeVa} showed that for C*-algebras $\mathcal{A}$,
Rieffel's topological stable rank coincides with the Bass stable rank used in
algebraic K-theory. So we will denote $\operatorname{tsr}\left(
\mathcal{A}\right)  $ simply as $\operatorname{sr}\left(  \mathcal{A}\right)
$ in our discussion.

In this section, we review an estimate of the stable ranks of the Toeplitz
algebras $\mathcal{T}^{\otimes n}$ and quantum spheres $C\left(
\mathbb{S}_{H}^{2n-1}\right)  $, which will be used in our study of their
finitely generated projective modules. For the case of $n=1$, it is known
\cite{Ri:dsr} that $\operatorname{sr}\left(  \mathcal{T}\right)
=\operatorname{csr}\left(  C\left(  \mathbb{T}\right)  \right)  =2$.

As an illustration of the groupoid approach to C*-algebras, we first establish
some composition sequence structure for $\mathcal{T}^{\otimes n}$ and
$C\left(  \mathbb{S}_{H}^{2n-1}\right)  $, which leads to an easy estimate of
their stable ranks.

\textbf{Proposition 1}. There is a finite composition sequence of closed
ideals
\[
\mathcal{T}^{\otimes n}\equiv C^{\ast}\left(  \mathfrak{T}_{n}\right)
\equiv\mathcal{I}_{n}\rhd\mathcal{I}_{n-1}\rhd\cdots\rhd\mathcal{I}_{1}%
\rhd\mathcal{I}_{0}\rhd\mathcal{I}_{-1}\equiv\left\{  0\right\}
\]
such that $\mathcal{T}^{\otimes n}/\mathcal{I}_{0}\cong C\left(
\mathbb{S}_{H}^{2n-1}\right)  $, and for $0\leq j\leq n$,
\[
\mathcal{I}_{j}/\mathcal{I}_{j-1}\cong\oplus^{\frac{n!}{j!\left(  n-j\right)
!}}\left(  \mathcal{K}\left(  \ell^{2}\left(  \mathbb{Z}_{\geq}^{n-j}\right)
\right)  \otimes C\left(  \mathbb{T}^{j}\right)  \right)  ,
\]
where $\mathbb{T}^{0}$ and $\mathbb{Z}_{\geq}^{0}$ denote a singleton.

Proof. For $0\leq j\leq n$, let $X_{j}$ be the set consisting of
$z\in\overline{\mathbb{Z}_{\geq}}^{n}$ with exactly $j$ of the components
$z_{1},z_{2},...,z_{n}$ being equal to $\infty$, and hence $X_{n}=\left\{
\infty^{n}\right\}  $. Then the sets%
\[
Y_{j}:=X_{0}\sqcup X_{1}\sqcup\cdots\sqcup X_{j}%
\]
are open invariant subsets of the unit space $\overline{\mathbb{Z}_{\geq}}%
^{n}$ of $\mathfrak{T}_{n}$ with
\[
\mathbb{Z}_{\geq}^{n}=Y_{0}\subset Y_{1}\subset\cdots\subset Y_{n}%
=\overline{\mathbb{Z}_{\geq}}^{n}%
\]
which determines an increasing chain of closed ideals $\mathcal{I}%
_{0}\vartriangleleft\mathcal{I}_{1}\vartriangleleft\cdots\vartriangleleft
\mathcal{I}_{n}$ of $C^{\ast}\left(  \mathfrak{T}_{n}\right)  $ defined by
\[
\mathcal{I}_{j}:=C^{\ast}\left(  \left.  \left(  \mathbb{Z}^{n}\ltimes
\overline{\mathbb{Z}}^{n}\right)  \right\vert _{Y_{j}}\right)  \equiv C^{\ast
}\left(  \left.  \mathfrak{T}_{n}\right\vert _{Y_{j}}\right)  .
\]

Note that $Y_{j}\backslash Y_{j-1}=X_{j}$ with $Y_{-1}:=\emptyset$ is a
disjoint union of $\frac{n!}{j!\left(  n-j\right)  !}$ copies of
$\mathbb{Z}_{\geq}^{n-j}\times\left\{  \infty^{j}\right\}  $ each of which is
gotten from one of the $\frac{n!}{j!\left(  n-j\right)  !}$ possible
selections of exactly $j$ of the $n$ components of $\overline{\mathbb{Z}%
_{\geq}}^{n}$. With each such copy of $\mathbb{Z}_{\geq}^{n-j}\times\left\{
\infty^{j}\right\}  $ clearly a closed invariant subset of $Y_{j}\backslash
Y_{j-1}$, these $\frac{n!}{j!\left(  n-j\right)  !}$ copies of $\mathbb{Z}%
_{\geq}^{n-j}\times\left\{  \infty^{j}\right\}  $ are open invariant subsets
of $Y_{j}\backslash Y_{j-1}$, and hence
\[
C^{\ast}\left(  \left.  \mathfrak{T}_{n}\right\vert _{Y_{j}\backslash Y_{j-1}%
}\right)  =\oplus^{\frac{n!}{j!\left(  n-j\right)  !}}C^{\ast}\left(  \left.
\left(  \mathbb{Z}^{n}\ltimes\overline{\mathbb{Z}}^{n}\right)  \right\vert
_{\mathbb{Z}_{\geq}^{n-j}\times\left\{  \infty^{j}\right\}  }\right)
\]%
\[
=\oplus^{\frac{n!}{j!\left(  n-j\right)  !}}C^{\ast}\left(  \left(  \left.
\left(  \mathbb{Z}^{n-j}\ltimes\mathbb{Z}^{n-j}\right)  \right\vert
_{\mathbb{Z}_{\geq}^{n-j}}\right)  \times\mathbb{Z}^{j}\right)  =\oplus
^{\frac{n!}{j!\left(  n-j\right)  !}}\left(  \mathcal{K}\left(  \ell
^{2}\left(  \mathbb{Z}_{\geq}^{n-j}\right)  \right)  \otimes C\left(
\mathbb{T}^{j}\right)  \right)  .
\]
Thus with $\mathcal{I}_{j}=C^{\ast}\left(  \left.  \mathfrak{T}_{n}\right\vert
_{Y_{j}}\right)  $ and $\mathcal{I}_{j-1}=C^{\ast}\left(  \left.
\mathfrak{T}_{n}\right\vert _{Y_{j-1}}\right)  $, we get
\[
\mathcal{I}_{j}/\mathcal{I}_{j-1}\cong C^{\ast}\left(  \left.  \mathfrak{T}%
_{n}\right\vert _{Y_{j}\backslash Y_{j-1}}\right)  \cong\oplus^{\frac
{n!}{j!\left(  n-j\right)  !}}\left(  \mathcal{K}\left(  \ell^{2}\left(
\mathbb{Z}_{\geq}^{n-j}\right)  \right)  \otimes C\left(  \mathbb{T}%
^{j}\right)  \right)  .
\]

$\square$

\textbf{Corollary 1}. There is a finite composition sequence of closed ideals
\[
C\left(  \mathbb{S}_{H}^{2n-1}\right)  \equiv C^{\ast}\left(  \mathfrak{G}%
_{n}\right)  \equiv\mathcal{J}_{n}\rhd\mathcal{J}_{n-1}\rhd\cdots
\rhd\mathcal{J}_{1}\rhd\mathcal{J}_{0}\equiv\left\{  0\right\}
\]
such that for $1\leq j\leq n$,
\[
\mathcal{J}_{j}/\mathcal{J}_{j-1}\cong\oplus^{\frac{n!}{j!\left(  n-j\right)
!}}\left(  \mathcal{K}\left(  \ell^{2}\left(  \mathbb{Z}_{\geq}^{n-j}\right)
\right)  \otimes C\left(  \mathbb{T}^{j}\right)  \right)  .
\]

Proof. With $\mathcal{I}_{0}=\mathcal{K}\left(  \ell^{2}\left(  \mathbb{Z}%
_{\geq}^{n}\right)  \right)  $ and hence $C^{\ast}\left(  \mathfrak{T}%
_{n}\right)  /\mathcal{I}_{0}\cong C\left(  \mathbb{S}_{H}^{2n-1}\right)  $,
we simply take $\mathcal{J}_{j}:=\mathcal{I}_{j}/\mathcal{I}_{0}$. $\square$

The above composition sequences lead to the straightforward estimates
\[
\left\lfloor \frac{n}{2}\right\rfloor +1\leq\operatorname{sr}\left(  C\left(
\mathbb{S}_{H}^{2n-1}\right)  \right)  \leq\operatorname{sr}\left(
\mathcal{T}^{\otimes n}\right)  \leq\left\lfloor \frac{n+1}{2}\right\rfloor
+1
\]
and
\[
\operatorname{csr}\left(  \mathcal{T}^{\otimes n}\right)  \leq
\operatorname{csr}\left(  C\left(  \mathbb{S}_{H}^{2n-1}\right)  \right)
\leq\left\lfloor \frac{n+1}{2}\right\rfloor +1
\]
for all $n\geq1$, based on the general rules established in \cite{Ri:dsr} that
(i) $\operatorname{sr}\left(  \mathcal{A}\otimes\mathcal{K}\right)
=\min\left\{  2,\operatorname{sr}\left(  \mathcal{A}\right)  \right\}  $, (ii)
for any closed ideal $\mathcal{I}$ of a C*-algebra $\mathcal{A}$,
\[
\max\left\{  \operatorname{sr}\left(  \mathcal{A}/\mathcal{I}\right)
,\operatorname{sr}\left(  \mathcal{I}\right)  \right\}  \leq\operatorname{sr}%
\left(  \mathcal{A}\right)  \leq\max\left\{  \operatorname{sr}\left(
\mathcal{A}/\mathcal{I}\right)  ,\operatorname{sr}\left(  \mathcal{I}\right)
,\operatorname{csr}\left(  \mathcal{A}/\mathcal{I}\right)  \right\}  ,
\]
and (iii) $\operatorname{sr}\left(  C\left(  X\right)  \right)  =\left\lfloor
\frac{n}{2}\right\rfloor +1$ for any $n$-dimensional CW-complex $X$, and the
rule \cite{Sh:ct,Na,Ni} that for any closed ideal $\mathcal{I}$ of a
C*-algebra $\mathcal{A}$, (iv) $\operatorname{csr}\left(  \mathcal{A}%
\otimes\mathcal{K}\right)  \leq2$ (with $\operatorname{csr}\left(
\mathcal{K}\right)  =1$) and (v)
\[
\operatorname{csr}\left(  \mathcal{A}\right)  \leq\max\left\{
\operatorname{csr}\left(  \mathcal{A}/\mathcal{I}\right)  ,\operatorname{csr}%
\left(  \mathcal{I}\right)  \right\}  .
\]

Indeed, for $n>1$, applying (i)-(ii) and (iv)-(v) to the short exact sequences%
\[
0\rightarrow\mathcal{I}_{j-1}\rightarrow\mathcal{I}_{j}\rightarrow
\mathcal{I}_{j}/\mathcal{I}_{j-1}\cong\oplus^{\frac{n!}{j!\left(  n-j\right)
!}}\left(  \mathcal{K}\left(  \ell^{2}\left(  \mathbb{Z}_{\geq}^{n-j}\right)
\right)  \otimes C\left(  \mathbb{T}^{j}\right)  \right)  \rightarrow0
\]
inductively for $j$ increasing from $1$ to $n-1$, starting with the exact
sequence
\[
0\rightarrow\mathcal{K}\left(  \ell^{2}\left(  \mathbb{Z}_{\geq}^{n}\right)
\right)  \cong\mathcal{I}_{0}\rightarrow\mathcal{I}_{1}\rightarrow
\mathcal{I}_{1}/\mathcal{I}_{0}\cong\oplus^{n}\left(  \mathcal{K}\left(
\ell^{2}\left(  \mathbb{Z}_{\geq}^{n-1}\right)  \right)  \otimes C\left(
\mathbb{T}\right)  \right)  \rightarrow0
\]
for $j=1$, we get $\operatorname{csr}\left(  \mathcal{I}_{j}\right)
,\operatorname{sr}\left(  \mathcal{I}_{j}\right)  \leq2$ for all $1\leq j\leq
n-1$. In particular, $\operatorname{csr}\left(  \mathcal{I}_{n-1}\right)
,\operatorname{sr}\left(  \mathcal{I}_{n-1}\right)  \leq2$, which is also
valid for $n=1$ since $\mathcal{I}_{0}\cong\mathcal{K}\left(  \ell^{2}\left(
\mathbb{Z}_{\geq}^{n}\right)  \right)  $. Then with $\operatorname{csr}\left(
C\left(  \mathbb{T}^{n}\right)  \right)  \leq\left\lfloor \frac{n+1}%
{2}\right\rfloor +1$ by homotopy theory \cite{Wh}, we get $\operatorname{csr}%
\left(  \mathcal{T}^{\otimes n}\right)  \leq\left\lfloor \frac{n+1}%
{2}\right\rfloor +1$ and%
\[
\left\lfloor \frac{n}{2}\right\rfloor +1\leq\operatorname{sr}\left(
\mathcal{T}^{\otimes n}\right)  \leq\left\lfloor \frac{n+1}{2}\right\rfloor
+1
\]
by further applying (ii)-(iii) and (v) to the short exact sequence
\[
0\rightarrow\mathcal{I}_{n-1}\rightarrow\mathcal{I}_{n}\equiv\mathcal{T}%
^{\otimes n}\rightarrow\mathcal{I}_{n}/\mathcal{I}_{n-1}\cong C\left(
\mathbb{T}^{n}\right)  \rightarrow0.
\]
Similar argument yields $\operatorname{csr}\left(  C\left(  \mathbb{S}%
_{H}^{2n-1}\right)  \right)  \leq\left\lfloor \frac{n+1}{2}\right\rfloor +1$
and $\left\lfloor \frac{n}{2}\right\rfloor +1\leq\operatorname{sr}\left(
C\left(  \mathbb{S}_{H}^{2n-1}\right)  \right)  \leq\left\lfloor \frac{n+1}%
{2}\right\rfloor +1$, with the inequality $\operatorname{sr}\left(  C\left(
\mathbb{S}_{H}^{2n-1}\right)  \right)  \leq\operatorname{sr}\left(
\mathcal{T}^{\otimes n}\right)  $ obviously valid by (ii). Also
$\operatorname{csr}\left(  \mathcal{T}^{\otimes n}\right)  \leq
\operatorname{csr}\left(  C\left(  \mathbb{S}_{H}^{2n-1}\right)  \right)  $ by (iv)-(v).

Such an estimate determining $\operatorname{sr}\left(  \mathcal{T}^{\otimes
n}\right)  $ sharply for even $n$ and up to an error of $1$ for odd $n>1$ as
stated above was first obtained by G. Nagy in \cite{Na} and then sharpened to
the exact value
\[
\operatorname{sr}\left(  \mathcal{T}^{\otimes n}\right)  =\left\lfloor
\frac{n}{2}\right\rfloor +1\text{\ \ (and hence }\operatorname{sr}\left(
C\left(  \mathbb{S}_{H}^{2n-1}\right)  \right)  =\left\lfloor \frac{n}%
{2}\right\rfloor +1\text{)}%
\]
for general $n>1$ by Nistor in \cite{Ni} which also gives $\operatorname{csr}%
\left(  \mathcal{T}^{\otimes n}\right)  \leq\left\lfloor \frac{n+1}%
{2}\right\rfloor +1 $. We summarize these results as follows.

\textbf{Proposition 2}. For all $n>1$,
\[
\operatorname{sr}\left(  C\left(  \mathbb{S}_{H}^{2n-1}\right)  \right)
=\operatorname{sr}\left(  \mathcal{T}^{\otimes n}\right)  =\left\lfloor
\frac{n}{2}\right\rfloor +1
\]
and
\[
\operatorname{csr}\left(  \mathcal{T}^{\otimes n}\right)  \leq
\operatorname{csr}\left(  C\left(  \mathbb{S}_{H}^{2n-1}\right)  \right)
\leq\left\lfloor \frac{n+1}{2}\right\rfloor +1.
\]

\textbf{Corollary 2}. For any $n>1$ and any $k\geq\left\lfloor \frac{n}%
{2}\right\rfloor +3$, the topological group $GL_{k}\left(  \mathcal{T}%
^{\otimes n}\right)  $ is connected.

Proof. By the K\"{u}nneth formula \cite{Blac} for $K$-groups, we get
$K_{1}\left(  \mathcal{T}^{\otimes n}\right)  =0$ since $K_{1}\left(
\mathcal{T}\right)  =0$ is well known. So by the theorem \cite{Ri:dsr} that
$K_{1}\left(  \mathcal{A}\right)  \cong GL_{k}\left(  \mathcal{A}\right)
/GL_{k}^{0}\left(  \mathcal{A}\right)  $ for any unital C*-algebra
$\mathcal{A}$ with $k\geq\operatorname{sr}\left(  \mathcal{A}\right)  +2$, we
get $GL_{k}\left(  \mathcal{T}^{\otimes n}\right)  =GL_{k}^{0}\left(
\mathcal{T}^{\otimes n}\right)  $ for any $k\geq\left\lfloor \frac{n}%
{2}\right\rfloor +3\geq\operatorname{sr}\left(  \mathcal{T}^{\otimes
n}\right)  +2$.

$\square$

Note that the above statement holds for the case of $n=1$, since
$GL_{k}\left(  \mathcal{T}\right)  $ is connected for all $k\geq1$ in the case
of $n=1$ by the index theory of Toeplitz operators for the unit disk
$\mathbb{D}$.

\section{Projective modules over $\mathcal{T}^{\otimes n}$}

Before proceeding to study finitely generated projective modules over
$\mathcal{T}^{\otimes n}$, we now point out a structure of $\mathcal{T}%
^{\otimes n}$ which facilitates some inductive procedures for the study of
such modules.

For all $n\in\mathbb{N}$, the topological groupoid $\left.  \mathfrak{T}%
_{n}\right\vert _{\overline{\mathbb{Z}_{\geq}}^{n-1}\times\left\{
\infty\right\}  }$ is isomorphic to the product topological groupoid
$\mathfrak{T}_{n-1}\times\mathbb{Z}$, while the topological groupoid $\left.
\mathfrak{T}_{n}\right\vert _{\overline{\mathbb{Z}_{\geq}}^{n-1}%
\times\mathbb{Z}_{\geq}}$ is isomorphic to the product topological groupoid
$\mathfrak{T}_{n-1}\times\left.  \left(  \mathbb{Z}\ltimes\mathbb{Z}\right)
\right\vert _{\mathbb{Z}_{\geq}}$, where the closed subset $\overline
{\mathbb{Z}_{\geq}}^{n-1}\times\left\{  \infty\right\}  $ and its open
complement $\overline{\mathbb{Z}_{\geq}}^{n-1}\times\mathbb{Z}_{\geq}$ in the
unit space $\overline{\mathbb{Z}_{\geq}}^{n}$ of $\mathfrak{T}_{n}$ are
invariant. (Here it is understood that when $n-1=0$, the first factor
$\overline{\mathbb{Z}_{\geq}}^{n-1}$ is dropped.) Hence we get the short exact
sequence of C*-algebras%
\begin{align*}
0  &  \rightarrow C^{\ast}\left(  \left.  \mathfrak{T}_{n}\right\vert
_{\overline{\mathbb{Z}_{\geq}}^{n-1}\times\mathbb{Z}_{\geq}}\right)
\cong\mathcal{T}^{\otimes n-1}\otimes\mathcal{K}\left(  \ell^{2}\left(
\mathbb{Z}_{\geq}\right)  \right)  \rightarrow\\
C^{\ast}\left(  \mathfrak{T}_{n}\right)   &  \equiv\mathcal{T}^{\otimes
n}\overset{\kappa_{n}}{\rightarrow}C^{\ast}\left(  \left.  \mathfrak{T}%
_{n}\right\vert _{\overline{\mathbb{Z}_{\geq}}^{n-1}\times\left\{
\infty\right\}  }\right)  \cong\mathcal{T}^{\otimes n-1}\otimes C\left(
\mathbb{T}\right)  \rightarrow0
\end{align*}
with $\mathcal{T}^{\otimes0}:=\mathbb{C}$. Furthermore the quotient maps
$\kappa_{n}$ for $n\in\mathbb{N}$ resulting from a groupoid restriction
satisfy the commuting diagram
\[%
\begin{array}
[c]{ccccc}%
M_{k}\left(  \mathcal{T}^{\otimes n}\right)  & \overset{\kappa_{n}%
}{\rightarrow} & M_{k}\left(  \mathcal{T}^{\otimes n-1}\otimes C\left(
\mathbb{T}\right)  \right)  & \equiv & M_{k}\left(  \mathcal{T}^{\otimes
n-1}\right)  \otimes C\left(  \mathbb{T}\right) \\
_{{}}\downarrow_{\sigma_{n}} & \circlearrowright & _{{}}\downarrow
_{\sigma_{n-1}\otimes\operatorname{id}} &  & \downarrow_{\sigma_{n-1}%
\otimes\operatorname{id}}\\
M_{k}\left(  C\left(  \mathbb{T}^{n}\right)  \right)  & \overset{\equiv
}{\rightarrow} & M_{k}\left(  C\left(  \mathbb{T}^{n-1}\right)  \otimes
C\left(  \mathbb{T}\right)  \right)  & \equiv & M_{k}\left(  C\left(
\mathbb{T}^{n-1}\right)  \right)  \otimes C\left(  \mathbb{T}\right)
\end{array}
\]
where $\equiv$ stands for a canonical isomorphism and $\sigma_{0}%
:=\operatorname{id}_{\mathbb{C}}$.

To classify the isomorphism classes of finitely generated projective
$\mathcal{T}^{\otimes n}$-modules $E$ or equivalently the equivalence classes
of idempotents $P\in M_{\infty}\left(  \mathcal{T}^{\otimes n}\right)  $ over
$\mathcal{T}^{\otimes n}$, we first define the rank of (the class of) $E$ or
$P$ as the classical rank of (the isomorphism class of) the vector bundle
corresponding to (the class of) the $C\left(  \mathbb{T}^{n}\right)  $-module
$C\left(  \mathbb{T}^{n}\right)  \otimes_{\mathcal{T}^{\otimes n}}E$ or the
projection $\sigma_{n}\left(  P\right)  $ over $C\left(  \mathbb{T}%
^{n}\right)  $.

The set of equivalence classes of idempotents $P\in M_{\infty}\left(
\mathcal{T}^{\otimes n}\right)  $ equipped with the binary operation
$\boxplus$ becomes an abelian graded monoid
\[
\mathfrak{P}\left(  \mathcal{T}^{\otimes n}\right)  =\sqcup_{m=0}^{\infty
}\mathfrak{P}_{m}\left(  \mathcal{T}^{\otimes n}\right)
\]
where $\mathfrak{P}_{m}\left(  \mathcal{T}^{\otimes n}\right)  $ is the set of
all (equivalence classes of) idempotents over $\mathcal{T}^{\otimes n}$ of
rank $m$, and
\[
\mathfrak{P}_{m}\left(  \mathcal{T}^{\otimes n}\right)  \boxplus
\mathfrak{P}_{l}\left(  \mathcal{T}^{\otimes n}\right)  \subset\mathfrak{P}%
_{m+l}\left(  \mathcal{T}^{\otimes n}\right)
\]
for $m,l\geq0$. Clearly $\mathfrak{P}_{0}\left(  \mathcal{T}^{\otimes
n}\right)  $ is a submonoid of $\mathfrak{P}\left(  \mathcal{T}^{\otimes
n}\right)  $.

Next we define a submonoid of $\mathfrak{P}\left(  \mathcal{T}^{\otimes
n}\right)  $ generated by \textquotedblleft standard\textquotedblright\ type
of idempotents, which turns out to contain (equivalence classes of) all
idempotents of sufficiently high ranks, and then classify its elements.

Note that each permutation $\Theta$ on $\left\{  1,2,...,n\right\}  $ induces
canonically a C*-algebra automorphism, still denoted as $\Theta$ by abuse of
notation, on $\mathcal{T}^{\otimes n}$ by permuting the indices of the factors
in $a_{1}\otimes a_{2}\otimes\cdots\otimes a_{n}\in\mathcal{T}^{\otimes n}$
for $a_{i}\in\mathcal{T}$. A permutation $\Theta$ on $\left\{
1,2,...,n\right\}  $ is called a $\left(  j,n-j\right)  $-shuffle on $\left\{
1,2,...,n\right\}  $ if $\Theta\left(  1\right)  <\Theta\left(  2\right)
<\cdots<\Theta\left(  j\right)  $ and $\Theta\left(  j+1\right)
<\Theta\left(  j+2\right)  <\cdots<\Theta\left(  n\right)  $.

Some basic projections over $\mathcal{T}^{\otimes n}$ are given by
$\Theta\left(  P_{j,l}\right)  $ where
\[
P_{j,l}:=\boxplus^{l}\left(  \left(  \otimes^{j}I\right)  \otimes\left(
\otimes^{n-j}P_{1}\right)  \right)  \in M_{l}\left(  \mathcal{T}^{\otimes
n}\right)
\]
for $l\geq0$ and $0\leq j\leq n$ (in particular, $P_{n,m}\equiv\boxplus
^{m}\left(  \otimes^{n}I\right)  \equiv\boxplus^{m}\tilde{I}$ for the unit
$\tilde{I}$ of $\mathcal{T}^{\otimes n}$), and $\Theta$ is (the automorphism
defined by) a $\left(  j,n-j\right)  $-shuffle on $\left\{  1,2,...,n\right\}
$. Note that $\Theta\left(  P_{j,l}\right)  =\Theta\left(  \boxplus^{l}%
P_{j,1}\right)  =\boxplus^{l}\Theta\left(  P_{j,1}\right)  $,
\[
\Theta\left(  P_{j,l}\right)  \boxplus\Theta\left(  P_{j,l^{\prime}}\right)
\sim\Theta\left(  P_{j,l+l^{\prime}}\right)  ,
\]
and $\left(  \otimes^{j}I\right)  \otimes\left(  \otimes^{n-j-1}P_{1}\right)
\otimes P_{l}\sim P_{j,l}$ over $\mathcal{T}^{\otimes n}$ since $P_{l}%
\sim\boxplus^{l}P_{1}$ over $\mathcal{K}^{+}\subset\mathcal{T}$. Furthermore
\[
\sigma_{n}\left(  \Theta\left(  P_{j,l}\right)  \right)  =\left\{
\begin{array}
[c]{lll}%
0, & \text{if } & 0\leq j\leq n-1\\
\boxplus^{l}1, & \text{if } & j=n
\end{array}
\right.  ,
\]
and hence $\Theta\left(  P_{j,l}\right)  \in\mathfrak{P}_{0}\left(
\mathcal{T}^{\otimes n}\right)  $ if $j<n$ and $\Theta\left(  P_{n,l}\right)
=P_{n,l}\in\mathfrak{P}_{l}\left(  \mathcal{T}^{\otimes n}\right)  $, where
$1\in C\left(  \mathbb{T}^{n}\right)  $ is the constant function $1$ on
$\mathbb{T}^{n}$. So the set $\mathfrak{P}_{0}^{\prime}\left(  \mathcal{T}%
^{\otimes n}\right)  \subset\mathfrak{P}_{0}\left(  \mathcal{T}^{\otimes
n}\right)  $ consisting of (the equivalence classes of) all possible
$\boxplus$-sums of $\Theta\left(  P_{j,l}\right)  $ with $l\geq0$ and $\Theta$
a $\left(  j,n-j\right)  $-shuffle on $\left\{  1,2,...,n\right\}  $ for
$0\leq j\leq n-1$ is a submonoid of $\mathfrak{P}_{0}\left(  \mathcal{T}%
^{\otimes n}\right)  $. For $m\geq1$, we define a singleton%
\[
\mathfrak{P}_{m}^{\prime}\left(  \mathcal{T}^{\otimes n}\right)  :=\left\{
P_{n,m}\equiv\boxplus^{m}\tilde{I}\right\}  \subset\mathfrak{P}_{m}\left(
\mathcal{T}^{\otimes n}\right)
\]
where $\tilde{I}$ denotes the identity element of $\mathcal{T}^{\otimes n}$.
Clearly $\sqcup_{m=1}^{\infty}\mathfrak{P}_{m}^{\prime}\left(  \mathcal{T}%
^{\otimes n}\right)  $ is also a submonoid of $\mathfrak{P}\left(
\mathcal{T}^{\otimes n}\right)  $.

We define a partial ordering $\prec$ on the collection
\[
\Omega:=\left\{  \left(  j,\Theta\right)  :\ 0\leq j\leq n\text{\ and}%
\ \Theta\text{\ is a }\left(  j,n-j\right)  \text{-shuffle}\right\}
\]
by the condition that $\left(  j^{\prime},\Theta^{\prime}\right)  \prec\left(
j,\Theta\right)  $ if and only if $\Theta\left(  \left\{  1,2,...,j\right\}
\right)  \supsetneqq\Theta^{\prime}\left(  \left\{  1,2,...,j^{\prime
}\right\}  \right)  $ (and hence $j>j^{\prime}$). Here $\left\{
1,2,...,0\right\}  \equiv\emptyset$ is understood. Note that
$\operatorname{id}_{\left\{  1,2,...,n\right\}  }$ is a $\left(  j,n-j\right)
$-shuffle for every $j$, and $\left(  n,\operatorname{id}_{\left\{
1,2,...,n\right\}  }\right)  $ is the greatest element while $\left(
0,\operatorname{id}_{\left\{  1,2,...,n\right\}  }\right)  $ is the smallest
element in $\Omega$ with respect to $\prec$.

\textbf{Proposition 3}. $\mathfrak{P}^{\prime}\left(  \mathcal{T}^{\otimes
n}\right)  =\sqcup_{m=0}^{\infty}\mathfrak{P}_{m}^{\prime}\left(
\mathcal{T}^{\otimes n}\right)  $ is a graded submonoid of $\mathfrak{P}%
\left(  \mathcal{T}^{\otimes n}\right)  $ and its monoid structure is
explicitly determined by that for any $l,l^{\prime}>0$ and any $\left(
j^{\prime},\Theta^{\prime}\right)  \prec\left(  j,\Theta\right)  $ in $\Omega
$,%
\[
\Theta\left(  P_{j,l}\right)  \boxplus\Theta^{\prime}\left(  P_{j^{\prime
},l^{\prime}}\right)  \sim\Theta\left(  P_{j,l}\right)  .
\]

Proof. Note that since $\mathfrak{P}_{0}^{\prime}\left(  \mathcal{T}^{\otimes
n}\right)  $ and $\boxplus_{m=1}^{\infty}\mathfrak{P}_{m}^{\prime}\left(
\mathcal{T}^{\otimes n}\right)  $ are submonoids of $\mathfrak{P}\left(
\mathcal{T}^{\otimes n}\right)  $, the set $\mathfrak{P}^{\prime}\left(
\mathcal{T}^{\otimes n}\right)  $ is a submonoid if $\Theta\left(
P_{n,m}\right)  \boxplus\Theta^{\prime}\left(  P_{j^{\prime},l^{\prime}%
}\right)  \sim\Theta\left(  P_{n,m}\right)  $ holds for all $m>0$ and all
$\Theta^{\prime}\left(  P_{j^{\prime},l^{\prime}}\right)  $ with $j^{\prime
}\leq n-1$. Since $\left(  n,\operatorname{id}_{\left\{  1,2,...,n\right\}
}\right)  $ is the greatest element in $\Omega$, it remains to show that
$\Theta\left(  P_{j,l}\right)  \boxplus\Theta^{\prime}\left(  P_{j^{\prime
},l^{\prime}}\right)  \sim\Theta\left(  P_{j,l}\right)  $ for $n\geq
j>j^{\prime}\geq0$ with $\Theta\left(  \left\{  1,2,...,j\right\}  \right)
\supset\Theta^{\prime}\left(  \left\{  1,2,...,j^{\prime}\right\}  \right)  $
and $l,l^{\prime}>0$.

Note that for $\Theta\left(  \left\{  1,2,...,j\right\}  \right)
\supset\Theta^{\prime}\left(  \left\{  1,2,...,j^{\prime}\right\}  \right)  $,
there exists a permutation $\Theta^{\prime\prime}$ (not necessarily a shuffle)
on $\left\{  1,2,...,n\right\}  $ such that $\Theta^{\prime\prime}\left(
\Theta\left(  P_{j,l}\right)  \right)  =P_{j,l}$ and $\Theta^{\prime\prime
}\left(  \Theta^{\prime}\left(  P_{j^{\prime},l^{\prime}}\right)  \right)
=P_{j^{\prime},l^{\prime}}$. (In fact, one can find a permutation
$\Theta^{\prime\prime}$ such that $\Theta^{\prime\prime}\Theta$ fixes each of
$j+1,..,n$, and $\Theta^{\prime\prime}\Theta^{\prime}$ is each of
$1,2,...,j^{\prime}$.) So it suffices to prove that
\[
P_{j,l}\boxplus P_{j^{\prime},l^{\prime}}\sim P_{j,l}%
\]
whenever $j>j^{\prime}$ and $l,l^{\prime}>0$. Furthermore since $P_{j,l}%
=\boxplus^{l}P_{j,1}$, we only need to show that $P_{j,1}\boxplus
P_{j^{\prime},1}\sim P_{j,1}$ for $j>j^{\prime}$.

Note that $U\left(  P_{1}\boxplus I\right)  U^{\ast}=0\boxplus I$ in
$M_{2}\left(  \mathcal{T}\right)  $ for the unitary
\[
U:=e_{11}\otimes\mathcal{S}^{\ast}+e_{22}\otimes\mathcal{S}+e_{21}\otimes
e_{11}\in M_{2}\left(  \mathbb{C}\right)  \otimes\mathcal{T}\equiv
M_{2}\left(  \mathcal{T}\right)
\]
where $\mathcal{S}\in\mathcal{T}$ is the (forward) unilateral shift on
$\ell^{2}\left(  \mathbb{Z}_{\geq}\right)  $. So
\[
P_{j,1}\boxplus P_{j-1,1}=\left(  \left(  \otimes^{j}I\right)  \otimes\left(
\otimes^{n-j}P_{1}\right)  \right)  \boxplus\left(  \left(  \otimes
^{j-1}I\right)  \otimes\left(  \otimes^{n-j+1}P_{1}\right)  \right)
\]%
\[
=\left(  \otimes^{j-1}I\right)  \otimes\left(  I\boxplus P_{1}\right)
\otimes\left(  \otimes^{n-j}P_{1}\right)  \sim\left(  \otimes^{j-1}I\right)
\otimes I\otimes\left(  \otimes^{n-j}P_{1}\right)  =P_{j,1}.
\]
Thus by iteration of this result, we can \textquotedblleft
expand\textquotedblright\ $P_{j,1}$ to get for any $0\leq k<j$,
\[
P_{j,1}\sim P_{j,1}\boxplus P_{j-1,1}\boxplus\cdots\boxplus P_{k,1},
\]
and hence
\[
P_{j,1}\boxplus P_{j^{\prime},1}\sim P_{j,1}\boxplus P_{j-1,1}\boxplus
\cdots\boxplus P_{j^{\prime}+1,1}\boxplus P_{j^{\prime},1}\sim P_{j,1}.
\]

$\square$

For each $\left(  j,\Theta\right)  \in\Omega$, let $X_{\Theta}\subset
\overline{\mathbb{Z}_{\geq}}^{n}$ be the invariant closed subset of the unit
space of $\mathfrak{T}_{n}$ consisting of $z\in\overline{\mathbb{Z}_{\geq}%
}^{n}$ with $z_{k}=\infty$ for all $k\in\Theta\left(  \left\{
1,2,...,j\right\}  \right)  $, and let
\[
\sigma_{\left(  j,\Theta\right)  }:C^{\ast}\left(  \mathfrak{T}_{n}\right)
\rightarrow C^{\ast}\left(  \left.  \mathfrak{T}_{n}\right\vert _{X_{\Theta}%
}\right)  \cong C\left(  \mathbb{T}^{j}\right)  \otimes\mathcal{T}^{\otimes
n-j}\subset C\left(  \mathbb{T}^{j}\right)  \otimes\mathcal{B}\left(  \ell
^{2}\left(  \mathbb{Z}_{\geq}^{n-j}\right)  \right)
\]
be the canonical quotient map, where the isomorphism implicitly involves a
rearrangement of factors by the inverse permutation $\Theta^{-1}$. Here as
before, $\mathbb{T}^{0}$ is a singleton. Defining $\rho_{\left(
j,\Theta\right)  }\left(  P\right)  $ for an idempotent $P$ over $C^{\ast
}\left(  \mathfrak{T}_{n}\right)  $ as the rank of the projection operator
$\sigma_{\left(  j,\Theta\right)  }\left(  P\right)  \left(  t\right)
\in\mathcal{B}\left(  \ell^{2}\left(  \mathbb{Z}_{\geq}^{n-j}\right)  \right)
$ for any $t\in\mathbb{T}^{j}$, which depends only on the equivalence class of
$P$, we get a well-defined monoid homomorphism
\[
\rho_{\left(  j,\Theta\right)  }:\left(  \mathfrak{P}\left(  \mathcal{T}%
^{\otimes n}\right)  ,\boxplus\right)  \rightarrow\left(  \mathbb{Z}_{\geq
}\cup\left\{  \infty\right\}  ,+\right)  .
\]

A (finite) $\boxplus$-sum of (the equivalence classes of) projections
$\Theta\left(  P_{j,l}\right)  $ indexed by some $\left(  j,\Theta\right)
\in\Omega$ that are mutually unrelated by $\prec$ with $l\equiv l_{\left(
j,\Theta\right)  }>0$ depending on $\left(  j,\Theta\right)  $ is called a
reduced $\boxplus$-sum of standard projections over $\mathcal{T}^{\otimes n}$.
It is understood that an \textquotedblleft empty\textquotedblright\ $\boxplus
$-sum represents the zero projection and is a reduced $\boxplus$-sum. Two
reduced $\boxplus$-sums are called different when they have different sets of
(mutually $\prec$-unrelated) indices $\left(  j,\Theta\right)  \in\Omega$ or
have different weight functions $l$ of $\left(  j,\Theta\right)  $. We are
going to show that different reduced $\boxplus$-sums are inequivalent
projections. Clearly each projection $\Theta\left(  P_{j,l}\right)  $ with
$\left(  j,\Theta\right)  \in\Omega$ and $l>0$ is a reduced $\boxplus$-sum.

\textbf{Theorem 1}. The submonoid $\mathfrak{P}^{\prime}\left(  \mathcal{T}%
^{\otimes n}\right)  =\sqcup_{m=0}^{\infty}\mathfrak{P}_{m}^{\prime}\left(
\mathcal{T}^{\otimes n}\right)  $ of $\mathfrak{P}\left(  \mathcal{T}^{\otimes
n}\right)  $ consists exactly of reduced $\boxplus$-sums of standard
projections over $\mathcal{T}^{\otimes n}$, and different reduced $\boxplus
$-sums are mutually inequivalent projections. Furthermore the monoid
homomorphism
\[
\rho:P\in\mathfrak{P}^{\prime}\left(  \mathcal{T}^{\otimes n}\right)
\mapsto\prod_{\left(  j,\Theta\right)  \in\Omega}\rho_{\left(  j,\Theta
\right)  }\left(  P\right)  \in\prod_{\left(  j,\Theta\right)  \in\Omega
}\overline{\mathbb{Z}_{\geq}}%
\]
is injective, with $\rho_{\left(  j,\Theta\right)  }\left(  \Theta\left(
P_{j,l}\right)  \right)  =l\in\mathbb{N}$.

Proof. By definition, $\mathfrak{P}^{\prime}\left(  \mathcal{T}^{\otimes
n}\right)  $ consists of $\boxplus$-sums of (the equivalence classes of)
projections $\Theta\left(  P_{j,l}\right)  $ with $\left(  j,\Theta\right)
\in\Omega$ and $l>0$. Since $\Theta\left(  P_{j,l}\right)  +\Theta\left(
P_{j,l^{\prime}}\right)  \sim\Theta\left(  P_{j,l+l^{\prime}}\right)  $, we
only need to consider in the following those $\boxplus$-sums, in which all
summands $\Theta\left(  P_{j,l}\right)  $ are indexed by distinct $\left(
j,\Theta\right)  \in\Omega$ with $l$ depending on $\left(  j,\Theta\right)  $.
For any such a $\boxplus$-sum, using the property that $\Theta\left(
P_{j,l}\right)  \boxplus\Theta^{\prime}\left(  P_{j^{\prime},l^{\prime}%
}\right)  \sim\Theta\left(  P_{j,l}\right)  $ for any $\left(  j^{\prime
},\Theta^{\prime}\right)  \prec\left(  j,\Theta\right)  $, we can remove one
by one those $\boxplus$-summands $\Theta^{\prime}\left(  P_{j^{\prime
},l^{\prime}}\right)  $ with $\left(  j^{\prime},\Theta^{\prime}\right)  $
dominated by the index of another summand, without changing the equivalence
class, until we reach a $\boxplus$-sum of $\Theta\left(  P_{j,l}\right)  $
with $\left(  j,\Theta\right)  \in\Omega$ mutually unrelated by $\prec$, i.e.
a reduced $\boxplus$-sum. So $\mathfrak{P}^{\prime}\left(  \mathcal{T}%
^{\otimes n}\right)  $ consists of the reduced $\boxplus$-sums.

Note that for $\left(  j,\Theta\right)  \in\Omega$ and $l>0$,
\begin{align*}
\sigma_{\left(  j,\Theta\right)  }\left(  \Theta\left(  P_{j,l}\right)
\right)   &  =\sigma_{\left(  j,\Theta\right)  }\left(  \boxplus^{l}%
\Theta\left(  \left(  \otimes^{j}I\right)  \otimes\left(  \otimes^{n-j}%
P_{1}\right)  \right)  \right) \\
&  =1\otimes\left(  \boxplus^{l}\left(  \otimes^{n-j}P_{1}\right)  \right)
\in C\left(  \mathbb{T}^{j}\right)  \otimes\left(  \boxplus^{l}\mathcal{B}%
\left(  \ell^{2}\left(  \mathbb{Z}_{\geq}^{n-j}\right)  \right)  \right)
\end{align*}
and hence $\rho_{\left(  j,\Theta\right)  }\left(  \Theta\left(
P_{j,l}\right)  \right)  =l\in\mathbb{N}$ the operator rank of $\boxplus
^{l}\left(  \otimes^{n-j}P_{1}\right)  \in\mathcal{B}\left(  \oplus^{l}%
\ell^{2}\left(  \mathbb{Z}_{\geq}^{n-j}\right)  \right)  $. But for $\left(
j^{\prime},\Theta^{\prime}\right)  \neq\left(  j,\Theta\right)  $,%
\[
\rho_{\left(  j,\Theta\right)  }\left(  \Theta^{\prime}\left(  P_{j^{\prime
},l^{\prime}}\right)  \right)  :=\left\{
\begin{array}
[c]{lll}%
\infty, & \text{if } & \left(  j,\Theta\right)  \prec\left(  j^{\prime}%
,\Theta^{\prime}\right) \\
0, & \text{if } & \text{otherwise}%
\end{array}
\right.
\]
because either $\sigma_{\left(  j,\Theta\right)  }\left(  \Theta^{\prime
}\left(  P_{j^{\prime},l^{\prime}}\right)  \right)  =0$ when $\Theta\left(
\left\{  1,2,...,j\right\}  \right)  \backslash\Theta^{\prime}\left(  \left\{
1,2,...,j^{\prime}\right\}  \right)  \neq\emptyset$, or $\sigma_{\left(
j,\Theta\right)  }\left(  \Theta^{\prime}\left(  P_{j^{\prime},l^{\prime}%
}\right)  \right)  $ is an infinite-dimensional projection when $\Theta
^{\prime}\left(  \left\{  1,2,...,j^{\prime}\right\}  \right)  \supset
\Theta\left(  \left\{  1,2,...,j\right\}  \right)  $ (but $\Theta\left(
\left\{  1,2,...,j\right\}  \right)  \neq\Theta^{\prime}\left(  \left\{
1,2,...,j^{\prime}\right\}  \right)  $ since $\left(  j^{\prime}%
,\Theta^{\prime}\right)  \neq\left(  j,\Theta\right)  $), i.e. when $\left(
j,\Theta\right)  \prec\left(  j^{\prime},\Theta^{\prime}\right)  $.

For a reduced $\boxplus$-sum $P$ of $\Theta^{\prime}\left(  P_{j^{\prime
},l^{\prime}}\right)  $ indexed by $\left(  j^{\prime},\Theta^{\prime}\right)
$ in some subset $A\subset\Omega$, the $\left(  j,\Theta\right)  $-component
of $\rho\left(  P\right)  $ is%
\[
\sum_{\left(  j^{\prime},\Theta^{\prime}\right)  \in A}\rho_{\left(
j,\Theta\right)  }\left(  \Theta^{\prime}\left(  P_{j^{\prime},l^{\prime}%
}\right)  \right)  \left\{
\begin{array}
[c]{lll}%
=l\in\mathbb{N} & \text{if } & \left(  j,\Theta\right)  \in A\text{\ with
}\Theta\left(  P_{j,l}\right)  \text{ a summand of }P\\
\in\left\{  0,\infty\right\}  & \text{if } & \text{otherwise}%
\end{array}
\right.
\]
for any $\left(  j,\Theta\right)  \in\Omega$, since if $\left(  j,\Theta
\right)  \in A$ then $\left(  j,\Theta\right)  $ is $\prec$-unrelated to any
other $\left(  j^{\prime},\Theta^{\prime}\right)  \in A$. So $\rho\left(
P\right)  $ completely determines the summands of a reduced $\boxplus$-sum
$P$, namely, $P\ $is the $\boxplus$-sum of exactly those $\Theta\left(
P_{j,l}\right)  $ with $l$ equal to the $\left(  j,\Theta\right)  $-component
of $\rho\left(  P\right)  $ that is a strictly positive integer. Since
$\mathfrak{P}^{\prime}\left(  \mathcal{T}^{\otimes n}\right)  $ consists of
reduced $\boxplus$-sums, this also shows that the clearly well-defined monoid
homomorphism $\rho$ is injective.

Thus if $P\sim P^{\prime}$ for two reduced $\boxplus$-sums $P$ and $P^{\prime
}$ and hence $\rho\left(  P\right)  =\rho\left(  P^{\prime}\right)  $, then
the summands of $P$ and $P^{\prime}$ are exactly the same, i.e. $P$ and
$P^{\prime}$ are the same reduced $\boxplus$-sum. So different reduced
$\boxplus$-sums are mutually inequivalent projections.

$\square$

\textbf{Proposition 4}. $\mathfrak{P}\left(  \mathcal{T}\right)
=\mathfrak{P}^{\prime}\left(  \mathcal{T}\right)  $. More concretely,
\[
\mathfrak{P}\left(  \mathcal{T}\right)  \cong\left\{  \left(  0,l\right)
:l\in\mathbb{Z}_{\geq}\right\}  \cup\left\{  \left(  m,\infty\right)
:m>0\right\}  \subset\overline{\mathbf{Z}_{\geq}}^{2}%
\]
where $\overline{\mathbf{Z}_{\geq}}^{2}$ is equipped with the canonical monoid structure.

Proof. It suffices to show that any element of $\mathfrak{P}_{0}\left(
\mathcal{T}\right)  \equiv\mathfrak{P}_{0}\left(  \mathcal{T}^{\otimes
1}\right)  $ is of the form $P_{0,l}$ (realized as $\left(  0,l\right)
\in\overline{\mathbf{Z}_{\geq}}^{2}$) and any element of $\mathfrak{P}%
_{m}\left(  \mathcal{T}\right)  \equiv\mathfrak{P}_{m}\left(  \mathcal{T}%
^{\otimes1}\right)  $ for $m\in\mathbb{N}$ is of the form $P_{1,m}$ (realized
as $\left(  m,\infty\right)  \in\overline{\mathbf{Z}_{\geq}}^{2}$).

The argument sketched below is similar to one used in \cite{Sh:pmqpl}.

Since any complex vector bundle over $\mathbb{T}$ is trivial, any idempotent
over $C\left(  \mathbb{T}\right)  $ is equivalent to the standard projection
$1\otimes P_{m}\in C\left(  \mathbb{T}\right)  \otimes M_{\infty}\left(
\mathbb{C}\right)  $ for some $m\in\mathbb{Z}_{\geq}$. So for any idempotent
$P\in M_{\infty}\left(  \mathcal{T}\right)  $ over $\mathcal{T}$, there is
some $U\in GL_{\infty}\left(  C\left(  \mathbb{T}\right)  \right)  $ such
that
\[
U\sigma\left(  P\right)  U^{-1}=1\otimes P_{m}=\sigma\left(  \boxplus
^{m}I\right)
\]
for some $m\in\mathbb{Z}_{\geq}$ where $I$ is the identity of $\mathcal{K}%
^{+}\subset\mathcal{T}$, and hence $VPV^{-1}-\boxplus^{m}I\in M_{\infty
}\left(  \mathcal{K}\right)  $ for any lift $V\in GL_{\infty}\left(
\mathcal{T}\right)  $ (which exists) of $U\boxplus U^{-1}\in GL_{\infty}%
^{0}\left(  C\left(  \mathbb{T}\right)  \right)  $ along $\sigma$. Replacing
$P$ by the equivalent $VPV^{-1}$, we may assume that $P\in\left(  \boxplus
^{m}I\right)  +M_{k-1}\left(  \mathcal{K}\right)  \subset M_{k-1}\left(
\mathcal{K}^{+}\right)  $ for some large $k\geq m+1$. Now since $M_{\infty
}\left(  \mathbb{C}\right)  $ is dense in $\mathcal{K}$, there is an
idempotent $Q\in\left(  \boxplus^{m}I\right)  +M_{k-1}\left(  M_{N}\left(
\mathbb{C}\right)  \right)  $ sufficiently close to and hence equivalent to
$P$ for some large $N $. So replacing $P$ by $Q$, we may assume that
$K:=P-\boxplus^{m}I\in M_{k-1}\left(  M_{N}\left(  \mathbb{C}\right)  \right)
$.

Rearranging the entries of $P\equiv K+\boxplus^{m}I\in M_{k-1}\left(
\mathcal{T}\right)  \subset M_{k}\left(  \mathcal{T}\right)  $ via conjugation
by the unitary
\[
U_{k,N}:=\sum_{j=1}^{k-1}\left(  e_{jj}\otimes\left(  \mathcal{S}^{\ast
}\right)  ^{N}+e_{kj}\otimes\left(  \mathcal{S}^{\left(  j-1\right)  N}%
P_{N}\right)  \right)  +e_{kk}\otimes\mathcal{S}^{\left(  k-1\right)  N}\in
M_{k}\left(  \mathbb{C}\right)  \otimes\mathcal{T}\equiv M_{k}\left(
\mathcal{T}\right)
\]
we get%
\[
U_{k,N}PU_{k,N}^{-1}\equiv U_{k,N}\left(  P\boxplus0\right)  U_{k,N}%
^{-1}=\left(  \left(  \boxplus^{m}I\right)  \boxplus\left(  \boxplus
^{k-1-m}0\right)  \right)  \boxplus R
\]
for some $R\in M_{\left(  k-1\right)  N}\left(  \mathbb{C}\right)
\subset\mathcal{K}\subset\mathcal{T}$ which must be an idempotent. Since any
idempotent in $\mathcal{K}$ is equivalent over $\mathcal{K}^{+}$ to a standard
projection $P_{l}$, we get
\[
P\sim\left(  \left(  \boxplus^{m}I\right)  \boxplus\left(  \boxplus
^{k-1-m}0\right)  \right)  \boxplus P_{l}%
\]
for some $l\in\mathbb{Z}_{\geq}$.

If $m=0$, then clearly $P\sim P_{l}$. Since it is well known that $P_{l}$ and
$\boxplus^{l}P_{1}\equiv P_{0,l}$ are equivalent over $\mathcal{K}^{+}$ and
hence over $\mathcal{T}\supset\mathcal{K}^{+}$, we get $P\sim P_{0,l}$.

If $m\in\mathbb{N}$, then we can rearrange entries via conjugation by the
unitary%
\[
U_{l}:=e_{11}\otimes\mathcal{S}^{l}+e_{1k}\otimes P_{l}+\sum_{j=2}^{k-1}%
e_{jj}\otimes I+e_{kk}\otimes\left(  \mathcal{S}^{\ast}\right)  ^{l}\in
M_{k}\left(  \mathbb{C}\right)  \otimes\mathcal{T}\equiv M_{k}\left(
\mathcal{T}\right)
\]
to get
\[
U_{l}\left(  \left(  \left(  \boxplus^{m}I\right)  \boxplus\left(
\boxplus^{k-1-m}0\right)  \right)  \boxplus P_{l}\right)  U_{l}^{-1}=\left(
\boxplus^{m}I\right)  \boxplus\left(  \boxplus^{k-m}0\right)  \equiv
\boxplus^{m}I\equiv P_{1,m}.
\]

$\square$

\textbf{Theorem 2}. For $n>1$ and $m>0$, if $\mathfrak{P}_{m}\left(
\mathcal{T}^{\otimes n-1}\right)  =\mathfrak{P}_{m}^{\prime}\left(
\mathcal{T}^{\otimes n-1}\right)  \equiv\left\{  \boxplus^{m}\left(
\otimes^{n-1}I\right)  \right\}  $ and $GL_{m}\left(  \mathcal{T}^{\otimes
n-1}\right)  $ is connected, then $\mathfrak{P}_{m}\left(  \mathcal{T}%
^{\otimes n}\right)  =\mathfrak{P}_{m}^{\prime}\left(  \mathcal{T}^{\otimes
n}\right)  $.

Proof. In this proof, we use $I$ and $\tilde{I}$ to denote respectively the
identity elements of $\mathcal{T}^{\otimes n-1}$ and $\mathcal{T}^{\otimes n}
$.

Let $P\in\mathfrak{P}_{m}\left(  \mathcal{T}^{\otimes n}\right)  $. The
idempotent $\kappa_{n}\left(  P\right)  $ over $\mathcal{T}^{\otimes
n-1}\otimes C\left(  \mathbb{T}\right)  $ satisfies that for any
$z\in\mathbb{T}$,
\[
\sigma_{n-1}\left(  \kappa_{n}\left(  P\right)  \left(  z\right)  \right)
=\sigma_{n}\left(  P\right)  \left(  \cdot,z\right)  \in M_{\infty}\left(
C\left(  \mathbb{T}^{n-1}\right)  \right)
\]
which is of rank $m$ pointwise, and hence
\[
\kappa_{n}\left(  P\right)  \left(  z\right)  \in\mathfrak{P}_{m}\left(
\mathcal{T}^{\otimes n-1}\right)  =\mathfrak{P}_{m}^{\prime}\left(
\mathcal{T}^{\otimes n-1}\right)  ,
\]
i.e. $\kappa_{n}\left(  P\right)  \left(  z\right)  \sim\boxplus^{m}I$ over
$\mathcal{T}^{\otimes n-1}$. In particular, there is a continuous
idempotent-valued path $\gamma:\left[  0,1\right]  \rightarrow M_{k}\left(
\mathcal{T}^{\otimes n-1}\right)  $ for $k$ sufficiently large going from the
idempotent $\kappa_{n}\left(  P\right)  \left(  1\right)  $ to $\left(
\boxplus^{m}I\right)  \boxplus\left(  \boxplus^{k-m}0\right)  $. Clearly we
may assume that $\gamma$ is locally constant at $1$, say, $\gamma\left(
t\right)  =\boxplus^{m}I$ for $t\geq1/2$. The concatenation of the path
$\gamma^{-1}$, the loop $\kappa_{n}\left(  P\right)  $, and the path $\gamma$
defines an idempotent-valued continuous loop $\Gamma:\mathbb{T}\rightarrow
M_{k}\left(  \mathcal{T}^{\otimes n-1}\right)  $ starting and ending at
$\boxplus^{m}I$ with $\Gamma\left(  e^{i\theta}\right)  =\left(  \boxplus
^{m}I\right)  \boxplus\left(  \boxplus^{k-m}0\right)  $, say, for all
$\theta\in\left[  3\pi/2,2\pi\right]  $ (and $\left[  0,\pi/2\right]  $), and
is homotopic to the loop $\kappa_{n}\left(  P\right)  $ via idempotents, i.e.
there is a path of idempotents in $M_{k}\left(  \mathcal{T}^{\otimes
n-1}\otimes C\left(  \mathbb{T}\right)  \right)  $ from $\kappa_{n}\left(
P\right)  $ to $\Gamma$. Consequently, there is a continuous path of
invertibles $U_{t}\in GL_{k}\left(  \mathcal{T}^{\otimes n-1}\otimes C\left(
\mathbb{T}\right)  \right)  $ with $U_{0}=I_{k}$ such that $U_{1}\kappa
_{n}\left(  P\right)  U_{1}^{-1}=\Gamma$, which can be lifted along
$\kappa_{n}$ to a continuous path of invertible $V_{t}\in GL_{k}\left(
\mathcal{T}^{\otimes n}\right)  $ with $V_{0}=I_{k}$ such that $\kappa
_{n}\left(  V_{1}PV_{1}^{-1}\right)  =\Gamma$.

Replacing $P$ by the equivalent idempotent $V_{1}PV_{1}^{-1}$, we may now
assume directly that the idempotent $\kappa_{n}\left(  P\right)  $ over
$\mathcal{T}^{\otimes n-1}\otimes C\left(  \mathbb{T}\right)  $ is a
continuous loop of idempotents in $M_{k}\left(  \mathcal{T}^{\otimes
n-1}\right)  $ such that $\kappa_{n}\left(  P\right)  \left(  e^{i\theta
}\right)  =\left(  \boxplus^{m}I\right)  \boxplus\left(  \boxplus
^{k-m}0\right)  $ for all $\theta\in\left[  3\pi/2,2\pi\right]  $. So there is
a continuous path
\[
\theta\in\left[  0,3\pi/2\right]  \mapsto W_{\theta}\in GL_{k}\left(
\mathcal{T}^{\otimes n-1}\right)
\]
with $W_{0}=I_{k}$ such that
\[
W_{\theta}\left(  \kappa_{n}\left(  P\right)  \left(  e^{i\theta}\right)
\right)  W_{\theta}^{-1}=\kappa_{n}\left(  P\right)  \left(  1\right)
=\left(  \boxplus^{m}I\right)  \boxplus\left(  \boxplus^{k-m}0\right)
\]
for all $\theta\in\left[  0,3\pi/2\right]  $. In particular,
\[
W_{3\pi/2}\left(  \left(  \boxplus^{m}I\right)  \boxplus\left(  \boxplus
^{k-m}0\right)  \right)  =\left(  \left(  \boxplus^{m}I\right)  \boxplus
\left(  \boxplus^{k-m}0\right)  \right)  W_{3\pi/2}%
\]
and hence $W_{3\pi/2}=W^{\prime}\boxplus W^{\prime\prime}$ for some
invertibles $W^{\prime}\in GL_{m}\left(  \mathcal{T}^{\otimes n-1}\right)  $
and $W^{\prime\prime}\in GL_{k-m}\left(  \mathcal{T}^{\otimes n-1}\right)  $.

By the connectedness assumption on $GL_{m}\left(  \mathcal{T}^{\otimes
n-1}\right)  $, there is a continuous path $\alpha:\left[  3\pi/2,2\pi\right]
\rightarrow GL_{m}\left(  \mathcal{T}^{\otimes n-1}\right)  $ with
$\alpha\left(  3\pi/2\right)  =W^{\prime}$ and $\alpha\left(  2\pi\right)
=I_{m}$. Since by K\"{u}nneth formula, $K_{1}\left(  \mathcal{T}^{\otimes
n-1}\right)  =0$ and hence $GL_{N}\left(  \mathcal{T}^{\otimes n-1}\right)  $
is connected for $N$ sufficiently large, we may suitably increase the value of
$k$ by adding diagonal $\boxplus$-summands $0$ to idempotents and diagonal
$\boxplus$-summands $I$ to invertibles, so that $GL_{k-m}\left(
\mathcal{T}^{\otimes n-1}\right)  $ is also connected and hence there is a
continuous path $\beta:\left[  3\pi/2,2\pi\right]  \rightarrow GL_{k-m}\left(
\mathcal{T}^{\otimes n-1}\right)  $ with $\beta\left(  3\pi/2\right)
=W^{\prime\prime}$ and $\beta\left(  2\pi\right)  =I_{k-m}$. Now the function
$\theta\mapsto W_{\theta}$ can be continuously extended to the whole interval
$\left[  0,2\pi\right]  $ by setting
\[
W_{\theta}:=\alpha\left(  \theta\right)  \boxplus\beta\left(  \theta\right)
\in GL_{k}\left(  \mathcal{T}^{\otimes n-1}\right)
\]
for $\theta\in\left[  3\pi/2,2\pi\right]  $, giving rise to a well-defined
continuous loop
\[
W:e^{i\theta}\in\mathbb{T}\mapsto W_{\theta}\in GL_{k}\left(  \mathcal{T}%
^{\otimes n-1}\right)  ,
\]
i.e. $W\in GL_{k}\left(  \mathcal{T}^{\otimes n-1}\otimes C\left(
\mathbb{T}\right)  \right)  $, satisfying
\[
W\left(  \kappa_{n}\left(  P\right)  \right)  W^{-1}=\left(  \boxplus
^{m}I\right)  \boxplus\left(  \boxplus^{k-m}0\right)  .
\]
So the idempotent $\kappa_{n}\left(  P\right)  $ over $\mathcal{T}^{\otimes
n-1}\otimes C\left(  \mathbb{T}\right)  $ is equivalent to the idempotent
$\boxplus^{m}I$.

Replacing $P$ by the equivalent idempotent $\tilde{W}\left(  P\boxplus\left(
\boxplus^{k}0\right)  \right)  \tilde{W}^{-1}$ for any fixed lifting
$\tilde{W}\in GL_{2k}^{0}\left(  \mathcal{T}^{\otimes n}\right)  $ of
$W\boxplus W^{-1}\in GL_{2k}^{0}\left(  \mathcal{T}^{\otimes n-1}\otimes
C\left(  \mathbb{T}\right)  \right)  $ along $\kappa_{n}$, we may now assume
that
\[
\kappa_{n}\left(  P\right)  =\left(  \boxplus^{m}I\right)  \boxplus\left(
\boxplus^{2k-m}0\right)  =\kappa_{n}\left(  \left(  \boxplus^{m}\tilde
{I}\right)  \boxplus\left(  \boxplus^{2k-m}0\right)  \right)
\]
and proceed to show that $P\sim\boxplus^{m}\tilde{I}$.

Note that $P-\left(  \left(  \boxplus^{m}\tilde{I}\right)  \boxplus\left(
\boxplus^{2k-m}0\right)  \right)  \in M_{2k}\left(  \mathcal{T}^{\otimes
n-1}\otimes\mathcal{K}\right)  $. Since $M_{\infty}\left(  \mathbb{C}\right)
$ is dense in $\mathcal{K}$, we may replace $P$ by a suitable equivalent
idempotent and assume that
\[
K:=P-\left(  \left(  \boxplus^{m}\tilde{I}\right)  \boxplus\left(
\boxplus^{2k-m}0\right)  \right)  \in M_{2k}\left(  \mathcal{T}^{\otimes
n-1}\otimes M_{N}\left(  \mathbb{C}\right)  \right)  \subset M_{2k}\left(
\mathcal{T}^{\otimes n}\right)  \text{\ }%
\]
for some $N\in\mathbb{N}$.

Rearranging the entries of $P\equiv P\boxplus0\in M_{2k+1}\left(
\mathcal{T}^{\otimes n-1}\otimes M_{N}\left(  \mathbb{C}\right)  \right)  $
via conjugation by the unitary
\begin{align*}
U_{k,N}  &  :=\sum_{j=1}^{2k}\left(  e_{jj}\otimes\left(  I\otimes
\mathcal{S}^{\ast}\right)  ^{N}+e_{2k+1,j}\otimes\left(  I\otimes
\mathcal{S}^{\left(  j-1\right)  N}P_{N}\right)  \right)  +e_{2k+1,2k+1}%
\otimes\left(  I\otimes\mathcal{S}^{2kN}\right) \\
&  \in M_{2k+1}\left(  \mathbb{C}\right)  \otimes\mathcal{T}^{\otimes
n-1}\otimes\mathcal{T}\equiv M_{2k+1}\left(  \mathcal{T}^{\otimes n}\right)
\end{align*}
we get%
\[
P\sim U_{k,N}PU_{k,N}^{-1}\equiv U_{k,N}\left(  P\boxplus0\right)
U_{k,N}^{-1}=\left(  \left(  \boxplus^{m}\tilde{I}\right)  \boxplus\left(
\boxplus^{2k-m}0\right)  \right)  \boxplus R
\]
for some
\[
R\in M_{2kN}\left(  \mathcal{T}^{\otimes n-1}\right)  \equiv\mathcal{T}%
^{\otimes n-1}\otimes M_{2kN}\left(  \mathbb{C}\right)  \subset\mathcal{T}%
^{\otimes n-1}\otimes\mathcal{T}\equiv\mathcal{T}^{\otimes n}%
\]
which must be an idempotent over $\mathcal{T}^{\otimes n-1}$.

Since $K_{0}\left(  \mathcal{T}^{\otimes n-1}\right)  =\mathbb{Z}$ by
K\"{u}nneth formula, $R\boxplus\left(  \boxplus^{r}I\right)  \sim\left(
\boxplus^{r+\left[  R\right]  }I\right)  $ for a sufficiently large
$r\in\mathbb{N}$ where $\left[  R\right]  \in\mathbb{Z}$ denotes the class of
$R$ in $K_{0}\left(  \mathcal{T}^{\otimes n-1}\right)  $. So there is an
invertible $U\in GL_{d}\left(  \mathcal{T}^{\otimes n-1}\right)  $ for some
large $d\geq\max\left\{  2kN+r,r+\left[  R\right]  \right\}  $ such that
\[
U\left(  R\boxplus\left(  \boxplus^{d-2kN-r}0\right)  \boxplus\left(
\boxplus^{r}I\right)  \right)  U^{-1}=\left(  \boxplus^{d-r-\left[  R\right]
}0\right)  \boxplus\left(  \boxplus^{r+\left[  R\right]  }I\right)  .
\]

With $m>0$, we can rearrange entries via conjugation by the unitary%
\begin{align*}
U_{d-r}  &  :=e_{11}\otimes\left(  I\otimes\mathcal{S}^{d-r}\right)
+e_{1,2k+1}\otimes I\otimes P_{d-r}+\sum_{j=2}^{2k}e_{jj}\otimes\tilde
{I}+e_{2k+1,2k+1}\otimes\left(  I\otimes\mathcal{S}^{\ast}\right)  ^{d-r}\\
&  \in M_{2k+1}\left(  \mathbb{C}\right)  \otimes\mathcal{T}^{\otimes
n-1}\otimes\mathcal{T}\equiv M_{2k+1}\left(  \mathcal{T}^{\otimes n}\right)
\end{align*}
to get
\[
P\sim U_{d-r}\left(  \left(  \boxplus^{m}\tilde{I}\right)  \boxplus\left(
\boxplus^{2k-m}0\right)  \boxplus R\right)  U_{d-r}^{-1}=R^{\prime}%
\boxplus\left(  \boxplus^{m-1}\tilde{I}\right)  \boxplus\left(  \boxplus
^{2k+1-m}0\right)
\]
where
\begin{align*}
R^{\prime}  &  =\left(  R\boxplus\left(  \boxplus^{d-2kN-r}0\right)  \right)
+\left(  \tilde{I}-I\otimes P_{d-r}\right)  \in\tilde{I}+\left(
\mathcal{T}^{\otimes n-1}\otimes M_{d-r}\left(  \mathbb{C}\right)  \right) \\
&  \subset\tilde{I}+\left(  \mathcal{T}^{\otimes n-1}\otimes\mathcal{K}%
\right)  \subset\mathcal{T}^{\otimes n-1}\otimes\mathcal{T}=\mathcal{T}%
^{\otimes n}.
\end{align*}
Note that $R^{\prime}$ can be interpreted as $R\boxplus\left(  \boxplus
^{d-2kN-r}0\right)  \boxplus\left(  \boxplus^{\infty}I\right)  \in
\mathcal{T}^{\otimes n-1}\otimes\mathcal{K}^{+}\subset\mathcal{T}^{\otimes n}%
$, which when conjugated by the invertible $U\equiv U\boxplus\left(
\boxplus^{\infty}I\right)  \in\mathcal{T}^{\otimes n-1}\otimes\mathcal{K}%
^{+}\subset\mathcal{T}^{\otimes n}$ becomes
\[
\left(  \boxplus^{d-r-\left[  R\right]  }0\right)  \boxplus\left(
\boxplus^{r+\left[  R\right]  }I\right)  \boxplus\left(  \boxplus^{\infty
}I\right)  =\tilde{I}-I\otimes P_{d-r-\left[  R\right]  }\in\tilde{I}+\left(
\mathcal{T}^{\otimes n-1}\otimes\mathcal{K}\right)  \subset\mathcal{T}%
^{\otimes n}.
\]

So we get
\[
P\sim\left(  \tilde{I}-I\otimes P_{d-r-\left[  R\right]  }\right)
\boxplus\left(  \boxplus^{m-1}\tilde{I}\right)  \boxplus\left(  \boxplus
^{2k+1-m}0\right)  ,
\]
the latter of which when conjugated by $U_{d-r-\left[  R\right]  }^{-1}$
yields $\tilde{I}\boxplus\left(  \boxplus^{m-1}\tilde{I}\right)
\boxplus\left(  \boxplus^{2k+1-m}0\right)  $, where $U_{d-r-\left[  R\right]
}$ is defined as $U_{d-r}$ by replacing $d-r$ by $d-r-\left[  R\right]  $.
Thus we get $P\sim\left(  \boxplus^{m}\tilde{I}\right)  \boxplus\left(
\boxplus^{2k+1-m}0\right)  \equiv\boxplus^{m}\tilde{I}$.

$\square$

\textbf{Corollary 3}. $\mathfrak{P}_{m}\left(  \mathcal{T}^{\otimes n}\right)
=\mathfrak{P}_{m}^{\prime}\left(  \mathcal{T}^{\otimes n}\right)
\equiv\left\{  \boxplus^{m}\tilde{I}\right\}  $ for all $m\geq\left\lfloor
\frac{n-1}{2}\right\rfloor +3$ and any $n\in\mathbb{N}$, where $\tilde{I}$ is
the identity element of $\mathcal{T}^{\otimes n}$.

Proof. We prove by induction on $n\in\mathbb{N}$. For $n=1$, we already know
that $\mathfrak{P}_{m}^{\prime}\left(  \mathcal{T}^{\otimes n}\right)
\equiv\mathfrak{P}_{m}\left(  \mathcal{T}^{\otimes n}\right)  $ for all $m>0$.

Now assume as the induction hypothesis that $\mathfrak{P}_{m}^{\prime}\left(
\mathcal{T}^{\otimes n}\right)  =\mathfrak{P}_{m}\left(  \mathcal{T}^{\otimes
n}\right)  $ for all $m\geq\left\lfloor \frac{n-1}{2}\right\rfloor +3$ for an
$n\in\mathbb{N}$.

Since we know that $GL_{m}\left(  \mathcal{T}^{\otimes n}\right)  $ is
connected for all $m\geq\left\lfloor \frac{n}{2}\right\rfloor +3$, the above
theorem implies that $\mathfrak{P}_{m}^{\prime}\left(  \mathcal{T}^{\otimes
n+1}\right)  =\mathfrak{P}_{m}\left(  \mathcal{T}^{\otimes n+1}\right)  $ for
all $m\geq\left\lfloor \frac{n}{2}\right\rfloor +3$.

$\square$

It remains open the problem of classification of low-rank idempotents over
$\mathcal{T}^{\otimes n}$. In particular, it is not clear whether there are
idempotents of non-standard (equivalence) type.

\section{Projective modules over $C\left(  \mathbb{S}_{H}^{2n-1}\right)  $}

Most of the arguments and results in the above study of projective modules
over $\mathcal{T}^{\otimes n}$ can be adapted to the case of the quantum
spheres $C\left(  \mathbb{S}_{H}^{2n-1}\right)  $.

Let $\partial_{n}:\mathcal{T}^{\otimes n}\rightarrow C\left(  \mathbb{S}%
_{H}^{2n-1}\right)  $ be the canonical quotient map by restricting the
groupoid $\mathfrak{T}_{n}$ to the closed invariant set $\overline
{\mathbb{Z}_{\geq}}^{n}\backslash\mathbb{Z}_{\geq}^{n}$ in its unit space.

First we note that there is a short exact sequence of C*-algebras%
\[
0\rightarrow C\left(  \mathbb{S}_{H}^{2n-3}\right)  \otimes\mathcal{K}\left(
\ell^{2}\left(  \mathbb{Z}_{\geq}\right)  \right)  \rightarrow C\left(
\mathbb{S}_{H}^{2n-1}\right)  \overset{\lambda_{n}}{\rightarrow}C^{\ast
}\left(  \left.  \mathfrak{T}_{n}\right\vert _{\overline{\mathbb{Z}_{\geq}%
}^{n-1}\times\left\{  \infty\right\}  }\right)  \cong\mathcal{T}^{\otimes
n-1}\otimes C\left(  \mathbb{T}\right)  \rightarrow0
\]
for all $n>1$. Indeed, since $\left(  \overline{\mathbb{Z}_{\geq}}%
^{n-1}\backslash\mathbb{Z}_{\geq}^{n-1}\right)  \times\mathbb{Z}_{\geq}$ is an
open invariant subset of the unit space $\overline{\mathbb{Z}_{\geq}}%
^{n}\backslash\mathbb{Z}_{\geq}^{n}$ of the groupoid $\mathfrak{G}_{n}%
\equiv\left.  \left(  \mathbb{Z}^{n}\ltimes\overline{\mathbb{Z}}^{n}\right)
\right\vert _{\overline{\mathbb{Z}_{\geq}}^{n}\backslash\mathbb{Z}_{\geq}^{n}%
}$ with the invariant complement
\[
\left(  \overline{\mathbb{Z}_{\geq}}^{n}\backslash\mathbb{Z}_{\geq}%
^{n}\right)  \backslash\left(  \left(  \overline{\mathbb{Z}_{\geq}}%
^{n-1}\backslash\mathbb{Z}_{\geq}^{n-1}\right)  \times\mathbb{Z}_{\geq
}\right)  =\overline{\mathbb{Z}_{\geq}}^{n-1}\times\left\{  \infty\right\}  ,
\]
the groupoid C*-algebra
\[
C^{\ast}\left(  \left.  \mathfrak{G}_{n}\right\vert _{\left(  \overline
{\mathbb{Z}_{\geq}}^{n-1}\backslash\mathbb{Z}_{\geq}^{n-1}\right)
\times\mathbb{Z}_{\geq}}\right)  =C^{\ast}\left(  \left.  \left(
\mathbb{Z}^{n-1}\ltimes\overline{\mathbb{Z}}^{n-1}\right)  \right\vert
_{\overline{\mathbb{Z}_{\geq}}^{n-1}\backslash\mathbb{Z}_{\geq}^{n-1}}%
\times\left.  \left(  \mathbb{Z}\ltimes\mathbb{Z}\right)  \right\vert
_{\mathbb{Z}_{\geq}}\right)
\]%
\[
\cong C^{\ast}\left(  \left.  \left(  \mathbb{Z}^{n-1}\ltimes\overline
{\mathbb{Z}}^{n-1}\right)  \right\vert _{\overline{\mathbb{Z}_{\geq}}%
^{n-1}\backslash\mathbb{Z}_{\geq}^{n-1}}\right)  \otimes C^{\ast}\left(
\left.  \left(  \mathbb{Z}\ltimes\mathbb{Z}\right)  \right\vert _{\mathbb{Z}%
_{\geq}}\right)  =C\left(  \mathbb{S}_{H}^{2n-3}\right)  \otimes
\mathcal{K}\left(  \ell^{2}\left(  \mathbb{Z}_{\geq}\right)  \right)
\]
is a closed ideal of $C^{\ast}\left(  \mathfrak{G}_{n}\right)  =C\left(
\mathbb{S}_{H}^{2n-1}\right)  $ with quotient
\[
C^{\ast}\left(  \mathfrak{G}_{n}\right)  /C^{\ast}\left(  \left.
\mathfrak{G}_{n}\right\vert _{\left(  \overline{\mathbb{Z}_{\geq}}%
^{n-1}\backslash\mathbb{Z}_{\geq}^{n-1}\right)  \times\mathbb{Z}_{\geq}%
}\right)  \cong C^{\ast}\left(  \left.  \mathfrak{G}_{n}\right\vert
_{\overline{\mathbb{Z}_{\geq}}^{n-1}\times\left\{  \infty\right\}  }\right)
\]%
\[
=C^{\ast}\left(  \left.  \left(  \mathbb{Z}^{n-1}\ltimes\overline{\mathbb{Z}%
}^{n-1}\right)  \right\vert _{\overline{\mathbb{Z}_{\geq}}^{n-1}}%
\times\mathbb{Z}\right)  \cong C^{\ast}\left(  \left.  \left(  \mathbb{Z}%
^{n-1}\ltimes\overline{\mathbb{Z}}^{n-1}\right)  \right\vert _{\overline
{\mathbb{Z}_{\geq}}^{n-1}}\right)  \otimes C\left(  \mathbb{T}\right)
=\mathcal{T}^{\otimes n-1}\otimes C\left(  \mathbb{T}\right)  .
\]
So we get the above short exact sequence with $\lambda_{n}$ being the
canonical map from $C^{\ast}\left(  \mathfrak{G}_{n}\right)  $ to its quotient
$\mathcal{T}^{\otimes n-1}\otimes C\left(  \mathbb{T}\right)  $ resulting from
restricting the groupoid $\mathfrak{G}_{n}$ to the closed invariant set
$\overline{\mathbb{Z}_{\geq}}^{n-1}\times\left\{  \infty\right\}  $.

Clearly $\kappa_{n}=\lambda_{n}\circ\partial_{n}$. Furthermore all the
quotient maps $\sigma_{\left(  j,\Theta\right)  }$ on $\mathcal{T}^{\otimes
n}$ with $j>0$ factors through $\partial_{n}$ and induces a quotient map
\[
\tau_{\left(  j,\Theta\right)  }:C\left(  \mathbb{S}_{H}^{2n-1}\right)
\rightarrow C\left(  \mathbb{T}^{j}\right)  \otimes\mathcal{B}\left(  \ell
^{2}\left(  \mathbb{Z}_{\geq}^{n-j}\right)  \right)
\]
such that $\sigma_{\left(  j,\Theta\right)  }=\tau_{\left(  j,\Theta\right)
}\circ\partial_{n}$.

Note that the quotient maps $\lambda_{n}$ for $n\in\mathbb{N}$ satisfy the
commuting diagram
\[%
\begin{array}
[c]{ccccc}%
M_{k}\left(  C\left(  \mathbb{S}_{H}^{2n-1}\right)  \right)  &
\overset{\lambda_{n}}{\rightarrow} & M_{k}\left(  \mathcal{T}^{\otimes
n-1}\otimes C\left(  \mathbb{T}\right)  \right)  & \equiv & M_{k}\left(
\mathcal{T}^{\otimes n-1}\right)  \otimes C\left(  \mathbb{T}\right) \\
_{{}}\downarrow_{\tau_{n}} & \circlearrowright & _{{}}\downarrow_{\sigma
_{n-1}\otimes\operatorname{id}} &  & \downarrow_{\sigma_{n-1}\otimes
\operatorname{id}}\\
M_{k}\left(  C\left(  \mathbb{T}^{n}\right)  \right)  & \overset{\equiv
}{\rightarrow} & M_{k}\left(  C\left(  \mathbb{T}^{n-1}\right)  \otimes
C\left(  \mathbb{T}\right)  \right)  & \equiv & M_{k}\left(  C\left(
\mathbb{T}^{n-1}\right)  \right)  \otimes C\left(  \mathbb{T}\right)  .
\end{array}
\]

We define the rank of (the equivalence class of) an idempotent $Q\in
M_{\infty}\left(  C\left(  \mathbb{S}_{H}^{2n-1}\right)  \right)  $ over
$C\left(  \mathbb{S}_{H}^{2n-1}\right)  $ as the rank of the matrix value
$\tau_{n}\left(  Q\right)  \left(  z\right)  \in M_{\infty}\left(
\mathbb{C}\right)  $ at any $z\in\mathbb{T}^{n}$ (independent of $z$ since
$\mathbb{T}^{n}$ is connected). Then the set of equivalence classes of
idempotents $Q\in M_{\infty}\left(  C\left(  \mathbb{S}_{H}^{2n-1}\right)
\right)  $ equipped with the binary operation $\boxplus$ becomes an abelian
graded monoid
\[
\mathfrak{P}\left(  C\left(  \mathbb{S}_{H}^{2n-1}\right)  \right)
=\sqcup_{m=0}^{\infty}\mathfrak{P}_{m}\left(  C\left(  \mathbb{S}_{H}%
^{2n-1}\right)  \right)
\]
where $\mathfrak{P}_{m}\left(  C\left(  \mathbb{S}_{H}^{2n-1}\right)  \right)
$ is the set of all (equivalence classes of) idempotents over $C\left(
\mathbb{S}_{H}^{2n-1}\right)  $ of rank $m$, with clearly
\[
\mathfrak{P}_{m}\left(  C\left(  \mathbb{S}_{H}^{2n-1}\right)  \right)
\boxplus\mathfrak{P}_{l}\left(  C\left(  \mathbb{S}_{H}^{2n-1}\right)
\right)  \subset\mathfrak{P}_{m+l}\left(  C\left(  \mathbb{S}_{H}%
^{2n-1}\right)  \right)
\]
for $m,l\geq0$.

Since $\sigma_{n}=\tau_{n}\circ\partial_{n}$, the rank of an idempotent $P$
over $C\left(  \mathcal{T}^{\otimes n}\right)  $ equals the rank of the
idempotent $\partial_{n}P$ over $C\left(  \mathbb{S}_{H}^{2n-1}\right)  $. We
now define
\[
\mathfrak{P}_{m}^{\prime}\left(  C\left(  \mathbb{S}_{H}^{2n-1}\right)
\right)  :=\partial_{n}\left(  \mathfrak{P}_{m}^{\prime}\left(  \mathcal{T}%
^{\otimes n}\right)  \right)  \subset\mathfrak{P}_{m}\left(  C\left(
\mathbb{S}_{H}^{2n-1}\right)  \right)  ,
\]
and the projections
\[
Q_{j,\Theta,l}:=\partial_{n}\left(  \Theta\left(  P_{j,l}\right)  \right)
\]
over $C\left(  \mathbb{S}_{H}^{2n-1}\right)  $. Note that $\mathfrak{P}%
_{m}^{\prime}\left(  C\left(  \mathbb{S}_{H}^{2n-1}\right)  \right)  =\left\{
\boxplus^{m}\tilde{I}\right\}  $ for $m>0$, where $\tilde{I}$ denotes the
identity element of $C\left(  \mathbb{S}_{H}^{2n-1}\right)  $.

Also note that $Q_{0,\operatorname{id},l}=0$ for all $l$, where
$\operatorname{id}\equiv\operatorname{id}_{\left\{  1,2,...,n\right\}  }$ is
the only $\left(  0,n\right)  $-shuffle. The monoid homomorphism
\[
\rho_{0}:P\in\mathfrak{P}^{\prime}\left(  \mathcal{T}^{\otimes n}\right)
\mapsto\prod_{\left(  j,\Theta\right)  \in\Omega_{0}}\rho_{\left(
j,\Theta\right)  }\left(  P\right)  \in\prod_{\left(  j,\Theta\right)
\in\Omega_{0}}\overline{\mathbb{Z}_{\geq}}\ ,
\]
with
\[
\Omega_{0}:=\Omega\backslash\left\{  \left(  0,\operatorname{id}\right)
\right\}  \equiv\left\{  \left(  j,\Theta\right)  :\ 0<j\leq n\text{\ and}%
\ \Theta\text{\ is a }\left(  j,n-j\right)  \text{-shuffle}\right\}  ,
\]
\textquotedblleft truncated\textquotedblright\ from $\rho$ induces a
well-defined monoid homomorphism
\[
\rho_{\partial}:\mathfrak{P}^{\prime}\left(  C\left(  \mathbb{S}_{H}%
^{2n-1}\right)  \right)  \rightarrow\prod_{\left(  j,\Theta\right)  \in
\Omega_{0}}\overline{\mathbb{Z}_{\geq}},
\]
in the sense that $\rho=\rho_{\partial}\circ\partial_{n}$. Indeed for $\left(
j,\Theta\right)  \in\Omega_{0}$, i.e. with $j>0$, the quotient map
\[
\sigma_{\left(  j,\Theta\right)  }:\mathcal{T}^{\otimes n}\equiv C^{\ast
}\left(  \mathfrak{T}_{n}\right)  \rightarrow C^{\ast}\left(  \left.
\mathfrak{T}_{n}\right\vert _{X_{\Theta}}\right)
\]
factors through $\partial_{n}$ since the unit space $\overline{\mathbb{Z}%
_{\geq}^{n}}\backslash\mathbb{Z}_{\geq}^{n}$ of $\mathfrak{G}_{n}$ contains
$X_{\Theta}$, and hence the map $\rho_{\left(  j,\Theta\right)  }$ factors
through $\partial_{n}$.

We call a $\boxplus$-sum of $Q_{j,\Theta,l}$ indexed by $\prec$-unrelated
$\left(  j,\Theta\right)  \in\Omega_{0}$ (i.e. $1\leq j\leq n$) and $l\equiv
l_{\left(  j,\Theta\right)  }>0$ depending on $\left(  j,\Theta\right)  $ to
be a reduced $\boxplus$-sum of standard projections over $C\left(
\mathbb{S}_{H}^{2n-1}\right)  $. (The degenerate empty $\boxplus$-sum $0$ is
taken as a reduced $\boxplus$-sum.) Two such reduced $\boxplus$-sums are
called different when they have different sets of (mutually $\prec$-unrelated)
indices $\left(  j,\Theta\right)  \in\Omega_{0}$ or have different weight
functions $l$ of $\left(  j,\Theta\right)  $. Each $Q_{j,\Theta,l}$ with
$j,l>0 $ is a reduced $\boxplus$-sum of standard projections over $C\left(
\mathbb{S}_{H}^{2n-1}\right)  $.

\textbf{Proposition 5}. Different reduced $\boxplus$-sums of standard
projections over $C\left(  \mathbb{S}_{H}^{2n-1}\right)  $ are mutually
inequivalent projections over $C\left(  \mathbb{S}_{H}^{2n-1}\right)  $, and
they form a graded submonoid
\[
\mathfrak{P}^{\prime}\left(  C\left(  \mathbb{S}_{H}^{2n-1}\right)  \right)
=\sqcup_{m=0}^{\infty}\mathfrak{P}_{m}^{\prime}\left(  C\left(  \mathbb{S}%
_{H}^{2n-1}\right)  \right)
\]
of the monoid $\mathfrak{P}\left(  C\left(  \mathbb{S}_{H}^{2n-1}\right)
\right)  $, with its monoid structure explicitly determined by $Q_{j,\Theta
,l}\boxplus Q_{j^{\prime},\Theta^{\prime},l^{\prime}}\sim Q_{j,\Theta,l} $ for
$\left(  j^{\prime},\Theta^{\prime}\right)  \prec\left(  j,\Theta\right)  $
with $j,j^{\prime},l,l^{\prime}>0$. Furthermore the monoid homomorphism
\[
\rho_{\partial}:\mathfrak{P}^{\prime}\left(  C\left(  \mathbb{S}_{H}%
^{2n-1}\right)  \right)  \rightarrow\prod_{\left(  j,\Theta\right)  \in
\Omega_{0}}\overline{\mathbb{Z}_{\geq}}%
\]
is injective.

Proof. The submonoid $\mathfrak{P}^{\prime}\left(  C\left(  \mathbb{S}%
_{H}^{2n-1}\right)  \right)  =\partial_{n}\left(  \mathfrak{P}^{\prime}\left(
C\left(  \mathcal{T}^{\otimes n}\right)  \right)  \right)  $ consists of
reduced $\boxplus$-sums of $Q_{j,\Theta,l}=\partial_{n}\left(  \Theta\left(
P_{j,l}\right)  \right)  $ with $j>0$, since $Q_{0,\operatorname{id},l}=0$.

Let $\mathfrak{M}$ be the subset of $\mathfrak{P}^{\prime}\left(  C\left(
\mathcal{T}^{\otimes n}\right)  \right)  $ consisting of all reduced
$\boxplus$-sums $P$ of $\Theta\left(  P_{j,l}\right)  $ with $j>0$. Then
$\partial_{n}|_{\mathfrak{M}}:\mathfrak{M}\rightarrow\mathfrak{P}^{\prime
}\left(  C\left(  \mathbb{S}_{H}^{2n-1}\right)  \right)  $ is still
surjective, and $\rho_{0}|_{\mathfrak{M}}$ still factors through
$\rho_{\partial}$, i.e. $\rho_{0}|_{\mathfrak{M}}=\rho_{\partial}\circ
\partial_{n}|_{\mathfrak{M}}$. These imply that $\rho_{\partial}$ is injective
if $\rho_{0}|_{\mathfrak{M}}$ is injective.

For any reduced $\boxplus$-sum $P\in\mathfrak{M}$ of $\Theta\left(
P_{j,l}\right)  $ with $j>0$, the $\left(  j,\Theta\right)  $-component of
$\rho\left(  P\right)  $ is the same as that of $\rho_{0}\left(  P\right)  $
for all $\left(  j,\Theta\right)  \in\Omega_{0}$, while the only other
component, namely, the $\left(  0,\operatorname{id}\right)  $-component of
$\rho\left(  P\right)  $ is $\infty$ since $\rho_{\left(  0,\operatorname{id}%
\right)  }\left(  \Theta\left(  P_{j,l}\right)  \right)  =\infty$ for any
$j>0$. Thus we get $\rho\left(  P\right)  =\left(  \infty,\rho_{0}\left(
P\right)  \right)  $ for all $P\in\mathfrak{M}$. Hence the injectivity of
$\rho|_{\mathfrak{M}}$ implies the injectivity of $\rho_{0}|_{\mathfrak{M}}$
on $\mathfrak{M}$, and hence the injectivity of $\rho_{\partial}$.

Since two different reduced $\boxplus$-sums $Q,Q^{\prime}$ over $C\left(
\mathbb{S}_{H}^{2n-1}\right)  $ are of the form $\partial_{n}\left(  P\right)
,\partial_{n}\left(  P^{\prime}\right)  $ respectively for two different
reduced $\boxplus$-sums $P,P^{\prime}\in\mathfrak{M}$ over $C\left(
\mathcal{T}^{\otimes n}\right)  $ which are inequivalent over $C\left(
\mathcal{T}^{\otimes n}\right)  $ and hence $\rho_{0}\left(  P\right)
\neq\rho_{0}\left(  P^{\prime}\right)  $, we get $\rho_{\partial}\left(
Q\right)  \neq\rho_{\partial}\left(  Q^{\prime}\right)  $ showing that
$Q,Q^{\prime}$ are different equivalence classes in $\mathfrak{P}^{\prime
}\left(  C\left(  \mathbb{S}_{H}^{2n-1}\right)  \right)  $.

The property that $\Theta\left(  P_{j,l}\right)  \boxplus\Theta^{\prime
}\left(  P_{j^{\prime},l^{\prime}}\right)  \sim\Theta\left(  P_{j,l}\right)  $
over $\mathcal{T}^{\otimes n}$ for $\left(  j^{\prime},\Theta^{\prime}\right)
\prec\left(  j,\Theta\right)  $ is clearly preserved under the quotient map
$\partial_{n}$, i.e. $Q_{j,\Theta,l}\boxplus Q_{j^{\prime},\Theta^{\prime
},l^{\prime}}\sim Q_{j,\Theta,l}$ over $C\left(  \mathbb{S}_{H}^{2n-1}\right)
$.

$\square$

\textbf{Theorem 3}. For $n>1$ and $m\in\mathbb{N}$, if $\mathfrak{P}%
_{m}\left(  \mathcal{T}^{\otimes n-1}\right)  =\mathfrak{P}_{m}^{\prime
}\left(  \mathcal{T}^{\otimes n-1}\right)  $ and $GL_{m}\left(  \mathcal{T}%
^{\otimes n-1}\right)  $ is connected, then $\mathfrak{P}_{m}^{\prime}\left(
C\left(  \mathbb{S}_{H}^{2n-1}\right)  \right)  =\mathfrak{P}_{m}\left(
C\left(  \mathbb{S}_{H}^{2n-1}\right)  \right)  $.

Proof. Many arguments used to prove a similar theorem for $\mathcal{T}%
^{\otimes n}$ instead of $C\left(  \mathbb{S}_{H}^{2n-1}\right)  $ can be used
again here with minor modifications. In this proof, $I$ and $\tilde{I}$ denote
respectively the identity element of $\mathcal{T}^{\otimes n-1}$ and $C\left(
\mathbb{S}_{H}^{2n-1}\right)  $.

Let $P\in\mathfrak{P}_{m}\left(  C\left(  \mathbb{S}_{H}^{2n-1}\right)
\right)  $. The idempotent $\lambda_{n}\left(  P\right)  $ over $\mathcal{T}%
^{\otimes n-1}\otimes C\left(  \mathbb{T}\right)  $ satisfies that for any
$z\in\mathbb{T}$,
\[
\sigma_{n-1}\left(  \lambda_{n}\left(  P\right)  \left(  z\right)  \right)
=\tau_{n}\left(  P\right)  \left(  \cdot,z\right)  \in M_{\infty}\left(
C\left(  \mathbb{T}^{n-1}\right)  \right)
\]
which is of rank $m$ pointwise, and hence
\[
\lambda_{n}\left(  P\right)  \left(  z\right)  \in\mathfrak{P}_{m}\left(
\mathcal{T}^{\otimes n-1}\right)  =\mathfrak{P}_{m}^{\prime}\left(
\mathcal{T}^{\otimes n-1}\right)  ,
\]
i.e. $\lambda_{n}\left(  P\right)  \left(  z\right)  \sim\boxplus^{m}I$ over
$\mathcal{T}^{\otimes n-1}$. As before, for some large $k$, there is an
idempotent-valued continuous loop $\Gamma:\mathbb{T}\rightarrow M_{k}\left(
\mathcal{T}^{\otimes n-1}\right)  $ starting and ending at $\boxplus^{m}I$
with $\Gamma\left(  e^{i\theta}\right)  =\left(  \boxplus^{m}I\right)
\boxplus\left(  \boxplus^{k-m}0\right)  $, say, for all $\theta\in\left[
3\pi/2,2\pi\right]  $, and homotopic to the loop $\lambda_{n}\left(  P\right)
$ via idempotents. Consequently, there is a continuous path of invertibles
$U_{t}\in GL_{k}\left(  \mathcal{T}^{\otimes n-1}\otimes C\left(
\mathbb{T}\right)  \right)  $ with $U_{0}=I_{k}$ such that $U_{1}\lambda
_{n}\left(  P\right)  U_{1}^{-1}=\Gamma$, which can be lifted along
$\lambda_{n}$ to a continuous path of invertible $V_{t}\in GL_{k}\left(
C\left(  \mathbb{S}_{H}^{2n-1}\right)  \right)  $ with $V_{0}=I_{k}$ such that
$\lambda_{n}\left(  V_{1}PV_{1}^{-1}\right)  =\Gamma$.

Replacing $P$ by the equivalent idempotent $V_{1}PV_{1}^{-1}$, we may now
assume directly that the idempotent $\lambda_{n}\left(  P\right)  $ over
$\mathcal{T}^{\otimes n-1}\otimes C\left(  \mathbb{T}\right)  $ is a
continuous loop of idempotents in $M_{k}\left(  \mathcal{T}^{\otimes
n-1}\right)  $ such that $\lambda_{n}\left(  P\right)  \left(  e^{i\theta
}\right)  =\left(  \boxplus^{m}I\right)  \boxplus\left(  \boxplus
^{k-m}0\right)  $ for all $\theta\in\left[  3\pi/2,2\pi\right]  $. As before,
by the connectedness assumption on $GL_{m}\left(  \mathcal{T}^{\otimes
n-1}\right)  $, after suitably increasing the size $k$, we can find a
well-defined continuous loop
\[
W:e^{i\theta}\in\mathbb{T}\mapsto W_{\theta}\in GL_{k}\left(  \mathcal{T}%
^{\otimes n-1}\right)  ,
\]
i.e. $W\in GL_{k}\left(  \mathcal{T}^{\otimes n-1}\otimes C\left(
\mathbb{T}\right)  \right)  $, satisfying
\[
W\left(  \lambda_{n}\left(  P\right)  \right)  W^{-1}=\left(  \boxplus
^{m}I\right)  \boxplus\left(  \boxplus^{k-m}0\right)  .
\]
So the idempotent $\lambda_{n}\left(  P\right)  $ over $\mathcal{T}^{\otimes
n-1}\otimes C\left(  \mathbb{T}\right)  $ is equivalent to the idempotent
$\boxplus^{m}I$.

Replacing $P$ by the equivalent idempotent $\tilde{W}\left(  P\boxplus\left(
\boxplus^{k}0\right)  \right)  \tilde{W}^{-1}$ for any fixed lifting
$\tilde{W}\in GL_{2k}^{0}\left(  C\left(  \mathbb{S}_{H}^{2n-1}\right)
\right)  $ of $W\boxplus W^{-1}\in GL_{2k}^{0}\left(  \mathcal{T}^{\otimes
n-1}\otimes C\left(  \mathbb{T}\right)  \right)  $ along $\lambda_{n}$, we may
now assume that
\[
\lambda_{n}\left(  P\right)  =\left(  \boxplus^{m}I\right)  \boxplus\left(
\boxplus^{2k-m}0\right)  =\lambda_{n}\left(  \left(  \boxplus^{m}\tilde
{I}\right)  \boxplus\left(  \boxplus^{2k-m}0\right)  \right)
\]
and proceed to show that $P\sim\boxplus^{m}\tilde{I}$ over $C\left(
\mathbb{S}_{H}^{2n-1}\right)  $, where we use $\tilde{I}$ to denote the
identity element in $C\left(  \mathbb{S}_{H}^{2n-1}\right)  $ so as to
distinguish it from the identity element $I$ of $\mathcal{T}^{\otimes n-1}$.

With $P-\left(  \left(  \boxplus^{m}\tilde{I}\right)  \boxplus\left(
\boxplus^{2k-m}0\right)  \right)  \in M_{2k}\left(  C\left(  \mathbb{S}%
_{H}^{2n-3}\right)  \otimes\mathcal{K}\right)  $ and $M_{\infty}\left(
\mathbb{C}\right)  $ dense in $\mathcal{K}$, we may replace $P$ by a suitable
equivalent idempotent and assume that
\[
P=K+\left(  \left(  \boxplus^{m}\tilde{I}\right)  \boxplus\left(
\boxplus^{2k-m}0\right)  \right)  \in M_{2k}\left(  \left(  C\left(
\mathbb{S}_{H}^{2n-3}\right)  \otimes\mathcal{K}\right)  ^{+}\right)  \subset
M_{2k}\left(  C\left(  \mathbb{S}_{H}^{2n-1}\right)  \right)  \text{\ }%
\]
for some $K\in M_{2k}\left(  C\left(  \mathbb{S}_{H}^{2n-3}\right)  \otimes
M_{N}\left(  \mathbb{C}\right)  \right)  $\ and some $N\in\mathbb{N}$.

As before, by rearranging entries via conjugation, we get%
\[
P\sim\partial_{n}\left(  U_{k,N}\right)  P\partial_{n}\left(  U_{k,N}%
^{-1}\right)  \equiv\partial_{n}\left(  U_{k,N}\right)  \left(  P\boxplus
0\right)  \partial_{n}\left(  U_{k,N}^{-1}\right)  =\left(  \left(
\boxplus^{m}\tilde{I}\right)  \boxplus\left(  \boxplus^{2k-m}0\right)
\right)  \boxplus R
\]
for some
\[
R\in M_{2kN}\left(  C\left(  \mathbb{S}_{H}^{2n-3}\right)  \right)  \equiv
C\left(  \mathbb{S}_{H}^{2n-3}\right)  \otimes M_{2kN}\left(  \mathbb{C}%
\right)  \subset\left(  C\left(  \mathbb{S}_{H}^{2n-3}\right)  \otimes
\mathcal{K}\right)  ^{+}\subset C\left(  \mathbb{S}_{H}^{2n-1}\right)
\]
which must be an idempotent over $C\left(  \mathbb{S}_{H}^{2n-3}\right)  $.
More precisely, we can lift $P$ to
\[
\hat{P}=\hat{K}+\left(  \left(  \boxplus^{m}I_{\mathcal{T}^{\otimes n}%
}\right)  \boxplus\left(  \boxplus^{2k-m}0\right)  \right)  \in M_{2k}\left(
\left(  \mathcal{T}^{\otimes n-1}\otimes\mathcal{K}\right)  ^{+}\right)
\subset M_{2k}\left(  \mathcal{T}^{\otimes n}\right)
\]
for some $\hat{K}\in M_{2k}\left(  \mathcal{T}^{\otimes n-1}\otimes
M_{N}\left(  \mathbb{C}\right)  \right)  $ and conjugate it by the unitary
$U_{k,N}$ over $\mathcal{T}^{\otimes n}$ to get the form $\left(  \left(
\boxplus^{m}I_{\mathcal{T}^{\otimes n}}\right)  \boxplus\left(  \boxplus
^{2k-m}0\right)  \right)  \boxplus\hat{R}$ with $\hat{R}\in M_{2kN}\left(
\mathcal{T}^{\otimes n-1}\right)  $ as we did for the case of $\mathcal{T}%
^{\otimes n}$. Then the above $R$ is $\partial_{n}\left(  \hat{R}\right)  $.
Note that even though $\hat{P}$ and $\hat{R}$ are not necessarily idempotents,
$R$ is since it is the idempotent $P$ conjugated by the unitary $\partial
_{n}\left(  U_{k,N}\right)  $ over $C\left(  \mathbb{S}_{H}^{2n-1}\right)  $.

Since $K_{0}\left(  C\left(  \mathbb{S}_{H}^{2n-3}\right)  \right)
=\mathbb{Z}$ \cite{HaNePaSiZi}, $R\boxplus\left(  \boxplus^{r}\hat{I}\right)
\sim\left(  \boxplus^{r+\left[  R\right]  }\hat{I}\right)  $ for a
sufficiently large $r\in\mathbb{N}$ where $\left[  R\right]  \in\mathbb{Z}$
denotes the class of $R$ in $K_{0}\left(  C\left(  \mathbb{S}_{H}%
^{2n-3}\right)  \right)  $ and $\hat{I}$ is the identity element of $C\left(
\mathbb{S}_{H}^{2n-3}\right)  $. So there is an invertible $U\in GL_{d}\left(
C\left(  \mathbb{S}_{H}^{2n-3}\right)  \right)  $ for some large $d\geq
\max\left\{  2kN+r,r+\left[  R\right]  \right\}  $ such that
\[
U\left(  R\boxplus\left(  \boxplus^{d-2kN-r}0\right)  \boxplus\left(
\boxplus^{r}\hat{I}\right)  \right)  U^{-1}=\left(  \boxplus^{d-r-\left[
R\right]  }0\right)  \boxplus\left(  \boxplus^{r+\left[  R\right]  }\hat
{I}\right)  .
\]

As before, with $m>0$, by rearranging entries via conjugation, we can get
\[
P\sim R^{\prime}\boxplus\left(  \boxplus^{m-1}\tilde{I}\right)  \boxplus
\left(  \boxplus^{2k+1-m}0\right)
\]
where the idempotent
\[
R^{\prime}=\left(  R\boxplus\left(  \boxplus^{d-2kN-r}0\right)  \right)
+\left(  \tilde{I}-\hat{I}\otimes P_{d-r}\right)  \in\left(  C\left(
\mathbb{S}_{H}^{2n-3}\right)  \otimes\mathcal{K}\right)  ^{+}\subset C\left(
\mathbb{S}_{H}^{2n-1}\right)
\]
when conjugated by the invertible $U\equiv U\boxplus\left(  \tilde{I}-\hat
{I}\otimes P_{d}\right)  \in\left(  C\left(  \mathbb{S}_{H}^{2n-3}\right)
\otimes\mathcal{K}\right)  ^{+}$ becomes
\[
\left(  \boxplus^{d-r-\left[  R\right]  }0\right)  \boxplus\left(  \tilde
{I}-\hat{I}\otimes P_{d-r-\left[  R\right]  }\right)  \in\left(  C\left(
\mathbb{S}_{H}^{2n-3}\right)  \otimes\mathcal{K}\right)  ^{+}\subset C\left(
\mathbb{S}_{H}^{2n-1}\right)  .
\]

So we get
\[
P\sim\left(  \left(  \boxplus^{d-r-\left[  R\right]  }0\right)  \boxplus
\left(  \tilde{I}-\hat{I}\otimes P_{d-r-\left[  R\right]  }\right)  \right)
\boxplus\left(  \boxplus^{m-1}\tilde{I}\right)  \boxplus\left(  \boxplus
^{2k+1-m}0\right)  ,
\]
the latter of which as before is equivalent to $\tilde{I}\boxplus\left(
\boxplus^{m-1}\tilde{I}\right)  \boxplus\left(  \boxplus^{2k+1-m}0\right)  $
by a further conjugation by $U_{d-r-\left[  R\right]  }^{-1}$. Thus
$P\sim\left(  \boxplus^{m}\tilde{I}\right)  \boxplus\left(  \boxplus
^{2k+1-m}0\right)  \equiv\boxplus^{m}\tilde{I}$.

$\square$

\textbf{Corollary 4}. $\mathfrak{P}_{m}\left(  C\left(  \mathbb{S}_{H}%
^{2n-1}\right)  \right)  =\mathfrak{P}_{m}^{\prime}\left(  C\left(
\mathbb{S}_{H}^{2n-1}\right)  \right)  \equiv\left\{  \boxplus^{m}\tilde
{I}\right\}  $ for all $m\geq\left\lfloor \frac{n-1}{2}\right\rfloor +3$ and
any $n\in\mathbb{N}$, where $\tilde{I}$ is the identity element of $C\left(
\mathbb{S}_{H}^{2n-1}\right)  $.

Proof. The case of $n=1$ is well known. For $n>1$, since $\mathfrak{P}%
_{m}^{\prime}\left(  \mathcal{T}^{\otimes n-1}\right)  =\mathfrak{P}%
_{m}\left(  \mathcal{T}^{\otimes n-1}\right)  $ for all $m\geq\left\lfloor
\frac{n-2}{2}\right\rfloor +3$ and $GL_{m}\left(  \mathcal{T}^{\otimes
n-1}\right)  $ is connected for all $m\geq\left\lfloor \frac{n-1}%
{2}\right\rfloor +3$, the above theorem implies that $\mathfrak{P}_{m}%
^{\prime}\left(  C\left(  \mathbb{S}_{H}^{2n-1}\right)  \right)
=\mathfrak{P}_{m}\left(  C\left(  \mathbb{S}_{H}^{2n-1}\right)  \right)  $ for
all $m\geq\left\lfloor \frac{n-1}{2}\right\rfloor +3$.

$\square$

It is not clear whether there are (low-rank) idempotents over $C\left(
\mathbb{S}_{H}^{2n-1}\right)  $ of non-standard (equivalence) type and whether
the cancellation law holds for them.

\section{Projective modules over $C\left(  \mathbb{P}^{n-1}\left(
\mathcal{T}\right)  \right)  $}

In this section we study the problem of classification of finitely generated
projective modules over the multipullback quantum complex projective space
$\mathbb{P}^{n-1}\left(  \mathcal{T}\right)  $ that was introduced and studied
by Hajac, Kaygun, Zieli\'{n}ski in \cite{HaKaZi}.

In \cite{HaNePaSiZi}, $K_{0}\left(  C\left(  \mathbb{P}^{n-1}\left(
\mathcal{T}\right)  \right)  \right)  =\mathbb{Z}^{n}$ and $K_{1}\left(
C\left(  \mathbb{P}^{n-1}\left(  \mathcal{T}\right)  \right)  \right)  =0$ are
computed, and $\mathbb{P}^{n-1}\left(  \mathcal{T}\right)  $ is shown to be a
quantum quotient space of $\mathbb{S}_{H}^{2n-1}$. More precisely, the
C*-algebra $C\left(  \mathbb{P}^{n-1}\left(  \mathcal{T}\right)  \right)  $ is
isomorphic to the invariant C*-subalgebra $\left(  C\left(  \mathbb{S}%
_{H}^{2n-1}\right)  \right)  ^{U\left(  1\right)  }$ of $C\left(
\mathbb{S}_{H}^{2n-1}\right)  $ under the canonical diagonal $U\left(
1\right)  $-action on $C\left(  \mathbb{S}_{H}^{2n-1}\right)  \cong%
\mathcal{T}^{\otimes n}/\mathcal{K}^{\otimes n}$, which in the groupoid
context can be implemented by the multiplication operator
\[
U_{\zeta}:f\in C_{c}\left(  \mathfrak{G}_{n}\right)  \mapsto h_{\zeta}f\in
C_{c}\left(  \mathfrak{G}_{n}\right)
\]
for $\zeta\in U\left(  1\right)  \equiv\mathbb{T}$ where
\[
h_{\zeta}:\left(  m,p\right)  \in\mathfrak{G}_{n}\subset\mathbb{Z}^{n}%
\ltimes\overline{\mathbb{Z}}^{n}\mapsto\zeta^{\Sigma m}\in\mathbb{T}%
\text{\ \ with }\Sigma m:=\sum_{i=1}^{n}m_{i}%
\]
is a groupoid character. Then $C\left(  \mathbb{P}^{n-1}\left(  \mathcal{T}%
\right)  \right)  $ is realized as the groupoid C*-algebra $C^{\ast}\left(
\left(  \mathfrak{G}_{n}\right)  _{0}\right)  $ of the subgroupoid $\left(
\mathfrak{G}_{n}\right)  _{0}$ of $\mathfrak{G}_{n}$, where
\[
\left(  \mathfrak{G}_{n}\right)  _{k}:=\left\{  \left(  m,p\right)
\in\mathfrak{G}_{n}:\Sigma m=k\right\}
\]
for $k\in\mathbb{Z}$. Furthermore, $C^{\ast}\left(  \mathfrak{G}_{n}\right)  $
becomes a (completion of the) graded algebra $\oplus_{k\in\mathbb{Z}}%
\overline{C_{c}\left(  \left(  \mathfrak{G}_{n}\right)  _{k}\right)  }$ with
the component $\overline{C_{c}\left(  \left(  \mathfrak{G}_{n}\right)
_{k}\right)  } $ being the quantum line bundle $C\left(  \mathbb{S}_{H}%
^{2n-1}\right)  _{k}$ \cite{HaNePaSiZi} of degree $k$ over the quantum space
$\mathbb{P}^{n-1}\left(  \mathcal{T}\right)  $.

It is easy to see that the standard projections $Q_{j,\Theta,l}\equiv
\partial_{n}\left(  \Theta\left(  P_{j,l}\right)  \right)  $ over $C\left(
\mathbb{S}_{H}^{2n-1}\right)  $ with $j,l>0$ found in the previous section lie
in $M_{\infty}\left(  C^{\ast}\left(  \left(  \mathfrak{G}_{n}\right)
_{0}\right)  \right)  $ since $P_{j,l}=\boxplus^{l}\left(  \left(  \otimes
^{j}I\right)  \otimes\left(  \otimes^{n-j}P_{1}\right)  \right)  $ is in
$C^{\ast}\left(  \left(  \mathfrak{T}_{n}\right)  _{0}\right)  $, and hence
are also projections over $C^{\ast}\left(  \left(  \mathfrak{G}_{n}\right)
_{0}\right)  \equiv C\left(  \mathbb{P}^{n-1}\left(  \mathcal{T}\right)
\right)  $. Furthermore with $C\left(  \mathbb{P}^{n-1}\left(  \mathcal{T}%
\right)  \right)  \subset C\left(  \mathbb{S}_{H}^{2n-1}\right)  $,
inequivalent $\boxplus$-sums of standard projections over $C\left(
\mathbb{S}_{H}^{2n-1}\right)  $ must be inequivalent over $C\left(
\mathbb{P}^{n-1}\left(  \mathcal{T}\right)  \right)  $ as well.

\textbf{Proposition 6}. Different reduced $\boxplus$-sums of standard
projections $Q_{j,\Theta,l}$ over $C\left(  \mathbb{S}_{H}^{2n-1}\right)  $
with $j,l>0$ when viewed as projections over $C\left(  \mathbb{P}^{n-1}\left(
\mathcal{T}\right)  \right)  $ are mutually inequivalent over $C\left(
\mathbb{P}^{n-1}\left(  \mathcal{T}\right)  \right)  $, and they form a graded
submonoid
\[
\mathfrak{P}^{\prime}\left(  C\left(  \mathbb{P}^{n-1}\left(  \mathcal{T}%
\right)  \right)  \right)  =\sqcup_{m=0}^{\infty}\mathfrak{P}_{m}^{\prime
}\left(  C\left(  \mathbb{P}^{n-1}\left(  \mathcal{T}\right)  \right)
\right)
\]
of the monoid $\mathfrak{P}\left(  C\left(  \mathbb{P}^{n-1}\left(
\mathcal{T}\right)  \right)  \right)  $. Furthermore the monoid homomorphism
\[
\mathfrak{P}^{\prime}\left(  C\left(  \mathbb{P}^{n-1}\left(  \mathcal{T}%
\right)  \right)  \right)  \rightarrow\prod_{\left(  j,\Theta\right)
\in\Omega_{0}}\overline{\mathbb{Z}_{\geq}}%
\]
inherited from $\rho_{\partial}$ is injective.

However, for $\left(  j^{\prime},\Theta^{\prime}\right)  \prec\left(
j,\Theta\right)  $ with $j,j^{\prime},l,l^{\prime}>0$, it is no longer true in
general that $Q_{j,\Theta,l}\boxplus Q_{j^{\prime},\Theta^{\prime},l^{\prime}%
}\sim Q_{j,\Theta,l}$ over $C\left(  \mathbb{P}^{n-1}\left(  \mathcal{T}%
\right)  \right)  $, even though $Q_{j,\Theta,l}\boxplus Q_{j^{\prime}%
,\Theta^{\prime},l^{\prime}}\sim Q_{j,\Theta,l}$ over $C\left(  \mathbb{S}%
_{H}^{2n-1}\right)  $ since the invertible matrix over $C\left(
\mathbb{S}_{H}^{2n-1}\right)  $ intertwining $Q_{j,\Theta,l}\boxplus
Q_{j^{\prime},\Theta^{\prime},l^{\prime}}$ and $Q_{j,\Theta,l}$ may not be
replaced by one over the subalgebra $C\left(  \mathbb{P}^{n-1}\left(
\mathcal{T}\right)  \right)  $ of $C\left(  \mathbb{S}_{H}^{2n-1}\right)  $.

In the following, we show that the standard projections
$Q_{j,\operatorname{id},1}$ with $j>0$ provide a set of representatives of
$K_{0}$-classes that freely generate the abelian $K_{0}$-group of $C\left(
\mathbb{P}^{n-1}\left(  \mathcal{T}\right)  \right)  $.

The subgroupoid $\mathfrak{H}_{j}:=\mathfrak{G}_{j}\times\left(
\mathbb{Z}^{n-j}\ltimes\mathbb{Z}_{\geq}^{n-j}\right)  $ of $\mathfrak{G}_{n}$
for $1\leq j\leq n$ is the groupoid $\mathfrak{G}_{n}$ restricted to the open
invariant subset $\left(  \overline{\mathbb{Z}_{\geq}^{j}}\times
\mathbb{Z}_{\geq}^{n-j}\right)  \backslash\mathbb{Z}_{\geq}^{n}$ and inherits
the grading of $\mathfrak{G}_{n}$. The grade-$0$ part $\left(  \mathfrak{H}%
_{j}\right)  _{0}$ of $\mathfrak{H}_{j}$ is the groupoid $\left(
\mathfrak{G}_{n}\right)  _{0}$ restricted to $\left(  \overline{\mathbb{Z}%
_{\geq}^{j}}\times\mathbb{Z}_{\geq}^{n-j}\right)  \backslash\mathbb{Z}_{\geq
}^{n}$, and from the increasing chain of $\left(  \mathfrak{H}_{j}\right)
_{0}$, we get an increasing composition sequence of closed ideals of $C^{\ast
}\left(  \left(  \mathfrak{G}_{n}\right)  _{0}\right)  $ as
\[
0=:C^{\ast}\left(  \left(  \mathfrak{H}_{0}\right)  _{0}\right)
\vartriangleleft C^{\ast}\left(  \left(  \mathfrak{H}_{1}\right)  _{0}\right)
\vartriangleleft\cdots\vartriangleleft C^{\ast}\left(  \left(  \mathfrak{H}%
_{n-1}\right)  _{0}\right)  \vartriangleleft C^{\ast}\left(  \left(
\mathfrak{H}_{n}\right)  _{0}\right)  =C^{\ast}\left(  \left(  \mathfrak{G}%
_{n}\right)  _{0}\right)
\]
such that with $\left(  \overline{\mathbb{Z}_{\geq}^{j}}\times\mathbb{Z}%
_{\geq}^{n-j}\right)  \backslash\left(  \overline{\mathbb{Z}_{\geq}^{j-1}%
}\times\mathbb{Z}_{\geq}^{n-j+1}\right)  =\overline{\mathbb{Z}_{\geq}^{j-1}%
}\times\left\{  \infty\right\}  \times\mathbb{Z}_{\geq}^{n-j}$,
\[
C^{\ast}\left(  \left(  \mathfrak{H}_{j}\right)  _{0}\right)  /C^{\ast}\left(
\left(  \mathfrak{H}_{j-1}\right)  _{0}\right)  \cong C^{\ast}\left(  \left(
\left.  \mathfrak{G}_{n}\right\vert _{\overline{\mathbb{Z}_{\geq}^{j-1}}%
\times\left\{  \infty\right\}  \times\mathbb{Z}_{\geq}^{n-j}}\right)
_{0}\right)  \cong C^{\ast}\left(  \left.  \mathfrak{T}_{n-1}\right\vert
_{\overline{\mathbb{Z}_{\geq}^{j-1}}\times\mathbb{Z}_{\geq}^{n-j}}\right)
\cong\mathcal{T}^{\otimes j-1}\otimes\mathcal{K}\left(  \mathbb{Z}_{\geq
}^{n-j}\right)
\]
because the groupoid $\left(  \left.  \mathfrak{G}_{n}\right\vert
_{\overline{\mathbb{Z}_{\geq}^{j-1}}\times\left\{  \infty\right\}
\times\mathbb{Z}_{\geq}^{n-j}}\right)  _{0}$ is isomorphic to the groupoid
$\left.  \mathfrak{T}_{n-1}\right\vert _{\overline{\mathbb{Z}_{\geq}^{j-1}%
}\times\mathbb{Z}_{\geq}^{n-j}}$ via the groupoid isomorphism
\[
\left(  m,k,l,p,\infty,q\right)  \mapsto\left(  m,l,p,q\right)
\]
where
\[
\left(  m,k,l,p,\infty,q\right)  \in\left.  \mathfrak{G}_{n}\right\vert
_{\overline{\mathbb{Z}_{\geq}^{j-1}}\times\left\{  \infty\right\}
\times\mathbb{Z}_{\geq}^{n-j}}\subset\mathbb{Z}^{j-1}\times\mathbb{Z}%
\times\mathbb{Z}^{n-j}\times\overline{\mathbb{Z}_{\geq}^{j-1}}\times\left\{
\infty\right\}  \times\mathbb{Z}_{\geq}^{n-j}%
\]
with $\sum_{i=1}^{j-1}m_{i}+k+\sum_{i=1}^{n-j}l_{i}=0$ and hence $k=-\sum
m-\sum l$ determined by $m,l$.

Since $K_{1}\left(  \mathcal{T}^{\otimes j-1}\otimes\mathcal{K}\left(
\mathbb{Z}_{\geq}^{n-j}\right)  \right)  =0$ and $K_{0}\left(  \mathcal{T}%
^{\otimes j-1}\otimes\mathcal{K}\left(  \mathbb{Z}_{\geq}^{n-j}\right)
\right)  =\mathbb{Z}$, it is easy to conclude from the cyclic six-term exact
sequence of $K$-groups for the pair $C^{\ast}\left(  \left(  \mathfrak{H}%
_{j-1}\right)  _{0}\right)  \vartriangleleft C^{\ast}\left(  \left(
\mathfrak{H}_{j}\right)  _{0}\right)  $ that the following sequence is exact
and splits%
\[
0\rightarrow K_{0}\left(  C^{\ast}\left(  \left(  \mathfrak{H}_{j-1}\right)
_{0}\right)  \right)  \rightarrow K_{0}\left(  C^{\ast}\left(  \left(
\mathfrak{H}_{j}\right)  _{0}\right)  \right)  \rightarrow K_{0}\left(
\mathcal{T}^{\otimes j-1}\otimes\mathcal{K}\left(  \mathbb{Z}_{\geq}%
^{n-j}\right)  \right)  \cong\mathbb{Z}\rightarrow0
\]
where the projection $\left(  \otimes^{j-1}I\right)  \otimes\left(
\otimes^{n-j}P_{1}\right)  $ is a generator of $K_{0}\left(  \mathcal{T}%
^{\otimes j-1}\otimes\mathcal{K}\left(  \mathbb{Z}_{\geq}^{n-j}\right)
\right)  $. Note that this $\left(  \otimes^{j-1}I\right)  \otimes\left(
\otimes^{n-j}P_{1}\right)  $ lifts to the projection element $\chi_{A_{j}}\in
C_{c}\left(  \left(  \mathfrak{H}_{j}\right)  _{0}\right)  \subset C^{\ast
}\left(  \left(  \mathfrak{H}_{j}\right)  _{0}\right)  $ given by the
characteristic function of the set
\[
A_{j}:=\left\{  0\right\}  \times\left\{  0\right\}  \times\left(
\overline{\mathbb{Z}_{\geq}^{j}}\backslash\mathbb{Z}_{\geq}^{j}\right)
\times\left\{  0\right\}  \subset\mathfrak{H}_{j}\subset\mathbb{Z}^{j}%
\times\mathbb{Z}^{n-j}\times\overline{\mathbb{Z}_{\geq}^{j}}\times
\mathbb{Z}_{\geq}^{n-j}.
\]
Furthermore $\chi_{A_{j}}=Q_{j,\operatorname{id},1}$ in $C\left(
\mathbb{P}^{n-1}\left(  \mathcal{T}\right)  \right)  \subset C\left(
\mathbb{S}_{H}^{2n-1}\right)  $. So we get
\[
K_{0}\left(  C^{\ast}\left(  \left(  \mathfrak{H}_{j}\right)  _{0}\right)
\right)  \cong K_{0}\left(  C^{\ast}\left(  \left(  \mathfrak{H}_{j-1}\right)
_{0}\right)  \right)  \oplus\mathbb{Z}\left[  Q_{j,\Theta,1}\right]
\]
with $K_{0}\left(  C^{\ast}\left(  \left(  \mathfrak{H}_{j-1}\right)
_{0}\right)  \right)  $ canonically embedded in $K_{0}\left(  C^{\ast}\left(
\left(  \mathfrak{H}_{j}\right)  _{0}\right)  \right)  $.

Putting together these results for all $j$, we get
\[
K_{0}\left(  C\left(  \mathbb{P}^{n-1}\left(  \mathcal{T}\right)  \right)
\right)  \equiv K_{0}\left(  C^{\ast}\left(  \left(  \mathfrak{H}_{n}\right)
_{0}\right)  \right)  \cong\oplus_{j=1}^{n}\mathbb{Z}\left[
Q_{j,\operatorname{id},1}\right]  \cong\mathbb{Z}^{n}%
\]
and hence $Q_{j,\operatorname{id},1}$ freely generate the abelian group
$K_{0}\left(  C\left(  \mathbb{P}^{n-1}\left(  \mathcal{T}\right)  \right)
\right)  $. Note that $Q_{j,\operatorname{id},l}=\boxplus^{l}%
Q_{j,\operatorname{id},1}$ and $\left[  Q_{j,\operatorname{id},l}\right]
=l\left[  Q_{j,\operatorname{id},1}\right]  $ in $K_{0}\left(  C\left(
\mathbb{P}^{n-1}\left(  \mathcal{T}\right)  \right)  \right)  $ for any
$l\in\mathbb{N}$.

We now summarize the above discussion.

\textbf{Theorem 4}. The projections $Q_{j,\Theta,l}$ over $C\left(
\mathbb{P}^{n-1}\left(  \mathcal{T}\right)  \right)  $ with $l\in\mathbb{N}$
and $\Theta$ a $\left(  j,n-j\right)  $-shuffle for $0<j\leq n$ are mutually
inequivalent, and the projections $Q_{j,\operatorname{id},1}$ with $0<j\leq n$
freely generate the abelian group $K_{0}\left(  C\left(  \mathbb{P}%
^{n-1}\left(  \mathcal{T}\right)  \right)  \right)  $, such that if $\left[
p\right]  =\sum_{j=1}^{n}m_{j}\left[  Q_{j,\operatorname{id},1}\right]  $ for
a projection $p$ over $C\left(  \mathbb{P}^{n-1}\left(  \mathcal{T}\right)
\right)  $, then the coefficient $m_{n}$ of $\left[  Q_{n,\operatorname{id}%
,1}\right]  $ is the rank of $p$.

Proof. We only need to note that the rank of $Q_{n,\operatorname{id},1}$ is
$1$ and the rank of any other $Q_{j,\operatorname{id},1}$ is $0$. $\square$

\textbf{Remark}. Since any permutation $\Theta$ on $\left\{
1,2,...,n\right\}  $ canonically induces a $U\left(  1\right)  $-equivariant
(outer) C*-algebra automorphism of $\mathcal{T}^{\otimes n}$ permuting its
tensor factors and preserving its ideal $\mathcal{K}^{\otimes n}$, the above
expression of free generators $\left[  \partial_{n}\left(  \otimes^{j}%
I\otimes\otimes^{n-j}P_{1}\right)  \right]  $ with $0<j\leq n$ of
$K_{0}\left(  C\left(  \mathbb{P}^{n-1}\left(  \mathcal{T}\right)  \right)
\right)  $ can be changed by a permutation to yield some other free
generators. For example, both $\left\{  \left[  \partial_{3}\left(  1\otimes
P_{1}\otimes P_{1}\right)  \right]  ,\left[  \partial_{3}\left(
1\otimes1\otimes P_{1}\right)  \right]  ,\left[  \partial_{3}\left(
1\otimes1\otimes1\right)  \right]  \right\}  $ and\newline$\quad\left\{
\left[  \partial_{3}\left(  P_{1}\otimes P_{1}\otimes1\right)  \right]
,\left[  \partial_{3}\left(  P_{1}\otimes1\otimes1\right)  \right]  ,\left[
\partial_{3}\left(  1\otimes1\otimes1\right)  \right]  \right\}  $ are sets of
free generators of $K_{0}\left(  C\left(  \mathbb{P}^{2}\left(  \mathcal{T}%
\right)  \right)  \right)  $.

The above theorem shows that for $\left(  j^{\prime},\operatorname{id}\right)
\prec\left(  j,\operatorname{id}\right)  $ in $\Omega_{0}$, i.e. $0<j^{\prime
}<j$, it is not true that $Q_{j,\operatorname{id},1}\boxplus Q_{j^{\prime
},\operatorname{id},1}\sim Q_{j,\operatorname{id},1}$ over $C\left(
\mathbb{P}^{n-1}\left(  \mathcal{T}\right)  \right)  $ because
\[
\left[  Q_{j,\operatorname{id},1}\boxplus Q_{j^{\prime},\operatorname{id}%
,1}\right]  =\left[  Q_{j,\operatorname{id},1}\right]  +\left[  Q_{j^{\prime
},\operatorname{id},1}\right]  \neq\left[  Q_{j,\operatorname{id},1}\right]
\]
in $K_{0}\left(  C\left(  \mathbb{P}^{n-1}\left(  \mathcal{T}\right)  \right)
\right)  $.

Next we consider the positive cone of $K_{0}\left(  C\left(  \mathbb{P}%
^{n-1}\left(  \mathcal{T}\right)  \right)  \right)  $. In the following, we
use $\hat{I}$ and $\tilde{I}$ to denote the identity elements of
$\mathcal{T}^{\otimes n-1}$ and $\mathcal{T}^{\otimes n}$ respectively.

First, it is easy to see that for $k>0$, the projection $\hat{I}\otimes P_{k}
$ is a sum of $k$ mutually orthogonal projections $\hat{I}\otimes e_{jj}$,
each equivalent to $\hat{I}\otimes P_{1}$ over $\left(  \mathcal{T}^{\otimes
n-1}\otimes\mathcal{K}\right)  ^{+}\subset\mathcal{T}^{\otimes n}$, and hence
the projection $\partial_{n}\left(  \hat{I}\otimes P_{k}\right)  $ is a sum of
$k$ mutually orthogonal projections $\partial_{n}\left(  \hat{I}\otimes
e_{jj}\right)  $, each equivalent to $\partial_{n}\left(  \hat{I}\otimes
P_{1}\right)  $ over $C\left(  \mathbb{S}_{H}^{2n-1}\right)  $. So%
\[
\hat{I}\otimes P_{k}\sim\boxplus^{k}\left(  \hat{I}\otimes P_{1}\right)
\equiv\boxplus^{k}P_{n-1,1}\equiv P_{n-1,k}\text{ \ over }\left(
\mathcal{T}^{\otimes n-1}\otimes\mathcal{K}\right)  ^{+}\subset\mathcal{T}%
^{\otimes n}%
\]
and $\partial_{n}\left(  \hat{I}\otimes P_{k}\right)  \sim\boxplus
^{k}Q_{n-1,\operatorname{id},1}\equiv Q_{n-1,\operatorname{id},k}$ over
$C\left(  \mathbb{S}_{H}^{2n-1}\right)  $. Similarly, by rearranging entries
via conjugation by shifts, the projection $\hat{I}\otimes P_{-k}$ is
equivalent to $\tilde{I}$ over $\mathcal{T}^{\otimes n}$, and hence
$\partial_{n}\left(  \hat{I}\otimes P_{-k}\right)  \sim\partial_{n}\left(
\tilde{I}\right)  $ over $C\left(  \mathbb{S}_{H}^{2n-1}\right)  $. However
such equivalences no longer hold over the algebra $C\left(  \mathbb{P}%
^{n-1}\left(  \mathcal{T}\right)  \right)  \subset C\left(  \mathbb{S}%
_{H}^{2n-1}\right)  $. For example, $\partial_{n}\left(  \hat{I}\otimes
P_{-k}\right)  \boxplus\partial_{n}\left(  \hat{I}\otimes P_{k}\right)
\sim\partial_{n}\left(  \tilde{I}\right)  $ over $C\left(  \mathbb{P}%
^{n-1}\left(  \mathcal{T}\right)  \right)  $ since $\partial_{n}\left(
\hat{I}\otimes P_{-k}\right)  $ and $\partial_{n}\left(  \hat{I}\otimes
P_{k}\right)  $ are orthogonal projections in $C\left(  \mathbb{P}%
^{n-1}\left(  \mathcal{T}\right)  \right)  $ which add up to $\tilde{I}$. So
\[
\left[  \partial_{n}\left(  \hat{I}\otimes P_{-1}\right)  \right]  =\left[
\partial_{n}\left(  \tilde{I}\right)  \right]  -\left[  \partial_{n}\left(
\hat{I}\otimes P_{1}\right)  \right]  =\left[  Q_{n,\operatorname{id}%
,1}\right]  -\left[  Q_{n-1,\operatorname{id},1}\right]  \text{\ \ in }%
K_{0}\left(  C\left(  \mathbb{P}^{n-1}\left(  \mathcal{T}\right)  \right)
\right)
\]
showing that
\[
\left[  \partial_{n}\left(  \hat{I}\otimes P_{-1}\right)  \right]
\in\mathbb{Z}^{n-2}\times\left\{  -1\right\}  \times\left\{  1\right\}
\subset\mathbb{Z}^{n}\cong K_{0}\left(  C\left(  \mathbb{P}^{n-1}\left(
\mathcal{T}\right)  \right)  \right)
\]
and $\partial_{n}\left(  \hat{I}\otimes P_{-1}\right)  $ is not even stably
equivalent over $C\left(  \mathbb{P}^{n-1}\left(  \mathcal{T}\right)  \right)
$ to any $\boxplus$-sum of the $K_{0}$-generating projections
$Q_{j,\operatorname{id},1}$ with $0<j\leq n$.

From now on, we include all projections of the form $\partial_{n}\left(
\left(  \otimes^{j-1}I\right)  \otimes P_{k}\otimes\left(  \otimes^{n-j}%
P_{1}\right)  \right)  $ with $k\in\mathbb{Z}$ as elementary projections over
$C\left(  \mathbb{P}^{n-1}\left(  \mathcal{T}\right)  \right)  $, where it is
understood that for $k=0$, we take $P_{k}:=P_{-0}\equiv I$ instead of
$P_{0}\equiv0$.

\textbf{Theorem 5}. The positive cone of $K_{0}\left(  C\left(  \mathbb{P}%
^{n-1}\left(  \mathcal{T}\right)  \right)  \right)  \cong\mathbb{Z}^{n}%
\equiv\bigoplus_{j=1}^{n}\mathbb{Z}\left[  Q_{j,\operatorname{id},1}\right]  $
contains
\[
\mathbb{Z}^{n}\backslash\left\{  z\in\mathbb{Z}^{n}:z_{j}<0=z_{j+1}%
=\cdots=z_{n}\text{ for some }1\leq j\leq n\right\}
\]
which is the part of the cone generated/spanned by the equivalence classes of
the elementary projections $\partial_{n}\left(  \left(  \otimes^{j-1}I\right)
\otimes P_{k}\otimes\left(  \otimes^{n-j}P_{1}\right)  \right)  $ with
$k\in\mathbb{Z}$ and $1\leq j\leq n$, where for $k=0$, we take $P_{k}%
:=P_{-0}\equiv I$.

Proof. In \cite{Sh:pmqpl}, it has been established that in the case of $n=2$,
the positive cone of
\[
K_{0}\left(  C\left(  \mathbb{P}^{1}\left(  \mathcal{T}\right)  \right)
\right)  =\mathbb{Z}\left[  Q_{1,\operatorname{id},1}\right]  \oplus
\mathbb{Z}\left[  Q_{2,\operatorname{id},1}\right]  \equiv\mathbb{Z}\left[
\partial_{2}\left(  I\otimes P_{1}\right)  \right]  \oplus\mathbb{Z}\left[
\partial_{2}\left(  I\otimes I\right)  \right]  \cong\mathbb{Z}^{2}%
\]
consists of $\left(  k,m\right)  \in\mathbb{Z}^{2}$ with either $k\geq0$ or
the rank $m>0$, such that $\left[  \partial_{2}\left(  I\otimes P_{k}\right)
\right]  =k\left[  \partial_{2}\left(  I\otimes P_{1}\right)  \right]
=\left(  k,0\right)  $ and
\[
\left[  \partial_{2}\left(  I\otimes P_{-k}\right)  \right]  =\left[
\partial_{2}\left(  I\otimes I\right)  \right]  -k\left[  \partial_{2}\left(
I\otimes P_{1}\right)  \right]  =\left(  -k,1\right)
\]
in $K_{0}\left(  C\left(  \mathbb{P}^{1}\left(  \mathcal{T}\right)  \right)
\right)  $ for all $k>0$.

By induction on $n$, we can show that the positive cone of
\[
K_{0}\left(  C\left(  \mathbb{P}^{n-1}\left(  \mathcal{T}\right)  \right)
\right)  =\mathbb{Z}\left[  Q_{1,\operatorname{id},1}\right]  \oplus
\cdots\oplus\mathbb{Z}\left[  Q_{n,\operatorname{id},1}\right]  \cong%
\mathbb{Z}^{n}%
\]
contains the set $\left(  \mathbb{Z}_{\geq}^{n-1}\times\left\{  0\right\}
\right)  \cup\left(  \mathbb{Z}^{n-1}\times\mathbb{N}\right)  $ consisting of
$\left(  k_{1},..,k_{n-1},m\right)  \in\mathbb{Z}^{n}$ with either $k_{j}%
\geq0$ for all $j$ or the rank $m>0$.

Indeed, under the canonical embedding
\[
\iota:C\left(  \mathbb{P}^{n-2}\left(  \mathcal{T}\right)  \right)  \equiv
C^{\ast}\left(  \left(  \mathfrak{G}_{n-1}\right)  _{0}\right)  \rightarrow
C\left(  \mathbb{P}^{n-1}\left(  \mathcal{T}\right)  \right)  \equiv C^{\ast
}\left(  \left(  \mathfrak{G}_{n}\right)  _{0}\right)
\]
due to the degree-preserving groupoid embedding of $\left.  \left(
\mathbb{Z}^{n-1}\ltimes\overline{\mathbb{Z}}^{n-1}\right)  \right\vert
_{\overline{\mathbb{Z}_{\geq}}^{n-1}}$ in $\left.  \left(  \mathbb{Z}%
^{n}\ltimes\overline{\mathbb{Z}}^{n}\right)  \right\vert _{\overline
{\mathbb{Z}_{\geq}}^{n}}$ as $\left.  \left(  \left(  \mathbb{Z}^{n-1}%
\times\left\{  0\right\}  \right)  \ltimes\left(  \overline{\mathbb{Z}}%
^{n-1}\times\left\{  0\right\}  \right)  \right)  \right\vert _{\overline
{\mathbb{Z}_{\geq}}^{n-1}\times\left\{  0\right\}  }$, a projection $p$ (for
example, $\partial_{n-1}\left(  P_{k_{1}}\otimes\cdots\otimes P_{k_{n-1}%
}\right)  $) over $C\left(  \mathbb{P}^{n-2}\left(  \mathcal{T}\right)
\right)  $ becomes the projection $p\otimes P_{1}$ (for example, $\partial
_{n}\left(  P_{k_{1}}\otimes\cdots\otimes P_{k_{n-1}}\otimes P_{1}\right)  $)
over $C\left(  \mathbb{P}^{n-1}\left(  \mathcal{T}\right)  \right)  $.
Furthermore if $p\sim q$ over $C\left(  \mathbb{P}^{n-2}\left(  \mathcal{T}%
\right)  \right)  $, say, $upu^{-1}=q$ for some $u\in GL_{\infty}\left(
C\left(  \mathbb{P}^{n-2}\left(  \mathcal{T}\right)  \right)  \right)  $ then
the equivalence $p\otimes P_{1}\sim q\otimes P_{1}$ over $C\left(
\mathbb{P}^{n-1}\left(  \mathcal{T}\right)  \right)  $ can be explicitly
constructed as
\[
\left(  \left(  u\otimes P_{1}\right)  +\partial_{n}\left(  I\otimes
P_{-1}\right)  \right)  \left(  p\otimes P_{1}\right)  \left(  \left(
u\otimes P_{1}\right)  +\partial_{n}\left(  I\otimes P_{-1}\right)  \right)
^{-1}=q\otimes P_{1}%
\]
with $\left(  u\otimes P_{1}\right)  +\partial_{n}\left(  I\otimes
P_{-1}\right)  \in GL_{\infty}\left(  C\left(  \mathbb{P}^{n-1}\left(
\mathcal{T}\right)  \right)  \right)  $. Now consider the well-defined group
homomorphism%
\[
K_{0}\left(  \iota\right)  :K_{0}\left(  C\left(  \mathbb{P}^{n-2}\left(
\mathcal{T}\right)  \right)  \right)  \cong\mathbb{Z}^{n-1}\rightarrow
K_{0}\left(  C\left(  \mathbb{P}^{n-1}\left(  \mathcal{T}\right)  \right)
\right)  \cong\mathbb{Z}^{n}%
\]
mapping the positive cone of $K_{0}\left(  C\left(  \mathbb{P}^{n-2}\left(
\mathcal{T}\right)  \right)  \right)  $ into that of $K_{0}\left(  C\left(
\mathbb{P}^{n-1}\left(  \mathcal{T}\right)  \right)  \right)  $. Since under
$\iota$, the projection $Q_{j,\operatorname{id},1}$ over $C\left(
\mathbb{P}^{n-2}\left(  \mathcal{T}\right)  \right)  $ for $0<j\leq n-1$ is
sent to the projection $Q_{j,\operatorname{id},1}$ over $C\left(
\mathbb{P}^{n-1}\left(  \mathcal{T}\right)  \right)  $, by induction
hypothesis, we get that the positive cone of $K_{0}\left(  C\left(
\mathbb{P}^{n-1}\left(  \mathcal{T}\right)  \right)  \right)  \cong%
\mathbb{Z}^{n}$ contains $\left(  \mathbb{Z}_{\geq}^{n-2}\times\left\{
0\right\}  \times\left\{  0\right\}  \right)  \cup\left(  \mathbb{Z}%
^{n-2}\times\mathbb{N}\times\left\{  0\right\}  \right)  $, and hence $\left(
\mathbb{Z}_{\geq}^{n-1}\times\left\{  0\right\}  \right)  \cup\left(
\mathbb{Z}^{n-2}\times\mathbb{N}\times\mathbb{Z}_{\geq}\right)  $.

On the other hand, for $k>0$,
\[
\hat{I}\otimes P_{-k}=\left(  \hat{I}\otimes P_{-\left(  k+1\right)  }\right)
\boxplus\left(  \hat{I}\otimes e_{kk}\right)  \sim\left(  \hat{I}\otimes
P_{-\left(  k+1\right)  }\right)  \boxplus\left(  I^{\prime}\otimes
P_{-k}\otimes P_{1}\right)  \text{ over }\mathcal{T}^{\otimes n}%
\]
where $I^{\prime}$ denotes the identity element of $\mathcal{T}^{\otimes n-2}$
and $e_{ij}$ with $i,j\in\mathbb{Z}_{\geq}$ represents a matrix unit
projection, because $\hat{I}\otimes P_{-k}$ is the sum of orthogonal
projections $\left(  \hat{I}\otimes P_{-\left(  k+1\right)  }\right)  $ and
$\left(  \hat{I}\otimes e_{kk}\right)  $, and $\left(  \hat{I}\otimes
e_{kk}\right)  \boxplus0$ when conjugated by
\[
u_{k}:=\left(
\begin{array}
[c]{cc}%
I^{\prime}\otimes I\otimes P_{k} & I^{\prime}\otimes\left(  \mathcal{S}%
^{k}\right)  ^{\ast}\otimes\mathcal{S}^{k}\\
I^{\prime}\otimes\mathcal{S}^{k}\otimes\left(  \mathcal{S}^{k}\right)  ^{\ast}
& I^{\prime}\otimes P_{k}\otimes I
\end{array}
\right)  \in GL_{2}\left(  \mathcal{T}^{\otimes n}\right)
\]
becomes $0\boxplus\left(  I^{\prime}\otimes P_{-k}\otimes P_{1}\right)  $.
Since $\partial_{n}\left(  u_{k}\right)  $ of total degree $0$ is in
$M_{2}\left(  C\left(  \mathbb{P}^{n-1}\left(  \mathcal{T}\right)  \right)
\right)  $, we get
\[
\partial_{n}\left(  \hat{I}\otimes P_{-k}\right)  \sim\partial_{n}\left(
\left(  \hat{I}\otimes P_{-\left(  k+1\right)  }\right)  \right)
\boxplus\iota\left(  \partial_{n-1}\left(  I^{\prime}\otimes P_{-k}\right)
\right)  \text{ over }C\left(  \mathbb{P}^{n-1}\left(  \mathcal{T}\right)
\right)
\]
and hence
\[
\left[  \partial_{n}\left(  \hat{I}\otimes P_{-k}\right)  \right]  -\left[
\partial_{n}\left(  \left(  \hat{I}\otimes P_{-\left(  k+1\right)  }\right)
\right)  \right]  \in\mathbb{Z}^{n-2}\times\left\{  1\right\}  \times\left\{
0\right\}  \text{ in }\mathbb{Z}^{n}%
\]
because $\left[  \partial_{n-1}\left(  I^{\prime}\otimes P_{-k}\right)
\right]  \in\mathbb{Z}^{n-2}\times\left\{  1\right\}  $ for the rank-one
projection $I^{\prime}\otimes P_{-k}$ over $\mathcal{T}^{\otimes n-1}$. With
\[
\left[  \partial_{n}\left(  \hat{I}\otimes P_{-1}\right)  \right]  =\left[
\partial_{n}\left(  \tilde{I}\right)  \right]  -\left[  \partial_{n}\left(
\hat{I}\otimes P_{1}\right)  \right]  =\left(  0,...,0,-1,1\right)
\in\mathbb{Z}^{n},
\]
we get inductively
\[
\left[  \partial_{n}\left(  \hat{I}\otimes P_{-k}\right)  \right]
\in\mathbb{Z}^{n-2}\times\left\{  -k\right\}  \times\left\{  1\right\}
\subset\mathbb{Z}^{n}\cong K_{0}\left(  C\left(  \mathbb{P}^{n-1}\left(
\mathcal{T}\right)  \right)  \right)
\]
for all $k>0$. Thus the positive cone of $K_{0}\left(  C\left(  \mathbb{P}%
^{n-1}\left(  \mathcal{T}\right)  \right)  \right)  \cong\mathbb{Z}^{n}$
contains $\left(  \mathbb{Z}_{\geq}^{n-1}\times\left\{  0\right\}  \right)
\cup\left(  \mathbb{Z}^{n-1}\times\left\{  1\right\}  \right)  $ and hence
$\left(  \mathbb{Z}_{\geq}^{n-1}\times\left\{  0\right\}  \right)  \cup\left(
\mathbb{Z}^{n-1}\times\mathbb{N}\right)  $. On the other hand, the positive
cone of $K_{0}\left(  C\left(  \mathbb{P}^{n-2}\left(  \mathcal{T}\right)
\right)  \right)  $ is mapped into the positive cone of $K_{0}\left(  C\left(
\mathbb{P}^{n-1}\left(  \mathcal{T}\right)  \right)  \right)  $ by the
homomorphism $\cdot\times\left\{  0\right\}  \equiv K_{0}\left(  \iota\right)
$. So it is easy to get inductively the conclusion. $\square$

We note that for the case of $n=2$, the finitely generated projective modules
over $C\left(  \mathbb{P}^{1}\left(  \mathcal{T}\right)  \right)  $ are
completely classified with the positive cone of $K_{0}\left(  C\left(
\mathbb{P}^{1}\left(  \mathcal{T}\right)  \right)  \right)  $ explicitly
identified in \cite{Sh:pmqpl}.

\section{Quantum line bundles}

In this section, we identify the quantum line bundles $L_{k}:=C\left(
\mathbb{S}_{H}^{2n-1}\right)  _{k}$ of degree $k$ over $C\left(
\mathbb{P}^{n-1}\left(  \mathcal{T}\right)  \right)  $ with a concrete
(equivalence class of) projection described in terms of the elementary
projections defined in the previous section. We continue to use $\hat{I}$ and
$\tilde{I}$ to denote the identity elements of $\mathcal{T}^{\otimes n-1}$ and
$\mathcal{T}^{\otimes n}$ respectively, and we start to use $0^{\left(
l\right)  }$ to denote the zero of $\mathbb{Z}^{l}$.

To distinguish between ordinary function product and convolution product, we
denote the groupoid C*-algebraic (convolution) multiplication of elements in
$C_{c}\left(  \mathcal{G}\right)  \subset C^{\ast}\left(  \mathcal{G}\right)
$ by $\ast$, while omitting $\ast$ when the elements are presented as
operators or when they are multiplied together pointwise as functions on
$\mathcal{G}$. We also view $C_{c}\left(  \mathfrak{G}_{n}\right)  $ or
$C_{c}\left(  \left(  \mathfrak{G}_{n}\right)  _{k}\right)  $ (also
abbreviated as $C_{c}\left(  \mathfrak{G}_{n}\right)  _{k}$) as left
$C_{c}\left(  \mathfrak{G}_{n}\right)  _{0}$-modules with $C_{c}\left(
\mathfrak{G}_{n}\right)  $ carrying the convolution algebra structure as a
subalgebra of the groupoid C*-algebra $C^{\ast}\left(  \mathfrak{G}%
_{n}\right)  $. Similarly, for a closed subset $X$ of the unit space of
$\mathfrak{G}_{n}$, the inverse image $\mathfrak{G}_{n}\upharpoonright_{X}$ of
$X$ under the source map of $\mathfrak{G}_{n}$ or its grade-$k$ component
$\left(  \mathfrak{G}_{n}\upharpoonright_{X}\right)  _{k}$ gives rise to a
left $C_{c}\left(  \mathfrak{G}_{n}\right)  _{0}$-module $C_{c}\left(
\mathfrak{G}_{n}\upharpoonright_{X}\right)  $ or $C_{c}\left(  \mathfrak{G}%
_{n}\upharpoonright_{X}\right)  _{k}$.

We define a partial isometry in $C\left(  \mathbb{S}_{H}^{2n-1}\right)  \equiv
C^{\ast}\left(  \mathfrak{G}_{n}\right)  $ for each $k\in\mathbb{Z}$ as the
characteristic function $\chi_{B_{k}}$ of the compact open set
\[
B_{k}:=\left\{  \left(  0,k,p,q\right)  \in\mathfrak{G}_{n}\subset
\mathbb{Z}^{n-1}\times\mathbb{Z}\times\overline{\mathbb{Z}_{\geq}^{n-1}}%
\times\overline{\mathbb{Z}_{\geq}}:q+k\geq0\right\}  \subset\mathfrak{G}_{n}.
\]
It is easy to verify that $\chi_{B_{k}}\in C_{c}\left(  \mathfrak{G}%
_{n}\right)  _{k}$ and $\left(  \chi_{B_{k}}\right)  ^{\ast}\in C_{c}\left(
\mathfrak{G}_{n}\right)  _{-k}$ such that
\[
\left(  \chi_{B_{k}}\right)  ^{\ast}\ast\chi_{B_{k}}=\left\{
\begin{array}
[c]{lll}%
\chi_{\left\{  0^{\left(  n\right)  }\right\}  \times\left(  \overline
{\mathbb{Z}_{\geq}^{n}}\backslash\mathbb{Z}_{\geq}^{n}\right)  }=1_{C^{\ast
}\left(  \mathfrak{G}_{n}\right)  }\equiv1_{C^{\ast}\left(  \mathfrak{G}%
_{n}\right)  _{0}} & \text{if } & k\geq0\\
\chi_{\left\{  0^{\left(  n\right)  }\right\}  \times\left(  \left(
\overline{\mathbb{Z}_{\geq}^{n-1}}\times\overline{\mathbb{Z}_{\geq\left\vert
k\right\vert }}\right)  \backslash\mathbb{Z}_{\geq}^{n}\right)  }=\partial
_{n}\left(  \hat{I}\otimes P_{-\left\vert k\right\vert }\right)  & \text{if }
& k<0
\end{array}
\right.
\]
and
\[
\chi_{B_{k}}\ast\left(  \chi_{B_{k}}\right)  ^{\ast}=\left\{
\begin{array}
[c]{lll}%
\chi_{\left\{  0^{\left(  n\right)  }\right\}  \times\left(  \left(
\overline{\mathbb{Z}_{\geq}^{n-1}}\times\overline{\mathbb{Z}_{\geq k}}\right)
\backslash\mathbb{Z}_{\geq}^{n}\right)  }=\partial_{n}\left(  \hat{I}\otimes
P_{-k}\right)  & \text{if } & k\geq0\\
\chi_{\left\{  0^{\left(  n\right)  }\right\}  \times\left(  \overline
{\mathbb{Z}_{\geq}^{n}}\backslash\mathbb{Z}_{\geq}^{n}\right)  }=1_{C^{\ast
}\left(  \mathfrak{G}_{n}\right)  }\equiv1_{C^{\ast}\left(  \mathfrak{G}%
_{n}\right)  _{0}} & \text{if } & k<0
\end{array}
\right.  .
\]

For $k\geq0$, we have $C_{c}\left(  \mathfrak{G}_{n}\right)  _{k}\ast\left(
\chi_{B_{k}}\right)  ^{\ast}\subset C_{c}\left(  \mathfrak{G}_{n}\right)  _{0}
$ and
\[
\left(  C_{c}\left(  \mathfrak{G}_{n}\right)  _{k}\ast\left(  \chi_{B_{k}%
}\right)  ^{\ast}\right)  \ast\chi_{B_{k}}=C_{c}\left(  \mathfrak{G}%
_{n}\right)  _{k}%
\]
which implies that the convolution operator $\cdot\ast\chi_{B_{k}}$ maps
$C_{c}\left(  \mathfrak{G}_{n}\right)  _{0}$ onto $C_{c}\left(  \mathfrak{G}%
_{n}\right)  _{k}$. Since $\chi_{B_{k}}\ast\left(  \chi_{B_{k}}\right)
^{\ast}=\chi_{\left\{  0^{\left(  n\right)  }\right\}  \times\left(  \left(
\overline{\mathbb{Z}_{\geq}^{n-1}}\times\overline{\mathbb{Z}_{\geq k}}\right)
\backslash\mathbb{Z}_{\geq}^{n}\right)  }$, we get $\cdot\ast\chi_{B_{k}}$
mapping $C_{c}\left(  \mathfrak{G}_{n}\right)  _{0}\ast\chi_{\left\{
0^{\left(  n\right)  }\right\}  \times\left(  \left(  \overline{\mathbb{Z}%
_{\geq}^{n-1}}\times\overline{\mathbb{Z}_{\geq k}}\right)  \backslash
\mathbb{Z}_{\geq}^{n}\right)  }$ bijectively onto $C_{c}\left(  \mathfrak{G}%
_{n}\right)  _{k}$ with $\cdot\ast\left(  \chi_{B_{k}}\right)  ^{\ast}$ as the
inverse. Furthermore $\cdot\ast\chi_{B_{k}}$ is a left $C_{c}\left(
\mathfrak{G}_{n}\right)  _{0}$-module homomorphism. With $\chi_{B_{k}}$ being
a partial isometry, $\cdot\ast\chi_{B_{k}}$ and $\cdot\ast\left(  \chi_{B_{k}%
}\right)  ^{\ast}$ extend continuously to provide an isomorphism between the
$C^{\ast}\left(  \mathfrak{G}_{n}\right)  _{0}$-modules
\[
C^{\ast}\left(  \mathfrak{G}_{n}\right)  _{0}\chi_{\left\{  0^{\left(
n\right)  }\right\}  \times\left(  \left(  \overline{\mathbb{Z}_{\geq}^{n-1}%
}\times\overline{\mathbb{Z}_{\geq k}}\right)  \backslash\mathbb{Z}_{\geq}%
^{n}\right)  }\equiv C^{\ast}\left(  \mathfrak{G}_{n}\right)  _{0}\partial
_{n}\left(  \hat{I}\otimes P_{-k}\right)
\]
and $C^{\ast}\left(  \mathfrak{G}_{n}\right)  _{k}\equiv\overline{C_{c}\left(
\mathfrak{G}_{n}\right)  _{k}}$. So the quantum line bundle $C^{\ast}\left(
\mathfrak{G}_{n}\right)  _{k}$ is identified with the projection $\partial
_{n}\left(  \hat{I}\otimes P_{-k}\right)  $.

For $k<0$, we consider the direct sum decomposition as left $C_{c}\left(
\mathfrak{G}_{n}\right)  _{0}$-modules%
\begin{align*}
C_{c}\left(  \mathfrak{G}_{n}\right)  _{k}  &  =\left(  C_{c}\left(
\mathfrak{G}_{n}\right)  _{k}\ast\chi_{\left\{  0^{\left(  n\right)
}\right\}  \times\left(  \overline{\mathbb{Z}_{\geq}^{n-1}}\backslash
\mathbb{Z}_{\geq}^{n-1}\right)  \times\left\{  0,1,...,\left\vert k\right\vert
-1\right\}  }\right)  \oplus\left(  C_{c}\left(  \mathfrak{G}_{n}\right)
_{k}\ast\chi_{\left\{  0^{\left(  n\right)  }\right\}  \times\left(  \left(
\overline{\mathbb{Z}_{\geq}^{n-1}}\times\overline{\mathbb{Z}_{\geq\left\vert
k\right\vert }}\right)  \backslash\mathbb{Z}_{\geq}^{n}\right)  }\right) \\
&  =C_{c}\left(  \mathfrak{G}_{n}\upharpoonright_{\left(  \overline
{\mathbb{Z}_{\geq}^{n-1}}\backslash\mathbb{Z}_{\geq}^{n-1}\right)
\times\left\{  0,1,...,\left\vert k\right\vert -1\right\}  }\right)
_{k}\oplus\left(  C_{c}\left(  \mathfrak{G}_{n}\right)  _{k}\ast\chi_{\left\{
0^{\left(  n\right)  }\right\}  \times\left(  \left(  \overline{\mathbb{Z}%
_{\geq}^{n-1}}\times\overline{\mathbb{Z}_{\geq\left\vert k\right\vert }%
}\right)  \backslash\mathbb{Z}_{\geq}^{n}\right)  }\right)  .
\end{align*}
From
\[
C_{c}\left(  \mathfrak{G}_{n}\right)  _{0}\ast\chi_{B_{k}}\ast\left(
\chi_{B_{k}}\right)  ^{\ast}\equiv C_{c}\left(  \mathfrak{G}_{n}\right)
_{0}\ast1_{C^{\ast}\left(  \mathfrak{G}_{n}\right)  }=C_{c}\left(
\mathfrak{G}_{n}\right)  _{0}%
\]
and%
\[
C_{c}\left(  \mathfrak{G}_{n}\right)  _{k}\ast\left(  \chi_{B_{k}}\right)
^{\ast}\ast\chi_{B_{k}}=C_{c}\left(  \mathfrak{G}_{n}\right)  _{k}\ast
\chi_{\left\{  0^{\left(  n\right)  }\right\}  \times\left(  \left(
\overline{\mathbb{Z}_{\geq}^{n-1}}\times\overline{\mathbb{Z}_{\geq\left\vert
k\right\vert }}\right)  \backslash\mathbb{Z}_{\geq}^{n}\right)  }%
\]
we see that $\cdot\ast\chi_{B_{\left\vert k\right\vert }}$ is a left
$C_{c}\left(  \mathfrak{G}_{n}\right)  _{0}$-module isomorphism between
$C_{c}\left(  \mathfrak{G}_{n}\right)  _{0}$ and $C_{c}\left(  \mathfrak{G}%
_{n}\right)  _{k}\ast\chi_{\left\{  0^{\left(  n\right)  }\right\}
\times\left(  \left(  \overline{\mathbb{Z}_{\geq}^{n-1}}\times\overline
{\mathbb{Z}_{\geq\left\vert k\right\vert }}\right)  \backslash\mathbb{Z}%
_{\geq}^{n}\right)  }$ with $\cdot\ast\left(  \chi_{B_{k}}\right)  ^{\ast}$ as
its inverse.

On the other hand, in the $C_{c}\left(  \mathfrak{G}_{n}\right)  _{0}$-module
decomposition
\[
C_{c}\left(  \mathfrak{G}_{n}\upharpoonright_{\left(  \overline{\mathbb{Z}%
_{\geq}^{n-1}}\backslash\mathbb{Z}_{\geq}^{n-1}\right)  \times\left\{
0,1,...,\left\vert k\right\vert -1\right\}  }\right)  _{k}=\bigoplus
_{j=0}^{\left\vert k\right\vert -1}C_{c}\left(  \mathfrak{G}_{n}%
\upharpoonright_{\left(  \overline{\mathbb{Z}_{\geq}^{n-1}}\backslash
\mathbb{Z}_{\geq}^{n-1}\right)  \times\left\{  j\right\}  }\right)  _{k},
\]
each $C_{c}\left(  \mathfrak{G}_{n}\upharpoonright_{\left(  \overline
{\mathbb{Z}_{\geq}^{n-1}}\backslash\mathbb{Z}_{\geq}^{n-1}\right)
\times\left\{  j\right\}  }\right)  _{k}$ is isomorphic to the $C_{c}\left(
\mathfrak{G}_{n}\right)  _{0}$-module $C_{c}\left(  \mathfrak{G}%
_{n}\upharpoonright_{\left(  \overline{\mathbb{Z}_{\geq}^{n-1}}\backslash
\mathbb{Z}_{\geq}^{n-1}\right)  \times\left\{  0\right\}  }\right)  _{k+j}$
with $k+j<0$ via the homeomorphism
\[
\left(  m,l,p,j\right)  \in\left(  \mathfrak{G}_{n}\upharpoonright_{\left(
\overline{\mathbb{Z}_{\geq}^{n-1}}\backslash\mathbb{Z}_{\geq}^{n-1}\right)
\times\left\{  j\right\}  }\right)  _{k}\subset\mathbb{Z}^{n-1}\times
\mathbb{Z}\times\overline{\mathbb{Z}_{\geq}^{n-1}}\times\mathbb{Z}_{\geq
}\mapsto\left(  m,l+j,p,0\right)  \in\left(  \mathfrak{G}_{n}\upharpoonright
_{\left(  \overline{\mathbb{Z}_{\geq}^{n-1}}\backslash\mathbb{Z}_{\geq}%
^{n-1}\right)  \times\left\{  0\right\}  }\right)  _{k+j}%
\]
where the implicit condition $l+j\geq0$ is equivalent to $l\geq-j$. So we
focus on analyzing $C_{c}\left(  \mathfrak{G}_{n}\right)  _{0}$-modules of the
form
\[
C_{c}\left(  \left(  \mathfrak{G}_{n}\upharpoonright_{\left(  \overline
{\mathbb{Z}_{\geq}^{n-1}}\backslash\mathbb{Z}_{\geq}^{n-1}\right)
\times\left\{  0\right\}  }\right)  _{-r}\right)  =C_{c}\left(  \mathfrak{G}%
_{n}\right)  _{-r}\ast\chi_{\left\{  0^{\left(  n\right)  }\right\}
\times\left(  \left(  \overline{\mathbb{Z}_{\geq}^{n-1}}\backslash
\mathbb{Z}_{\geq}^{n-1}\right)  \times\left\{  0\right\}  \right)  }%
=C_{c}\left(  \mathfrak{G}_{n}\right)  _{-r}\partial_{n}\left(  \hat{I}\otimes
P_{1}\right)
\]
with $r\geq0$. Note that the $C^{\ast}\left(  \mathfrak{G}_{n}\right)  _{0}%
$-module
\[
\overline{C_{c}\left(  \mathfrak{G}_{n}\right)  _{0}\partial_{n}\left(
\hat{I}\otimes P_{1}\right)  }=C^{\ast}\left(  \mathfrak{G}_{n}\right)
_{0}\partial_{n}\left(  \hat{I}\otimes P_{1}\right)
\]
is identified with the projection $\partial_{n}\left(  \hat{I}\otimes
P_{1}\right)  \equiv Q_{n-1,\operatorname{id},1}$.

For $r>0$, similar to the argument used above, it can be checked that the
compact open subset
\[
B_{-r}^{\prime}:=\left\{  \left(  0,-r,0,p,q,0\right)  \in\mathfrak{G}%
_{n}\subset\mathbb{Z}^{n-2}\times\mathbb{Z}\times\mathbb{Z}\times
\overline{\mathbb{Z}_{\geq}^{n-2}}\times\overline{\mathbb{Z}_{\geq}}%
\times\overline{\mathbb{Z}_{\geq}}:q\geq r\right\}  \subset\mathfrak{G}_{n}%
\]
defines a partial isometry $\chi_{B_{-r}^{\prime}}\in C_{c}\left(
\mathfrak{G}_{n}\right)  _{-r}$ with $\left(  \chi_{B_{-r}^{\prime}}\right)
^{\ast}\in C_{c}\left(  \mathfrak{G}_{n}\right)  _{r}$ such that
\[
\left(  \chi_{B_{-r}^{\prime}}\right)  ^{\ast}\ast\chi_{B_{-r}^{\prime}}%
=\chi_{\left\{  0^{\left(  n\right)  }\right\}  \times\left(  \left(
\overline{\mathbb{Z}_{\geq}^{n-2}}\times\overline{\mathbb{Z}_{\geq r}}%
\times\left\{  0\right\}  \right)  \backslash\mathbb{Z}_{\geq}^{n}\right)
}=\partial_{n}\left(  I^{\otimes n-2}\otimes P_{-r}\otimes I\right)
\]
and
\[
\chi_{B_{-r}^{\prime}}\ast\left(  \chi_{B_{-r}^{\prime}}\right)  ^{\ast}%
=\chi_{\left\{  0^{\left(  n\right)  }\right\}  \times\left(  \overline
{\mathbb{Z}_{\geq}^{n-1}}\backslash\mathbb{Z}_{\geq}^{n-1}\right)
\times\left\{  0\right\}  }.
\]
In the decomposition
\[
C_{c}\left(  \mathfrak{G}_{n}\right)  _{-r}\ast\chi_{\left\{  0^{\left(
n\right)  }\right\}  \times\left(  \left(  \overline{\mathbb{Z}_{\geq}^{n-1}%
}\backslash\mathbb{Z}_{\geq}^{n-1}\right)  \times\left\{  0\right\}  \right)
}%
\]%
\[
=\left(  C_{c}\left(  \mathfrak{G}_{n}\right)  _{-r}\ast\left(  \chi_{\left\{
0^{\left(  n\right)  }\right\}  \times\left(  \overline{\mathbb{Z}_{\geq
}^{n-2}}\backslash\mathbb{Z}_{\geq}^{n-2}\right)  \times\left\{
0,1,...,r-1\right\}  \times\left\{  0\right\}  }\right)  \right)
\oplus\left(  C_{c}\left(  \mathfrak{G}_{n}\right)  _{-r}\ast\chi_{\left\{
0^{\left(  n\right)  }\right\}  \times\left(  \left(  \overline{\mathbb{Z}%
_{\geq}^{n-2}}\times\overline{\mathbb{Z}_{\geq r}}\times\left\{  0\right\}
\right)  \backslash\mathbb{Z}_{\geq}^{n}\right)  }\right)
\]%
\[
=C_{c}\left(  \mathfrak{G}_{n}\upharpoonright_{\left(  \overline
{\mathbb{Z}_{\geq}^{n-2}}\backslash\mathbb{Z}_{\geq}^{n-2}\right)
\times\left\{  0,1,...,r-1\right\}  \times\left\{  0\right\}  }\right)
_{-r}\oplus\left(  C_{c}\left(  \mathfrak{G}_{n}\right)  _{-r}\ast
\chi_{\left\{  0^{\left(  n\right)  }\right\}  \times\left(  \left(
\overline{\mathbb{Z}_{\geq}^{n-2}}\times\overline{\mathbb{Z}_{\geq r}}%
\times\left\{  0\right\}  \right)  \backslash\mathbb{Z}_{\geq}^{n}\right)
}\right)  ,
\]
the second summand is isomorphic, via the right convolution $\cdot\ast\left(
\chi_{B_{-r}^{\prime}}\right)  ^{\ast}$ by the partial isometry $\left(
\chi_{B_{-r}^{\prime}}\right)  ^{\ast}$, to the $C_{c}\left(  \mathfrak{G}%
_{n}\right)  _{0}$-module%
\[
C_{c}\left(  \mathfrak{G}_{n}\right)  _{0}\ast\chi_{\left\{  0^{\left(
n\right)  }\right\}  \times\left(  \left(  \overline{\mathbb{Z}_{\geq}^{n-1}%
}\times\left\{  0\right\}  \right)  \backslash\mathbb{Z}_{\geq}^{n}\right)
}=C_{c}\left(  \mathfrak{G}_{n}\right)  _{0}\partial_{n}\left(  \hat{I}\otimes
P_{1}\right)  .
\]

Now we introduce the notation of a $C_{c}\left(  \mathfrak{G}_{n}\right)
_{0}$-module
\[
A_{r,l}:=C_{c}\left(  \left(  \mathfrak{G}_{n}\upharpoonright_{\left(
\overline{\mathbb{Z}_{\geq}^{l}}\backslash\mathbb{Z}_{\geq}^{l}\right)
\times\left\{  0^{\left(  n-l\right)  }\right\}  }\right)  _{-r}\right)
\subset C_{c}\left(  \left(  \mathfrak{G}_{n}\right)  _{-r}\right)  \subset
C_{c}\left(  \mathfrak{G}_{n}\right)
\]
for $r\geq0$ and $1\leq l\leq n-1$. We note that the $C_{c}\left(
\mathfrak{G}_{n}\right)  _{0}$-module
\[
A_{r,1}=C_{c}\left(  \left(  \mathfrak{G}_{n}\upharpoonright_{\left\{
\infty\right\}  \times\left\{  0^{\left(  n-1\right)  }\right\}  }\right)
_{-r}\right)
\]
is isomorphic to
\[
C_{c}\left(  \left(  \mathfrak{G}_{n}\upharpoonright_{\left\{  \infty\right\}
\times\left\{  0^{\left(  n-1\right)  }\right\}  }\right)  _{0}\right)
=C_{c}\left(  \left(  \mathfrak{G}_{n}\upharpoonright_{\left(  \overline
{\mathbb{Z}_{\geq}}\times\left\{  0^{\left(  n-1\right)  }\right\}  \right)
\backslash\mathbb{Z}_{\geq}^{n}}\right)  _{0}\right)
\]%
\[
\equiv C_{c}\left(  \mathfrak{G}_{n}\right)  _{0}\ast\chi_{\left\{  0^{\left(
n\right)  }\right\}  \times\left(  \left(  \overline{\mathbb{Z}_{\geq}}%
\times\left\{  0^{\left(  n-1\right)  }\right\}  \right)  \backslash
\mathbb{Z}_{\geq}^{n}\right)  }=C_{c}\left(  \mathfrak{G}_{n}\right)
_{0}\partial_{n}\left(  I\otimes P_{1}^{\otimes n-1}\right)
\]
via the homeomorphism%
\begin{align*}
\left(  s,t,\infty,0^{\left(  n-1\right)  }\right)   &  \in\left(
\mathfrak{G}_{n}\upharpoonright_{\left\{  \infty\right\}  \times\left\{
0^{\left(  n-1\right)  }\right\}  }\right)  _{-r}\subset\mathbb{Z}%
\times\mathbb{Z}^{n-1}\times\left\{  \infty\right\}  \times\overline
{\mathbb{Z}_{\geq}^{n-1}}\\
&  \mapsto\left(  s+r,t,\infty,0^{\left(  n-1\right)  }\right)  \in\left(
\mathfrak{G}_{n}\upharpoonright_{\left\{  \infty\right\}  \times\left\{
0^{\left(  n-1\right)  }\right\}  }\right)  _{0}.
\end{align*}

Applying the same kind of arguments as shown above, we get the isomorphism of
$C_{c}\left(  \mathfrak{G}_{n}\right)  _{0}$-modules%
\[
A_{r,l}\cong C_{c}\left(  \mathfrak{G}_{n}\upharpoonright_{\left(
\overline{\mathbb{Z}_{\geq}^{l-1}}\backslash\mathbb{Z}_{\geq}^{l-1}\right)
\times\left\{  0,1,...,r-1\right\}  \times\left\{  0^{\left(  n-l\right)
}\right\}  }\right)  _{-r}\oplus\left(  C_{c}\left(  \mathfrak{G}_{n}\right)
_{0}\ast\chi_{\left\{  0^{\left(  n\right)  }\right\}  \times\left(  \left(
\overline{\mathbb{Z}_{\geq}^{l}}\backslash\mathbb{Z}_{\geq}^{l}\right)
\times\left\{  0^{\left(  n-l\right)  }\right\}  \right)  }\right)
\]%
\[
\cong\bigoplus_{j=0}^{r-1}C_{c}\left(  \mathfrak{G}_{n}\upharpoonright
_{\left(  \overline{\mathbb{Z}_{\geq}^{l-1}}\backslash\mathbb{Z}_{\geq}%
^{l-1}\right)  \times\left\{  j\right\}  \times\left\{  0^{\left(  n-l\right)
}\right\}  }\right)  _{-r}\oplus\left(  C_{c}\left(  \mathfrak{G}_{n}\right)
_{0}\partial_{n}\left(  I^{\otimes l}\otimes P_{1}^{\otimes n-l}\right)
\right)
\]%
\[
\cong\bigoplus_{j=0}^{r-1}C_{c}\left(  \mathfrak{G}_{n}\upharpoonright
_{\left(  \overline{\mathbb{Z}_{\geq}^{l-1}}\backslash\mathbb{Z}_{\geq}%
^{l-1}\right)  \times\left\{  0\right\}  \times\left\{  0^{\left(  n-l\right)
}\right\}  }\right)  _{-r+j}\oplus\left(  C_{c}\left(  \mathfrak{G}%
_{n}\right)  _{0}\partial_{n}\left(  I^{\otimes l}\otimes P_{1}^{\otimes
n-l}\right)  \right)
\]%
\[
=\bigoplus_{j=0}^{r-1}A_{r-j,l-1}\oplus\left(  C_{c}\left(  \mathfrak{G}%
_{n}\right)  _{0}\partial_{n}\left(  I^{\otimes l}\otimes P_{1}^{\otimes
n-l}\right)  \right)
\]
for $2\leq l\leq n-1$. This provides a recursive formula to reduce the index
$l$ of the module $A_{r,l}$.

For $n>2$, we define a combinatorial number $\nu_{n}\left(  m,l\right)  $
recursively by
\[
\nu_{n}\left(  m,l\right)  :=\sum_{s=0}^{m}\nu_{n}\left(  s,l-1\right)
\]
and $\nu_{n}\left(  m,1\right)  :=1$, for $m\geq0$ and $2\leq l\leq n-1$, to
be used in the following theorem.

Thanks to Thomas Timmermann, as he pointed out to the author, $\nu_{n}\left(
m,l\right)  $ can be identified with a familiar combinatorial number, namely,%
\[
\nu_{n}\left(  m,l\right)  =C_{m}^{m+l-1}%
\]
for all $m\geq0$ and $l\geq1$. Indeed, if either $l=1$ or $m=0$ (e.g. when
$m+l\leq2$), we get easily from the definition, $\nu_{n}\left(  m,l\right)
=1=C_{m}^{m+l-1}$. On the other hand, for $l\geq2$ and $m\geq1$, since
\[
\nu_{n}\left(  m,l\right)  =\sum_{s=0}^{m-1}\nu_{n}\left(  s,l-1\right)
+\nu_{n}\left(  m,l-1\right)  =\nu_{n}\left(  m-1,l\right)  +\nu_{n}\left(
m,l-1\right)  ,
\]
the identification can be proved by an induction on $m+l\geq3$ as shown in
\[
\nu_{n}\left(  m-1,l\right)  +\nu_{n}\left(  m,l-1\right)  =C_{m-1}%
^{m+l-2}+C_{m}^{m+l-2}=C_{m}^{m+l-1}%
\]
which is valid due to either the induction hypothesis for $m+l-1$ (in the case
of $m+l-1>2$) or the already established identification (for the case of
$m+l-1=2$).

\textbf{Theorem 6}. For $n>2$, the quantum line bundle $L_{k}\equiv C\left(
\mathbb{S}_{H}^{2n-1}\right)  _{k}$ of degree $k\in\mathbb{Z}$ over $C\left(
\mathbb{P}^{n-1}\left(  \mathcal{T}\right)  \right)  $ is isomorphic to the
finitely generated projective left module over $C\left(  \mathbb{P}%
^{n-1}\left(  \mathcal{T}\right)  \right)  $ determined by the projection
$\partial_{n}\left(  \otimes^{n-1}I\otimes P_{-k}\right)  $ if $k\geq0$, and
the projection%
\[
\left(  \boxplus^{\sum_{m=0}^{\left\vert k\right\vert -1}\left(  \left\vert
k\right\vert -m\right)  \nu_{n}\left(  m,n-2\right)  }\partial_{n}\left(
I\otimes P_{1}^{\otimes n-1}\right)  \right)  \boxplus\left(  \boxplus
_{l=1}^{n-1}\boxplus^{\nu_{n}\left(  \left\vert k\right\vert -1,l\right)
}\partial_{n}\left(  I^{\otimes n-l+1}\otimes P_{1}^{\otimes l-1}\right)
\right)
\]
if $k<0$.

Proof. Only the case of $k<0$ remains to be proved as follows.

For $k<0$, starting with the established isomorphism%
\[
C_{c}\left(  \mathfrak{G}_{n}\right)  _{k}=C_{c}\left(  \mathfrak{G}%
_{n}\upharpoonright_{\left(  \overline{\mathbb{Z}_{\geq}^{n-1}}\backslash
\mathbb{Z}_{\geq}^{n-1}\right)  \times\left\{  0,1,...,\left\vert k\right\vert
-1\right\}  }\right)  _{k}\oplus C_{c}\left(  \mathfrak{G}_{n}\right)  _{0}%
\]%
\[
=\bigoplus_{m=0}^{\left\vert k\right\vert -1}C_{c}\left(  \mathfrak{G}%
_{n}\upharpoonright_{\left(  \overline{\mathbb{Z}_{\geq}^{n-1}}\backslash
\mathbb{Z}_{\geq}^{n-1}\right)  \times\left\{  m\right\}  }\right)  _{k}\oplus
C_{c}\left(  \mathfrak{G}_{n}\right)  _{0}\cong\bigoplus_{m=0}^{\left\vert
k\right\vert -1}A_{\left\vert k\right\vert -m,n-1}\oplus C_{c}\left(
\mathfrak{G}_{n}\right)  _{0},
\]
we apply repeatedly the recursive formula
\[
A_{r,l}=\bigoplus_{j=0}^{r-1}A_{r-j,l-1}\oplus\left(  C_{c}\left(
\mathfrak{G}_{n}\right)  _{0}\partial_{n}\left(  I^{\otimes l}\otimes
P_{1}^{\otimes n-l}\right)  \right)
\]
reducing $l$ for $A_{r,l}$ with $2\leq l\leq n$ until $l$ reaches $2$ with
\[
A_{r,2}\cong\left(  \bigoplus_{j=0}^{r-1}\left(  C_{c}\left(  \mathfrak{G}%
_{n}\right)  _{0}\partial_{n}\left(  I\otimes P_{1}^{\otimes n-1}\right)
\right)  \right)  \oplus\left(  C_{c}\left(  \mathfrak{G}_{n}\right)
_{0}\partial_{n}\left(  I^{\otimes2}\otimes P_{1}^{\otimes n-2}\right)
\right)  ,
\]
in order to convert all terms to $C_{c}\left(  \mathfrak{G}_{n}\right)  _{0}%
$-modules of the form $C_{c}\left(  \mathfrak{G}_{n}\right)  _{0}\partial
_{n}\left(  I^{\otimes j}\otimes P_{1}^{\otimes n-j}\right)  $ for some
$0<j\leq n$.

In fact, we check inductively on $1\leq j\leq n-2$ that
\[
\text{(*)\ \ }C_{c}\left(  \mathfrak{G}_{n}\right)  _{k}\cong\bigoplus
_{m=0}^{\left\vert k\right\vert -1}\left(  \oplus^{\nu_{n}\left(  m,j\right)
}A_{\left\vert k\right\vert -m,n-j}\right)  \oplus\bigoplus_{l=1}^{j}\left(
\oplus^{\nu_{n}\left(  \left\vert k\right\vert -1,l\right)  }C_{c}\left(
\mathfrak{G}_{n}\right)  _{0}\partial_{n}\left(  I^{\otimes n-l+1}\otimes
P_{1}^{\otimes l-1}\right)  \right)  .
\]
The case of $j=1$ is our starting point already proved. Now assuming that it
holds for $j$, we get by the above recursive formula
\begin{align*}
C_{c}\left(  \mathfrak{G}_{n}\right)  _{k}  &  \cong\bigoplus_{m=0}%
^{\left\vert k\right\vert -1}\oplus^{\nu_{n}\left(  m,j\right)  }\left(
\left(  \bigoplus_{s=0}^{\left\vert k\right\vert -m-1}A_{\left\vert
k\right\vert -m-s,n-j-1}\right)  \oplus\left(  C_{c}\left(  \mathfrak{G}%
_{n}\right)  _{0}\partial_{n}\left(  I^{\otimes n-j}\otimes P_{1}^{\otimes
j}\right)  \right)  \right) \\
&  \oplus\left(  \bigoplus_{l=1}^{j}\left(  \oplus^{\nu_{n}\left(  \left\vert
k\right\vert -1,l\right)  }C_{c}\left(  \mathfrak{G}_{n}\right)  _{0}%
\partial_{n}\left(  I^{\otimes n-l+1}\otimes P_{1}^{\otimes l-1}\right)
\right)  \right)
\end{align*}%
\begin{align*}
&  \cong\left(  \bigoplus_{m=0}^{\left\vert k\right\vert -1}\oplus^{\nu
_{n}\left(  m,j\right)  }\bigoplus_{s=0}^{\left\vert k\right\vert
-m-1}A_{\left\vert k\right\vert -m-s,n-j-1}\right)  \oplus\left(  \oplus
^{\sum_{m=0}^{\left\vert k\right\vert -1}\nu_{n}\left(  m,j\right)  }%
C_{c}\left(  \mathfrak{G}_{n}\right)  _{0}\partial_{n}\left(  I^{\otimes
n-j}\otimes P_{1}^{\otimes j}\right)  \right) \\
&  \oplus\left(  \bigoplus_{l=1}^{j}\left(  \oplus^{\nu_{n}\left(  \left\vert
k\right\vert -1,l\right)  }C_{c}\left(  \mathfrak{G}_{n}\right)  _{0}%
\partial_{n}\left(  I^{\otimes n-l+1}\otimes P_{1}^{\otimes l-1}\right)
\right)  \right)
\end{align*}%
\[
\cong\left(  \bigoplus_{m=0}^{\left\vert k\right\vert -1}\oplus^{\nu
_{n}\left(  m,j\right)  }\bigoplus_{s=0}^{\left\vert k\right\vert
-m-1}A_{\left\vert k\right\vert -m-s,n-j-1}\right)  \oplus\bigoplus
_{l=1}^{j+1}\left(  \oplus^{\nu_{n}\left(  \left\vert k\right\vert
-1,l\right)  }C_{c}\left(  \mathfrak{G}_{n}\right)  _{0}\partial_{n}\left(
I^{\otimes n-l+1}\otimes P_{1}^{\otimes l-1}\right)  \right)
\]%
\[
=\left(  \bigoplus_{m=0}^{\left\vert k\right\vert -1}\bigoplus_{s=0}^{m}%
\oplus^{\nu_{n}\left(  s,j\right)  }A_{\left\vert k\right\vert -m,n-j-1}%
\right)  \oplus\bigoplus_{l=1}^{j+1}\left(  \oplus^{\nu_{n}\left(  \left\vert
k\right\vert -1,l\right)  }C_{c}\left(  \mathfrak{G}_{n}\right)  _{0}%
\partial_{n}\left(  I^{\otimes n-l+1}\otimes P_{1}^{\otimes l-1}\right)
\right)
\]%
\[
=\left(  \bigoplus_{m=0}^{\left\vert k\right\vert -1}\oplus^{\nu_{n}\left(
m,j+1\right)  }A_{\left\vert k\right\vert -m,n-j-1}\right)  \oplus
\bigoplus_{l=1}^{j+1}\left(  \oplus^{\nu_{n}\left(  \left\vert k\right\vert
-1,l\right)  }C_{c}\left(  \mathfrak{G}_{n}\right)  _{0}\partial_{n}\left(
I^{\otimes n-l+1}\otimes P_{1}^{\otimes l-1}\right)  \right)  ,
\]
which verifies (*) for $j+1$.

For $j=n-2$, (*) says
\[
C_{c}\left(  \mathfrak{G}_{n}\right)  _{k}\cong\bigoplus_{m=0}^{\left\vert
k\right\vert -1}\left(  \oplus^{\nu_{n}\left(  m,n-2\right)  }A_{\left\vert
k\right\vert -m,2}\right)  \oplus\bigoplus_{l=1}^{n-2}\left(  \oplus^{\nu
_{n}\left(  \left\vert k\right\vert -1,l\right)  }C_{c}\left(  \mathfrak{G}%
_{n}\right)  _{0}\partial_{n}\left(  I^{\otimes n-l+1}\otimes P_{1}^{\otimes
l-1}\right)  \right)
\]%
\begin{align*}
&  \cong\bigoplus_{m=0}^{\left\vert k\right\vert -1}\oplus^{\nu_{n}\left(
m,n-2\right)  }\left(  \left(  \bigoplus_{j=0}^{\left\vert k\right\vert
-m-1}\left(  C_{c}\left(  \mathfrak{G}_{n}\right)  _{0}\partial_{n}\left(
I\otimes P_{1}^{\otimes n-1}\right)  \right)  \right)  \oplus\left(
C_{c}\left(  \mathfrak{G}_{n}\right)  _{0}\partial_{n}\left(  I^{\otimes
2}\otimes P_{1}^{\otimes n-2}\right)  \right)  \right) \\
&  \oplus\bigoplus_{l=1}^{n-2}\left(  \oplus^{\nu_{n}\left(  \left\vert
k\right\vert -1,l\right)  }C_{c}\left(  \mathfrak{G}_{n}\right)  _{0}%
\partial_{n}\left(  I^{\otimes n-l+1}\otimes P_{1}^{\otimes l-1}\right)
\right)
\end{align*}%
\begin{align*}
&  \cong\left(  \oplus^{\sum_{m=0}^{\left\vert k\right\vert -1}\left(
\left\vert k\right\vert -m\right)  \nu_{n}\left(  m,n-2\right)  }C_{c}\left(
\mathfrak{G}_{n}\right)  _{0}\partial_{n}\left(  I\otimes P_{1}^{\otimes
n-1}\right)  \right) \\
&  \oplus\left(  \oplus^{\sum_{m=0}^{\left\vert k\right\vert -1}\nu_{n}\left(
m,n-2\right)  }\left(  C_{c}\left(  \mathfrak{G}_{n}\right)  _{0}\partial
_{n}\left(  I^{\otimes2}\otimes P_{1}^{\otimes n-2}\right)  \right)  \right)
\\
&  \oplus\bigoplus_{l=1}^{n-2}\left(  \oplus^{\nu_{n}\left(  \left\vert
k\right\vert -1,l\right)  }C_{c}\left(  \mathfrak{G}_{n}\right)  _{0}%
\partial_{n}\left(  I^{\otimes n-l+1}\otimes P_{1}^{\otimes l-1}\right)
\right)
\end{align*}%
\begin{align*}
&  \cong\left(  \oplus^{\sum_{m=0}^{\left\vert k\right\vert -1}\left(
\left\vert k\right\vert -m\right)  \nu_{n}\left(  m,n-2\right)  }C_{c}\left(
\mathfrak{G}_{n}\right)  _{0}\partial_{n}\left(  I\otimes P_{1}^{\otimes
n-1}\right)  \right) \\
&  \oplus\bigoplus_{l=1}^{n-1}\left(  \oplus^{\nu_{n}\left(  \left\vert
k\right\vert -1,l\right)  }C_{c}\left(  \mathfrak{G}_{n}\right)  _{0}%
\partial_{n}\left(  I^{\otimes n-l+1}\otimes P_{1}^{\otimes l-1}\right)
\right)  .
\end{align*}

After completion, we get the $C^{\ast}\left(  \left(  \mathfrak{G}_{n}\right)
_{0}\right)  $-module $L_{k}$ isomorphic to
\begin{align*}
&  \left(  \oplus^{\sum_{m=0}^{\left\vert k\right\vert -1}\left(  \left\vert
k\right\vert -m\right)  \nu_{n}\left(  m,n-2\right)  }C^{\ast}\left(  \left(
\mathfrak{G}_{n}\right)  _{0}\right)  \partial_{n}\left(  I\otimes
P_{1}^{\otimes n-1}\right)  \right) \\
&  \oplus\bigoplus_{l=1}^{n-1}\left(  \oplus^{\nu_{n}\left(  \left\vert
k\right\vert -1,l\right)  }C^{\ast}\left(  \left(  \mathfrak{G}_{n}\right)
_{0}\right)  \partial_{n}\left(  I^{\otimes n-l+1}\otimes P_{1}^{\otimes
l-1}\right)  \right)
\end{align*}
which corresponds to the projection%
\[
\left(  \boxplus^{\sum_{m=0}^{\left\vert k\right\vert -1}\left(  \left\vert
k\right\vert -m\right)  \nu_{n}\left(  m,n-2\right)  }\partial_{n}\left(
I\otimes P_{1}^{\otimes n-1}\right)  \right)  \boxplus\left(  \boxplus
_{l=1}^{n-1}\left(  \boxplus^{\nu_{n}\left(  \left\vert k\right\vert
-1,l\right)  }\partial_{n}\left(  I^{\otimes n-l+1}\otimes P_{1}^{\otimes
l-1}\right)  \right)  \right)  .
\]
$\square$

Little is known about the cancellation problem and hence the classification
problem for finitely generated projective modules over $C\left(
\mathbb{P}^{n-1}\left(  \mathcal{T}\right)  \right)  $. We expect that these
problems will be far more complicated than those for over $C\left(
\mathbb{S}_{H}^{2n-1}\right)  $ and $C\left(  \mathcal{T}^{\otimes n}\right)
$.

The recent work of Farsi, Hajac, Maszczyk, and Zieli\'{n}ski \cite{FHMZ}
identifies one of three free generators of $K_{0}\left(  C\left(
\mathbb{P}^{2}\left(  \mathcal{T}\right)  \right)  \right)  $ as $\left[
L_{1}\right]  +\left[  L_{-1}\right]  -2\left[  I\right]  $ (in addition to
$\left[  L_{1}\right]  -\left[  I\right]  $ and $\left[  I\right]  $)
constructed from a Milnor module and then expresses it in terms of elementary
projections, showing a perfect consistency with our result.

\end{document}